\documentclass[a4paper,onecolumn,11pt]{article}
\usepackage{geometry}
\usepackage{amsmath}
\usepackage{indentfirst}
\usepackage{ntheorem}
\usepackage{authblk}
\usepackage{bbm}
\usepackage{graphicx}
\usepackage{multirow}
\numberwithin{equation}{section}
\newtheorem{thm}{Theorem}[section]
\newtheorem{prop}[thm]{Proposition}
\newtheorem{lem}[thm]{Lemma}
\newtheorem{cor}[thm]{Corollary}

\newcommand{\RNum}[1]{\uppercase\expandafter{\romannumeral #1\relax}}
\usepackage{bm}
\usepackage{amssymb,amsmath,mathrsfs,txfonts}
\newtheorem{remark}{\indent Remark}[section]
\usepackage[colorlinks=true]{hyperref}
\hypersetup{ linkcolor = black}
\hypersetup{ citecolor= black}

\title{\bf Asymptotic stability of planar rarefaction wave to a 2D hyperbolic-elliptic coupled system of the radiating gas on half space}
\author[1]{Minyi Zhang}
\author[2]{Changjiang Zhu\thanks{Corresponding author. Email: machjzhu@scut.edu.cn}}
\affil[1,2]{School of Mathematics, South China University of Technology, Guangzhou 510641, P.R. China}
\date{} 

\begin{document}
\maketitle

\begin{abstract}
  \indent This paper studies the asymptotic stability of solution to an initial-boundary value problem for a hyperbolic-elliptic coupled system on two-dimensional half space, where the data on the boundary and at the far field are prescribed as $u_-$ and $u_+$. We show that the solution to the problem converges to the corresponding planar rarefaction wave for $0\le u_-<u_+$ as time tends to infinity under smallness assumptions on the initial perturbation and wave strength. These results are based on the analysis of div-curl decomposition, the standard $L^2$-energy method, $L^1$-estimate, and the monotonicity of profile is given by the maximum principle.\\

\vspace{3mm}

\end{abstract}

\vspace{3mm}


\tableofcontents

\section{Introdution}

In this paper, we consider the asymptotic stability of solution to an initial-boundary value problem for a hyperbolic-elliptic coupled system on two-dimensional half space:
\begin{equation}\label{eq-yfc}
   \begin{cases}
     u_t+f(u)_x+g(u)_y+{\rm{div}}{q}=0,\quad (x,y)\in \mathbbm{R}_+\times\mathbbm{R}, \quad t>0,\\
     -\nabla{\rm{div}}{q}+q+\nabla{u}=0,\quad (x,y)\in \mathbbm{R}_+\times\mathbbm{R}, \quad t>0,
   \end{cases}
\end{equation}
with initial data
\begin{equation}\label{eq-initial}
  u(x,y,0)=u_0(x,y),
\end{equation}
and boundary condition
\begin{align}
  &u(0,y,t)=u_-,\quad t>0,\label{eq-bjtj}\\
  &u(+\infty,y,t)=u_+,\quad q(+\infty,y,t)=0,\quad t>0,\label{eq-bjtjinf}
\end{align}
where $u$ and $q=(q_1,q_2)$ are dependent variables with values in $\mathbbm{R}$ and $\mathbbm{R}^2$, respectively. $f,g$ are smooth functions and $u_\pm$ are constants. We assume that $f$ is strictly convex, i.e., for a certain positive constant $\alpha$,
\begin{equation}\label{eq-f-tj1}
  f''(u)\ge \alpha >0,\qquad u\in\mathbbm{R},
\end{equation}
and that the characteristic speeds $f'(u_\pm)$ satisfy
\begin{equation}\label{eq-f-tj2}
  f'(u_-)< f'(u_+).
\end{equation}
Without loss of generality, we also assume
\begin{equation}\label{eq-f-tj3}
  f(0)=f'(0)=0.
\end{equation}

The study of \eqref{eq-yfc} is motivated by physical models or the so-called radiative gas model, which describes the dynamics of a gas in which radiation is present. The model consists of the compressible Euler equations coupled with an elliptic system representing the radiative flux and is given by Vincenti and Kruger in \cite{Vincenti1965}. The system \eqref{eq-yfc}, first mentioned by Hamer in \cite{Hamer1971}, simplifies the model for the motion of radiating gas on two-dimensional half space. For the deduction of system $\eqref{eq-yfc}$, we refer to \cite{Francesco, Gao2008, Hamer1971, Vincenti1965}. In the \emph{in-flow} case of $u_->0$, the boundary condition $\eqref{eq-bjtj}$ is necessary for the single hyperbolic equation $\eqref{eq-yfc}_1$. In addition, we also require
\begin{equation}\label{eq-divp0tj}
  {\rm{div}}q(0,y,t)=0,
\end{equation}
 for the solvability of the coupled elliptic equation \eqref{eq-yfc}$_2$. On the contrary, in the \emph{out-flow} case of $u_-<0$, the boundary condition \eqref{eq-bjtj} is enough and necessary for the solvability of the hyperbolic-elliptic coupled system \eqref{eq-yfc}.

The conditions \eqref{eq-f-tj1} and \eqref{eq-f-tj2} give $u_-<u_+$. Referring to \cite{Liu1998}, the asymptotic behavior of the solution to the one-dimensional scalar viscous conversation law is classified into three cases according to the signs of the characteristic speeds $f'(u_\pm)$:
\begin{equation*}
	\begin{aligned}
		&(a) \ \ f'(u_-)<f'(u_+)\le 0,\\
		&(b) \ \ 0\le f'(u_-)<f'(u_+),\\
		&(c) \ \ f'(u_-)<0<f'(u_+).\\
	\end{aligned}
\end{equation*}
In case $(a)$, the solution of $\eqref{eq-yfc}$ converges to the corresponding stationary solution $(\bar{u},\bar{q})$ which satisfies
\begin{equation}\label{eq-wtj}
	\begin{cases}
		f(\bar{u})_x+\bar{q}_x=0, &x\in \mathbbm{R}_+,\\
		-\bar{q}_{xx}+\bar{q}+\bar{u}_x=0,  &x\in \mathbbm{R}_+,\\
        \bar{u}(0,t)=u_-, &\lim\limits_{x\rightarrow +\infty} \bar{u}(x)=u_+.
	\end{cases}
\end{equation}
In case $(b)$, the solution behaviors as the rarefaction wave which satisfies
\begin{equation}
  \begin{cases}
    r_t+f(r)_x=0,\quad x\in\mathbbm{R},\quad t>0,\\
    r(x,0)=r_0(x):=
    \begin{cases}
      u_-,\quad x<0,\\
      u_+,\quad x>0.
    \end{cases}
  \end{cases}
\end{equation}
Here $r(x,t)=r(\frac{x}{t})$ is given explicitly by
\begin{equation}\label{eq-zxxsb}
  r(x,t)=
  \begin{cases}
    u_-,& x\le f'(u_-)t,\\
    (f')^{-1}(\frac{x}{t}), &f'(u_-)t\le x\le f'(u_+)t,\\
    u_+,& x\ge f'(u_+)t.
  \end{cases}
\end{equation}
In case $(c)$, the solution tends to the superposition of stationary solution and rarefaction wave.

There are many works concerning the asymptotic stability of solution for different physical systems. Most of them are in the case of one dimension. The pioneer work for stability of nonlinear wave to Cauchy problem for scalar viscous conversation law was done by Il'in and Oleinik in \cite{Il'in1960} in 1960. Then, the convergence rate toward the rarefaction wave was first obtained by Harabetian in \cite{Harabetian1988} for viscous Burgers equation and the further work has been investigated by many authors in \cite{Hattori1991, Matsumura1986, Matsumura1992}.  For the half space problem for scalar viscous conversation law, Liu and Nishihara in \cite{Liu1997} considered the asymptotic stability of a viscous shock wave for the case of $u_->u_+$. For the case where $u_-<u_+$, Liu, Matsumura and Nishihara in \cite{Liu1998} first proved the asymptotic stability of rarefaction wave or stationary solution or superposition of these two kind of waves. The convergence rate of the rarefaction wave on half space was gained by Nakamura in \cite{Nakamura2003}. By a combination of the weighted $L^p$ energy method and the $L^1$ estimate, Hashimoto, Ueda and Kawashima in \cite{Hashimoto2009} obtained the convergence rate of the  superposition of stationary solution and rarefaction wave. These problems was also considered for hyperbolic-elliptic coupled system \eqref{eq-yfc} on one-dimensional space. Tanaka in \cite{Tanaka1995} proved that the solution approaches the diffusion wave for the case where $u_-=u_+=0$. Kawashima and Nishibata in \cite{Kawashima1998} investigated the existence and the asymptotic stability of traveling wave and obtained the convergence rate for $u_->u_+$. The remaining case $u_-<u_+$ was studied by Kawashima and Tanaka in \cite{Kawashima2004}. They showed the asymptotic stability of the rarefaction wave and obtained the convergence rate. Recently, the initial boundary value problem of the one-dimensional system \eqref{eq-yfc} was also studied thoroughly. Ruan and Zhu in \cite{Ruan2008} considered the case $0=u_-<u_+$ for a hyperbolic-elliptic system on one-dimensional half space toward the corresponding rarefaction wave. Moreover, Ji, Zhang and Zhu in \cite{Ji2021} further showed the asymptotic stability of the stationary solution, rarefaction wave, and the superposition of these two kind of waves for the case $u_-<u_+\le 0$, $0\le u_-<u_+$ and $u_-<0<u_+$, respectively. In addition, the asymptotic stability of rarefaction waves or stationary solutions for the compressible Navier-Stokes equations was fully studied in \cite{Kagei2005, Kawashima2003NS, Kawashima2009, Nakamura2007}.

In the cases of multi-dimensional space, Xin in \cite{Xin1990} first investigated the asymptotic stability of planar rarefaction wave for viscous conversation laws in two dimensions, and then Ito in \cite{Ito1996} gained the decay rate. Nishikawa in \cite{Nishikawa2000} further improved their results without smallness conditions.
Kawashima, Nishibata and Nishikawa in \cite{Kawashima2003} first studied the asymptotic stability of planar stationary solution for viscous conversation laws on two-dimensional half space and obtained the convergence rate, and Ueda, Nakamura and Kawashima in \cite{Ueda2010} optimized the decay rate for the degenerate case by utilizing the time weighted $L^p$ energy method.
Recently, there are also many papers discuss the asymptotic stability for hyperbolic-elliptic coupled system in multi-dimensional whole space. Gao, Ruan and Zhu investigated the asymptotic decay rate to the Cauchy problem of the planar rarefaction wave for a hyperbolic-elliptic system in $\mathbbm{R}^n$ $(n=2,3,4,5)$, cf. \cite{Gao, Gao2008}. For the initial-boundary problem for system \eqref{eq-yfc}-\eqref{eq-bjtjinf}, the case $(a)$ corresponding to planar stationary solution was considered by Zhang and Zhu in \cite{Zhang2021}. However, for the case $(b)$ where $0\le f'(u_-)<f'(u_+)$, the stability of planar rarefaction wave has been left open.
For the initial-boundary value problem of other physically meaningful equations in multi-dimensional case, there are interesting results about the stability of planar stationary solution, we refer to \cite{Kagei2006, Nakamura2009} for the compressible Navier-Stokes equations. However, the corresponding stability results of planar rarefaction wave is few. Because it is quite difficult to show the monotonicity profile and decay rate of the corresponding one-dimensional equations for all equations. On the other hand, the large-time behavior of the solution to the compressible Navier Stokes equations is determined by the Riemann problem to the corresponding inviscid Euler system, which contains a planar rarefaction wave in the genuinely nonlinear characteristic fields. This will lead to the emergence of error terms composed of planar rarefaction wave in the perturbation equations, which are only related to $x$ and independent of $y$, and the integration of these terms over $y\in\mathbbm{R}$ is divergent. As far as we know, the stability results of planar rarefaction wave for the multi-dimensional Navier-Stokes equations only consider the case in an infinitely long flat nozzle domain, cf. \cite{Li2018, Li18ARMA, Wang2021}. Obviously, in this case ($y$ belongs to a bounded domain), the planar rarefaction wave is integrable in $y$-direction. Therefore, it is meaningful to study the stability of the planar rarefaction wave to an initial-boundary value problem for the system \eqref{eq-yfc} on two-dimensional half space.

Our aim of this paper is to prove the large-time behaviors of the solution of \eqref{eq-yfc}-\eqref{eq-divp0tj} for the case $(b)$ where $0\le f'(u_-)<f'(u_+)$. Compared with the one-dimensional stability results in \cite{Ji2021} and the Cauchy problem in \cite{ Gao2008}, the additional difficulties here lie in the boundary estimates of higher order derivatives, the integral of the planar rarefaction wave over $y$-direction and the $L^1$-estimate of the perturbation $V:=U- \tilde{u}$. To overcome the first difficulty, we give the relation between the boundary values in Lemma \ref{lem-bj-gj} and Lemma \ref{lem-bjv-gj}, and utilize the tangential derivative estimates in $y$-direction and in $t$-direction to control the boundary terms occurred in normal direction. In addition, by utilizing the analysis of div-curl decomposition, we convert the boundary term in the form of $p(0,y,t)$ to the form of ${\rm div}p(0,y,t)$ to complete the $H^3$-estimates of the perturbation $p:=q-Q$. For the second difficulty, this is a hard problem faced by all multi-dimensional equations. To solve the second difficulty, we consider the one-dimensional initial-boundary value problem \eqref{eq-UQ-fc}-\eqref{eq-UQ-boundary} to further approximate to the rarefaction wave. Fortunately, the one-dimensional system \eqref{eq-UQ-fc}-\eqref{eq-UQ-boundary} can be rewritten as the following scalar equation form with convolution
\begin{equation} \label{eq-scalar}
U_t+f(U)_x+U-KU=0, \quad x\in \mathbbm{R}_+,\quad t>0,
\end{equation}
where $K$ is the inverse of the elliptic operator $-\partial_x^2 +1$ in $\mathbbm{R}_+$, which is defined in \eqref{eq-def-K}. It is this scalar equation form that allows us to generalize the monotonicity result in \cite{Xin1990} to our problem \eqref{eq-UQ-fc}-\eqref{eq-UQ-boundary}. From the maximum principle, we prove that the solution $U$ of system \eqref{eq-UQ-fc}-\eqref{eq-UQ-boundary} satisfies $U_x(x,t)\ge 0$ for any $(x,t)\in\mathbbm{R}_+\times\mathbbm{R}_+$ provided the initial data $U_0(x)$ is a smooth non-decreasing function of $x$-direction (see details in Lemma \ref{lem-U-dz}). Moreover, we obtain the decay rate of $U$ and its derivatives which plays an important role in estimating the perturbation $v:=u- U$. Thanks to the structure of \eqref{eq-scalar}, we can convert the perturbation equations \eqref{eq-rd-VP-1} into the form of \eqref{eq-VKV} which helps us deal with the third difficulty (see in Lemma \ref{lem-VL1-xsb}).

This paper is organized as follows. In Section \ref{sec-2}, we summarize some basic properties of the planar rarefaction wave, which are given in \cite{Hattori1991,Kawashima2004}. Then we reformulate the initial-boundary problem \eqref{eq-yfc}-\eqref{eq-bjtj} and present our main theorem. In Section \ref{sec-3}, we show the asymptotic behavior for the case of $0\le f'(u_-)<f'(u_+)$, which corresponds to the planar rarefaction wave. More precisely, we show that if the rarefaction wave strength is suitably weak ($|u_+-u_-|\ll 1$) and initial data $u_0$ in \eqref{eq-initial} is suitably close to the planar rarefaction wave, then the initial-boundary value problem \eqref{eq-yfc}-\eqref{eq-bjtj} has a global-in-time solution which converges to the planar rarefaction wave $r(x,t)$ as time tends to infinity. The specific results can be found in Theorem \ref{thm-main}.

\textbf{Notations.} Throughout this paper, we denote generic constants by $C$ and $c$ unless they need to be distinguished. We denote $\mathbbm{R}_+\times\mathbbm{R}$ by $\mathbbm{R}_+^2$. Let $\Omega=\mathbbm{R}_+$ or $\mathbbm{R}_+^2$. For any nonnegative constant $p$ $(1\le p\le \infty)$, $L^p(\Omega)$ denotes usual Lebesgue space over $\Omega$, equipped with the norm $\|\cdot\|_{L^p(\Omega)}$. For any $l\ge0$, $H^l(\Omega)$ denotes the usual Sobolev space over $\Omega$ with the norm $\|\cdot\|_{H^l(\Omega)}$. We use the notation $\nabla^k f$ as in the meaning
$$\nabla^k f=(\partial_x^k f,\partial_x^{(k-1)}\partial_y f,\cdots,\partial_x\partial_y^{(k-1)} f,\partial_y^k f),$$
where $f=f(x,y,t)$ and $\nabla^0 f=f$. And the notation $\Delta:=\partial_x^2+\partial_y^2$ denotes the Laplacian.

\section{Preliminaries and main theorem}\label{sec-2}
\subsection{Preliminaries}
In order to prove the main Theorem \ref{thm-main}, we first give some inequalities which will be used later. The following Lemma \ref{lem-parxf} is given in \cite{Iguchi2002}.
\begin{lem}\label{lem-parxf} Assume that $N\ge 2$ is a integer, $l_1,l_2,\cdots,l_N$ are non-negative integers, $1\le p,q,r\le \infty$ with $1/p=1/q+1/r$. Let $l=l_1+l_2+\cdots+l_N$. Then there exists a positive constant $C=C(N,p,q,r,l)$ such that the inequality
\begin{equation}\label{eq-f-parx}
  \left\| \prod\limits_{j=1}^N (\partial_x^{l_j} f_j) \right\|_{L^p(\mathbbm{R}_+)}\le C \| f \|_{L^\infty(\mathbbm{R}_+)}^{N-2}\| f \|_{L^q(\mathbbm{R}_+)}\| \partial_x^l f \|_{L^r(\mathbbm{R}_+)},
\end{equation}
holds for any $f(x)=(f_1(x),f_2(x),\cdots,f_N(x))$.

As an application of inequality $\eqref{eq-f-parx}$, for any $1\le p\le \infty $, \cite{Gao2008} points out that
\begin{equation}\label{eq-fw-lp}
  \left\| \partial_x^k\left\{\frac{f'''(w)}{f''(w)}w_x^2\right\} \right\|_{L^p(\mathbbm{R}_+)}\le C \| \partial_x^{k+2} w \|_{L^p(\mathbbm{R}_+)},
\end{equation}
provided $\| w \|_{L^\infty(\mathbbm{R}_+)}$ is bounded.
\end{lem}

Throughout this paper, we denote the rarefaction wave strength $\delta$ as
$$\delta:=|u_--u_+|.$$

Since the rarefaction wave $r(x,t)$ is only Lipschitz continuous, we need to find a smooth approximation rarefaction wave through the viscous Burgers equation as in \cite{Hattori1991} and \cite{Kawashima2004}. Define $\tilde{w}(x,t)$ as a solution of the Cauchy problem
\begin{equation}\label{eq-rare-tituw}
  \begin{cases}
    \tilde{w}_t+\tilde{w}\tilde{w}_x=\tilde{w}_{xx},\quad x\in\mathbbm{R},\quad t>0,\\
    \tilde{w}(x,0)=w_0^R(x),\quad x\in\mathbbm{R}.
  \end{cases}
\end{equation}
For the case of $f'(u_-)>0$, the initial data define by
\begin{equation}\label{eq-w-initial0}
w_0^R(x):=
  \begin{cases}
    f'(u_-),& x<0,\\
    f'(u_+),& x>0.
  \end{cases}
\end{equation}
For the case of $f'(u_-)=0$, $\tilde{w}(x,t)$ does not converge to the corresponding rarefaction wave fast enough around the boundary $x=0$ under the initial condition $\eqref{eq-w-initial0}$. Thus, when $f'(u_-)=0$, $w_0^R(x)$ is given as
\begin{equation}
  w_0^R(x):=
  \begin{cases}
    -f'(u_+),& x<0,\\
    f'(u_+),& x>0,
  \end{cases}
\end{equation}
which yields $\tilde{w}(0,t)=0$. Using the Hopf-Cole transformation, we can obtain the explicit expression of $\tilde{w}(x,t)$. From $\eqref{eq-f-tj1}$, we define a smooth approximation $w(x,t)$ of the rarefaction wave as
\begin{equation}\label{eq-rare-w}
  w(x,t)=(f')^{-1}(\tilde{w}(x,t)).
\end{equation}
Substituting $\eqref{eq-rare-w}$ into $\eqref{eq-rare-tituw}$, $w(x,t)$ satisfy the equation
\begin{equation}\label{eq-w-fc}
  \begin{cases}
    w_t+f(w)_x=w_{xx}+\frac{f'''(w)}{f''(w)}w_x^2,\quad x\in \mathbbm{R},\quad t>0,\\[3mm]
    w(x,0)=w_0(x):=(f')^{-1}(\tilde{w}(x,0)),\quad x\in \mathbbm{R}.
  \end{cases}
\end{equation}
The monotonicity and decay rate of rarefaction wave have been thoroughly studied in \cite{Hattori1991} and \cite{Kawashima2004}. We directly give the properties of the smooth rarefaction wave $w_i(x,t)$, $(i=1,2)$ in Lemma \ref{lem-ghxsb}, where $w_1(x,t)$ corresponds to the case $f'(u_-)=0$ and $w_2(x,t)$ corresponds to the case $f'(u_-)>0$.
\begin{lem}\label{lem-ghxsb}
For $1\le p\le \infty$ and $t>0$, the smooth rarefaction wave $w_{i}(x,t) ~ (i=1,2)$ satisfies the following properties:
\\[3mm] \indent
$\mathrm{(i)}$ $0\le w_2(0,t)-u_-\le C \delta \mathrm{e}^{-c(1+t)}$ for $f'(u_-)>0$ and $w_1(0,t)=0$ for $f'(u_-)=0$;
\\[3mm] \indent
$\mathrm{(ii)}$ $|\partial_x^k\partial_t^l w_2(0,t)|\le C \delta \mathrm{e}^{-c(1+t)}$, \quad $|\partial_x^k\partial_t^l w_1(0,t)|\le C (1+t)^{-\frac{1}{2} (k+l+1)}$, \quad $k+ l=1,2,3,4, \quad k,l\in \mathbbm{N}$;
\\[3mm] \indent
$\mathrm{(iii)}$ $\|w_i(t)-r(t)\|_{L^p(\mathbbm{R})}\le C (1+t)^{-\frac{1}{2}+\frac{1}{2p}  }$;
\\[3mm] \indent
$\mathrm{(iv)}$ $\|w_{ix}(t)\|_{L^p(\mathbbm{R})}\le C \delta^{\frac{1}{p} }(1+t)^{-1+\frac{1}{p} }$, \quad $||w_{it}(t)||_{L^p}\le C \delta^{\frac{1}{p} }(1+t)^{-1+\frac{1}{p} }$;
\\[3mm] \indent
$\mathrm{(v)}$ $\|\partial_x^k \partial_t^l w_i(t)\|_{L^p(\mathbbm{R})}\le C \delta(1+t)^{-\frac{1}{2}(k+l -\frac{1}{p}) }$, \quad $k+ l=1,2,3,4, \quad k,l\in \mathbbm{N}$;
\\[3mm] \indent
$\mathrm{(vi)}$ $\|\partial_x^k \partial_t^l w_i(t)\|_{L^p(\mathbbm{R})}\le C (1+t)^{-\frac{1}{2}(k+l+1 -\frac{1}{p}) }$, \quad $k+ l=1,2,3,4, \quad k,l\in \mathbbm{N}$;
\\[3mm] \indent
$\mathrm{(vii)}$ $w_{ix}>0, \quad x \in \mathbbm{R}$.
\end{lem}
For the case of $f'(u_-)>0$, we know that for any $x\in(-\infty,+\infty)$, $w_2(x,t)$ satisfies $u_-< w_2(x,t)<u_+$ owning to Lemma \ref{lem-ghxsb} $\mathrm{(vii)}$, which yields $w_2(0,t)\not=u_-$. Therefore, we need to modify $w_2(x,t)$ around the boundary ${x=0}$. For simplicity, we still denote $w_2$ as $w$ in the following. By employing the idea of Nakamura in \cite{Nakamura2003}, we define the modified smooth approximation $(\tilde{u},\tilde{q})$ as
\begin{equation}\label{eq-tiu-tiq}
  \begin{cases}
    \tilde{u}(x,t)=w(x,t)-\hat{u}(x,t),\\
    \tilde{q}(x,t)=-w_x(x,t)-\hat{q}(x,t),
  \end{cases}
\end{equation}
where
\begin{equation}\label{eq-hatu-def}
  \begin{cases}
   \hat{u}(x,t)=(w(0,t)-u_-){\rm{e}}^{-x},\\
    \hat{q}(x,t)=w_{xx}(0,t){\rm{e}}^{-x}.
  \end{cases}
\end{equation}
Note that $\hat{u}(0,t)\equiv 0$, if $f'(u_-)=0$. Substituting $\eqref{eq-tiu-tiq}$ into $\eqref{eq-w-fc}$ and capturing $x\in \mathbbm{R}_+$, we have
\begin{equation}\label{eq-jzxsb-fc}
  \begin{cases}
    \tilde{u}_t+f(\tilde{u})_x=\tilde{u}_{xx}+\hat{u}_{xx}-\hat{u}_t-(f(\tilde{u}+\hat{u})-f(\tilde{u}))_x+\frac{f'''(w)}{f''(w)}w_x^2,\quad x\in \mathbbm{R}_+,\quad t>0,\\
    \tilde{q}=-\tilde{u}_x- \hat{u}_x- \hat{q},\quad x\in \mathbbm{R}_+,\quad t>0,\\
    \tilde{u}(0,t)=u_-,\quad \tilde{q}_x(0,t)=0,\quad t>0,\\
    \tilde{u}(x,0)=\tilde{u}_0(x):=w_0(x)-\hat{u}(x,0),\quad x\in \mathbbm{R}_+.
  \end{cases}
\end{equation}
By simple calculation and implying Lemma \ref{lem-ghxsb}, the estimates of $\tilde{u}(x,t)$ can be obtained. For details, we refer to \cite{Hattori1991,Kawashima2004,Nakamura2003}.
\begin{lem}\label{lem-jzxsb} Suppose that $f'(u_-)>0$. For $1\le p\le \infty$ and $t>0$, the smooth rarefaction wave $\tilde{u}(x,t)$ satisfies the following:
\\[3mm] \indent
$\mathrm{(i)}$ $|\partial_x^k\partial_t^l \tilde{u}(0,t)|\le C \delta \mathrm{e}^{-c(1+t)}$,\quad $k+ l=1,2,3,4, \quad k,l\in \mathbbm{N}$ ;
\\[3mm] \indent
$\mathrm{(ii)}$ $\|\tilde{u}(t)-r(t)\|_{L^p(\mathbbm{R}_+)}\le C(1+t)^{-\frac{1}{2} +\frac{1}{2p}}$;
\\[3mm] \indent
$\mathrm{(iii)}$ $\|\tilde{u}_{x}(t)\|_{L^p(\mathbbm{R}_+)}\le C \delta^{\frac{1}{p} }(1+t)^{-1+\frac{1}{p} }$, \quad $||\tilde{u}_{t}(t)||_{L^p}\le C \delta^{\frac{1}{p} }(1+t)^{-1+\frac{1}{p} }$;
\\[3mm] \indent
$\mathrm{(iv)}$ $\|\partial_x^k \partial_t^l \tilde{u}(t)\|_{L^p(\mathbbm{R}_+)}\le C \delta(1+t)^{-\frac{1}{2}(k+l -\frac{1}{p}) }$, \quad $k+ l=1,2,3,4, \quad k,l\in \mathbbm{N}$;
\\[3mm] \indent
$\mathrm{(v)}$ $\|\partial_x^k \partial_t^l \tilde{u}(t)\|_{L^p(\mathbbm{R}_+)}\le C (1+t)^{-\frac{1}{2}(k+l+1 -\frac{1}{p}) }$,\quad $k+ l=2,3,4, \quad k,l\in \mathbbm{N}$;
\\[3mm] \indent
$\mathrm{(vi)}$ $\tilde{u}_{x}>0,\quad x \in \mathbbm{R}$.
\end{lem}
Since the terms on the right side of $\eqref{eq-jzxsb-fc}_1$ are not integrable with respect to $y$, we consider the one-dimensional initial-boundary value problem corresponding to \eqref{eq-yfc}-\eqref{eq-divp0tj}, which can further approximate to $\tilde{u}$.
\begin{equation}\label{eq-UQ-fc}
  \begin{cases}
    U_t+f(U)_x+Q_x=0,\\
    -Q_{xx}+Q+U_x=0,
  \end{cases}
\end{equation}
with initial data
\begin{equation}\label{eq-U-ini}
  U(x,0)=U_0(x),
\end{equation}
and boundary condition
\begin{equation}\label{eq-UQ-boundary}
  U(0,t)=u_-,\quad Q_x(0,t)=0,\quad U(+ \infty,t)=u_+,\quad Q(+ \infty,t)=0.
\end{equation}

Referring to \cite{Gao2008}, we get the monotonicity of $U(x,t)$ in $x$-direction by assuming that $\frac{\mathrm{d}}{\mathrm{d}x}U_0(x)\geq 0$ for $x\in\mathbbm{R}_+$.

\begin{lem}\label{lem-U-dz}(Monotonicity of profile) Assume that $U_0(x)$ is monotonically non-decreasing, i.e.
\begin{equation}
	\frac{\mathrm{d}}{\mathrm{d}x}U_0(x)\ge 0,\quad x\in \mathbbm{R}_+.
\end{equation}
Then the solution $(U(x,t),Q(x,t))$ of \eqref{eq-UQ-fc}-\eqref{eq-UQ-boundary} satisfies
\begin{equation}\label{eq-UQ-monotonicity}
	\frac{\partial}{\partial x}U(x,t)\ge0,\quad Q(x,t)\le0, \qquad (x,t)\in\mathbbm{R}_+\times\mathbbm{R}_+.
\end{equation}
\end{lem}
 {\it\bfseries Proof.}
Differential \eqref{eq-UQ-fc}$_1$ with respect to $x$ and denote $U_x(x,t)$ by $W(x,t)$. Consequently, we get
\begin{equation}\label{eq-WQ}
	\begin{cases}
		W_t+f'(U)W_x+f''(U)W^2+Q_{xx}=0,\\
		-Q_{xx}+Q+W=0,
	\end{cases}
\end{equation}
and
\begin{equation}\label{eq-U-ass}
	W(x,0)=U_x(x,0)=\frac{\mathrm{d}}{\mathrm{d}x}U_0(x)\ge 0.
\end{equation}
We extend the function $W(x,t)$ such taht
\begin{equation}
\widetilde{W}(x,t)=
\begin{cases}
	W(x,t), &x\ge 0,\\
	W(-x,t), &x< 0.
\end{cases}
\end{equation}
Then $Q(x,t)$ in \eqref{eq-WQ} can be solve as
 \begin{equation*}
 	Q(x,t)=\widetilde{Q}(x,t)|_{x\ge 0},
 \end{equation*}
 and satisfies $Q(x,0)\le 0$ owning to \eqref{eq-U-ass}, where
\begin{equation*}
	\widetilde{Q}(x,t)=-\frac{1}{2} \int_{\mathbbm{R}} \mathrm{e}^{-|x-y|}\widetilde{W}(y,t) \,\mathrm{d}y
	= -\frac{1}{2} \int_{\mathbbm{R}_+} (\mathrm{e}^{-|x-y|}+\mathrm{e}^{-|x+y|})W(y,t) \,\mathrm{d}y.
\end{equation*}
Therefore, we can rewrite \eqref{eq-WQ} as
 \begin{equation}\label{eq-WQ-fc}
 \begin{cases}
 	W_t+f'(U)W_x+f''(U)W^2+Q+W=0, & (x,t)\in\mathbbm{R}_+\times\mathbbm{R}_+,\\ 
 	Q=-\frac{1}{2} \int_{\mathbbm{R}_+} (\mathrm{e}^{-|x-y|}+\mathrm{e}^{-|x+y|})W(y,t) \,\mathrm{d}y, & (x,t)\in\mathbbm{R}_+\times\mathbbm{R}_+.
 \end{cases}
 \end{equation}
 We make the transformations
 \begin{equation}\label{eq-WbarW}
 		W=\overline{W}-\frac{1}{L}{\rm e}^{ t},
 \end{equation}
 where $L$ is a positive constant. From \eqref{eq-WQ-fc}$_2$, we get
 \begin{equation}\label{eq-def-barQ}
 	Q(x,t)=-\frac{1}{2} \int_{\mathbbm{R}_+} (\mathrm{e}^{-|x-y|}+\mathrm{e}^{-|x+y|})\overline{W}(y,t) \,\mathrm{d}y+\frac{1}{L}{\rm e}^{ t}:=\overline{Q}+\frac{1}{L}{\rm e}^{ t}.
 \end{equation}
 Consequently, we get from \eqref{eq-U-ass} and \eqref{eq-WQ-fc} that
 \begin{equation}\label{eq-WQ-bhfc}
 	\overline{W}_t+f'(U)\overline{W}_x+\left(f''(U)\overline{W}-\frac{2}{L}{\rm e}^{ t}f''(U)+1\right)\overline{W}+\overline{Q}=\frac{1}{L}{\rm e}^{ t}\left(1-f''(U)\frac{1}{L}e^{ t}\right),
 \end{equation}
 and
 \begin{equation*}
 	\overline{W}(x,0)>0,\quad \overline{Q}(x,0)<0.
 \end{equation*}
 We claim that $\overline{W}(x,t)>0$ and $\overline{Q}(x,t)<0$ for any $(x,t)\in\mathbbm{R}_+\times\mathbbm{R}_+$. For any $T>0$, let
 $$t_0=\inf\limits_{t}\{t~ |~\overline{W}(x,t)=0 \text{ or } \overline{Q}(x,t)=0, \quad \forall x\in\mathbbm{R}_+, \ \ t\in(0,T]  \}.$$
 If $t_0$ does not exist, the proof is completed. Otherwise, $t_0\in(0,T]$ and there exists $x_0\in\mathbbm{R}_+$ such that $\overline{Q}(x_0,t_0)=0$ and $\overline{W}(x,t_0)\ge0$ for any $x\in\mathbbm{R}_+$, or $\overline{W}(x_0,t_0)=0$ and $\overline{Q}(x,t_0)\le0$ for any $x\in\mathbbm{R}_+$. We argue by way of contradiction for the above two cases.

 \begin{enumerate}
 	\item[Case 1.] $\overline{Q}(x_0,t_0)=0$ and $\overline{W}(x,t_0)\ge0$ for any $x\in\mathbbm{R}_+$. From \eqref{eq-def-barQ}, we can deduce that $\overline{W}(x,t_0)=0$ for any $x\in\mathbbm{R}_+$. It follows that $W(x,t_0)$ and $Q(x,t_0)$ attain their minimum and maximum respectively at the point $x=x_0$. Notice that $\overline{W}_t(x_0,t_0)\le0$, $\overline{W}_x(x_0,t_0)=0$. If we choose $L$ sufficiently large satisfying
 	$$1-f''(U)\frac{1}{L}{\rm e}^{ t}\ge 1-f''(U)\frac{1}{L}{\rm e}^{ T}>0 \qquad \text{on } \mathbbm{R}_+\times(0,T],$$
 	from \eqref{eq-WQ-bhfc} we get a contradiction at the point $(x_0,t_0)$.
 	\item[Case 2.] $\overline{W}(x_0,t_0)=0$ and $\overline{Q}(x_0,t_0)<0$. In this case, we see $\overline{W}_t(x_0,t_0)\le0$, $\overline{W}_x(x_0,t_0)=0$. Similarly, we can get a contradiction by \eqref{eq-WQ-bhfc}, if we choose $L$ sufficiently large.
 \end{enumerate}
 Therefore, we obtain
 $$\overline{W}(x,t)>0,\quad \overline{Q}(x,t)<0, \qquad  (x,t)\in \mathbbm{R}_+\times\mathbbm{R}_+.$$
 Letting $L\rightarrow+\infty$ in \eqref{eq-WbarW}, we have
 $$W(x,t)\ge0,\quad Q(x,t)\le0, \qquad  (x,t)\in \mathbbm{R}_+\times\mathbbm{R}_+.$$
 We complete the proof of Lemma \ref{lem-U-dz}.
 $ \hfill\Box $

\vspace{3mm}

Setting
\begin{equation*}
\begin{aligned}
  u(x,y,t)-\tilde{u}(x,t)&=\{U(x,t)-\tilde{u}(x,t)\}+\{u(x,y,t)-U(x,t)\}\\
                         &:= V(x,t)+v(x,y,t),
\end{aligned}
\end{equation*}
and
\begin{equation*}
  \begin{aligned}
    q(x,y,t)-
    \left(
\begin{array}{c}
 \tilde{q}(x,t) \\
 0
 \end{array}
\right)
&=\left\{
    \left(
\begin{array}{c}
 Q(x,t) \\
 0
 \end{array}
\right)
-
    \left(
\begin{array}{c}
 \tilde{q}(x,t) \\
 0
 \end{array}
\right)
\right\}
+
\left\{
q(x,y,t)-
    \left(
\begin{array}{c}
 Q(x,t) \\
 0
 \end{array}
\right)
\right\}\\
&:=
    \left(
\begin{array}{c}
 P(x,t) \\
 0
 \end{array}
\right)
+p(x,y,t),
  \end{aligned}
\end{equation*}
we get two reformulated problems:
\begin{equation}\label{eq-rd-VP-1}
  \begin{cases}
    V_t+(f(V+\tilde{u})-f(\tilde{u}))_x+P_x=R_1,\\[1mm]
    -P_{xx}+P+V_x=R_2,\\[1mm]
    V(0,t)=0,\quad P_x(0,t)=0,\\[1mm]
    V(x,0)=V_0(x)=U_0(x)-\tilde{u}_0(x),
  \end{cases}
\end{equation}
and
\begin{equation}\label{eq-rdfc-vp-1}
  \begin{cases}
    v_t+(f(v+U)-f(U))_x+g(v+U)_y+{\rm{div}}{p}=0,\\[1mm]
    -\nabla{\rm{div}}{p}+p+\nabla{v}=0,\\[1mm]
    v(0,y,t)=0,\quad {\rm{div}}{p}(0,y,t)=0,\\[1mm]
    v(x,y,0)=v_0(x,y)=u_0(x,y)-U_0(x),
  \end{cases}
\end{equation}
where
\begin{equation}\label{eq-def-R1R2}
  \begin{cases}
    &R_1:=\hat{q}_x+\hat{u}_t+(f(\tilde{u}+\hat{u})-f(\tilde{u}))_x-\frac{f'''(w)}{f''(w)}w_x^2,\\
    &R_2:=\hat{u}_x+\hat{q}-\tilde{u}_{xxx}-\hat{u}_{xxx}-\hat{q}_{xx}.
  \end{cases}
\end{equation}
By utilizing Lemma \ref{lem-ghxsb}-\ref{lem-jzxsb}, it is not difficult to get the following Corollary \ref{cor-R1R2}.
\begin{cor}\label{cor-R1R2} Suppose that $f'(u_-)>0$. For $1\le p\le \infty$ and $t>0$, $R_1(x,t)$ and $R_2(x,t)$ satisfy
  \\[3mm] \indent
$\mathrm{(i)}$ $\| R_1(t)\|_{L^p(\mathbbm{R}_+)}\le C \delta\mathrm{e}^{-c(1+t)}+C \delta^{\frac{1}{p} }(1+t)^{-2+\frac{1}{p}}$;
  \\[3mm] \indent
$\mathrm{(ii)}$ $\| \partial_x^k\partial_t^l R_{1}(t)\|_{L^p(\mathbbm{R}_+)}\le C \delta\mathrm{e}^{-c(1+t)}+C \min \left\{\delta(1+t)^{-\frac{k+l+2}{2} +\frac{1}{2p}},(1+t)^{-\frac{k+l+3}{2} +\frac{1}{2p}} \right\}$,\quad $k+l=1,2,3,4$;
  \\[3mm] \indent
$\mathrm{(iii)}$ $\| \partial_x^k\partial_t^l R_2(t)\|_{L^p(\mathbbm{R}_+)}\le C \delta\mathrm{e}^{-c(1+t)}+C \min \left\{\delta(1+t)^{-\frac{k+l+3}{2} +\frac{1}{2p}}, (1+t)^ {-\frac{k+l+4}{2}+\frac{1}{2p}} \right\}$,\quad $k+l=0,1,2,3,4$;
  \\[3mm] \indent
$\mathrm{(iv)}$ $|\partial_x^k\partial_t^l R_{1}(0,t)|\le C \delta\mathrm{e}^{-c(1+t)}$,\quad $k,l=0,1,2,3,4$.\\[3mm]
For the case of $f'(u_-)=0$, $R_1(x,t)=-\frac{f'''(w_1)}{f''(w_1)}w_{1x}^2$ and $R_2(x,t)=-\tilde{u}_{xxx}$ satisfying
  \\[3mm] \indent
$\mathrm{(v)}$ $\| R_1(t)\|_{L^p(\mathbbm{R}_+)}\le C \delta^{\frac{1}{p} }(1+t)^{-2+\frac{1}{p}}$;
  \\[3mm] \indent
$\mathrm{(vi)}$ $\| \partial_x^k\partial_t^l R_{1}(t)\|_{L^p(\mathbbm{R}_+)}\le C \min \left\{\delta(1+t)^{-\frac{k+l+2}{2} +\frac{1}{2p}},(1+t)^{-\frac{k+l+3}{2} +\frac{1}{2p}} \right\}$,\quad $k+l=1,2,3,4$;
  \\[3mm] \indent
$\mathrm{(vii)}$ $\| \partial_x^k\partial_t^l R_2(t)\|_{L^p(\mathbbm{R}_+)}\le C \min \left\{\delta(1+t)^{-\frac{k+l+3}{2} +\frac{1}{2p}}, (1+t)^ {-\frac{k+l+4}{2}+\frac{1}{2p}} \right\}$,\quad $k+l=0,1,2,3,4$;
  \\[3mm] \indent
$\mathrm{(viii)}$ $|\partial_x^k\partial_t^l R_{1}(0,t)|\le C \min \left\{\delta(1+t)^{-\frac{k+l+2}{2} }, (1+t)^{-\frac{k+l+3}{2}} \right\}$,\quad $k+l=1,2,3,4$.
\end{cor}
 {\it\bfseries Proof.} We only give the proof of $\mathrm{(i)}$-$\mathrm{(iii)}$ and $\mathrm{(viii)}$. The remaining estimates can be obtained by similar method. From Lemma \ref{lem-ghxsb} $\mathrm{(i)} \mathrm{(iv)}$, $\eqref{eq-hatu-def}_1$ and $\eqref{eq-def-R1R2}_1$, we have
 \begin{equation*}
 \begin{aligned}
    \| R_1(t)\|_{L^p(\mathbbm{R}_+)}
   &\le C \delta\mathrm{e}^{-c(1+t)}+C\left(\int_{\mathbbm{R}_+} |w_x|^{2p} \,\mathrm{d}x\right)^{1/p}\\
   &\le  C \delta\mathrm{e}^{-c(1+t)}+C \| w_x(t) \|_{L^\infty(\mathbbm{R}_+)}\| w_x(t) \|_{L^p(\mathbbm{R}_+)} \\
   &\le  C \delta\mathrm{e}^{-c(1+t)}+C(1+t)^{-1} \delta^{ \frac{1}{p}}(1+t)^{-1+ \frac{1}{p} }\\
   &\le  C \delta\mathrm{e}^{-c(1+t)}+C\delta^{ \frac{1}{p}}(1+t)^{-2+\frac{1}{p} }.
 \end{aligned}
 \end{equation*}
 Therefore, the desired estimate $\mathrm{(i)}$ is obtained.

  Next, we try to show the estimate $\mathrm{(ii)}$. By utilizing $\eqref{eq-fw-lp}$, we can get
 \begin{equation*}
   \| \partial_x^k R_{1}(t)\|_{L^p(\mathbbm{R}_+)}
   \le C \delta\mathrm{e}^{-c(1+t)}+C\left\| \partial_x^k\left\{\frac{f'''(w)}{f''(w)}w_x^2\right\} \right\|_{L^p(\mathbbm{R}_+)}\\
   \le C \delta\mathrm{e}^{-c(1+t)}+C\|\partial_x^{k+2} w(t) \|_{L^p(\mathbbm{R}_+)}.
 \end{equation*}
 From Lemma \ref{lem-ghxsb} $\mathrm{(v)}$, it follows that
 \begin{equation}\label{eq-R1lp1}
   \begin{aligned}
         \|\partial_x^k R_1(t)\|_{L^p(\mathbbm{R}_+)}\le C \delta\mathrm{e}^{-c(1+t)}+C \delta(1+t)^{-\frac{1}{2}(k+2- \frac{1}{p} ) }.
   \end{aligned}
 \end{equation}
On the other hand, from Lemma \ref{lem-ghxsb} $\mathrm{(vi)}$, it also holds that
\begin{equation}\label{eq-R1lp2}
  \begin{aligned}
    \|\partial_x^k R_1(t)\|_{L^p(\mathbbm{R}_+)}\le C \delta\mathrm{e}^{-c(1+t)}+C (1+t)^{-\frac{1}{2}(k+3- \frac{1}{p} ) }.
  \end{aligned}
\end{equation}
Combining \eqref{eq-R1lp1} and \eqref{eq-R1lp2}, we complete the proof of $\mathrm{(ii)}$ in Corollary \ref{cor-R1R2}.

To get $\mathrm{(iii)}$, we can get from \eqref{eq-def-R1R2}$_2$ that
\begin{equation*}
  \begin{aligned}
    \| \partial_x^kR_2(t)\|_{L^p(\mathbbm{R}_+)}\le C \delta\mathrm{e}^{-c(1+t)}+\| \partial_x^{k+3} \tilde{u}\|_{L^p(\mathbbm{R}_+)}.
  \end{aligned}
\end{equation*}
 The estimate $\mathrm{(iii)}$ can be obtained by applying Lemma \ref{lem-jzxsb} $\mathrm{(iv)}$ and $\mathrm{(v)}$. Finally, we show the estimate $\mathrm{(viii)}$. From $\mathrm{(vi)}$,
 \begin{equation*}
   |\partial_x^k\partial_t^l R_{1}(0,t)|\le \| \partial_x^k\partial_t^l R_{1}(\cdot,t) \|_{L^\infty(\mathbbm{R}_+)} \le C \min \left\{\delta(1+t)^{-\frac{k+l+2}{2}}, (1+t)^{-\frac{k+l+3}{2}} \right\}.
 \end{equation*}
 $\hfill\Box$

We state the main result in this paper as follows.
\subsection{Main theorem}
\begin{thm}\label{thm-main} Assume that $0\le f'(u_-)<f'(u_+)$ holds. Suppose that $u_0(x,y)-r_0\in L^2(\mathbbm{R}_+^2)\cap L^1(\mathbbm{R}_+^2) $ and $\nabla u_0\in H^2(\mathbbm{R}_+^2)$, then there exists a positive constant $\delta_0$ such that
\begin{equation}
  \| u_0(x,y)-r_0 \|_{L^2(\mathbbm{R}_+^2)}  +\| \nabla{u}_0 \|_{H^2(\mathbbm{R}_+^2)} +|u_+-u_-|\le \delta_0,
\end{equation}
then the initial-boundary value problem \eqref{eq-yfc}-\eqref{eq-divp0tj} has a unique global solution $(u(x,y,t),q(x,y,t))$ which satisfies
\begin{equation*}
\begin{cases}
    u-r\in C^0([0,\infty);H^3(\mathbbm{R}_+^2)),\ \ \nabla{u}-r_x\in L^2(0,\infty;H^2(\mathbbm{R}_+^2)),\\[1mm]
  q+r_x\in C^0([0,\infty);H^3(\mathbbm{R}_+^2))\cap L^2(0,\infty;H^3(\mathbbm{R}_+^2)), \ \ {\rm{div}}{q}+r_{xx}\in C^0([0,\infty);H^3(\mathbbm{R}_+^2))\cap L^2(0,\infty;H^3(\mathbbm{R}_+^2)),
\end{cases}
\end{equation*}
and
\begin{equation}\label{eq-uq-largetime}
\begin{aligned}[b]
  &\sup\limits_{(x,y)\in\mathbbm{R}_+^2}|\nabla^k (u(x,y,t)-r(x,t))|\rightarrow 0 \quad \text{as}\quad t\rightarrow\infty,\quad k=0,1,\\
  &\sup\limits_{(x,y)\in\mathbbm{R}_+^2}|\nabla^k (q(x,y,t)+r_x(x,t))|\rightarrow 0 \quad \text{as}\quad t\rightarrow\infty,\quad k=0,1,\\
  &\sup\limits_{(x,y)\in\mathbbm{R}_+^2}|\nabla ({\rm{div}}{q}(x,y,t)+r_{xx}(x,t))|\rightarrow 0 \quad \text{as}\quad t\rightarrow\infty.	
\end{aligned}
\end{equation}
\end{thm}
\section{Asymptotics to planar rarefaction wave}\label{sec-3}
In this section, we focus on the case $(b): ~ 0\le f'(u_-)<f'(u_+)$. We denvote ourselves to show the asymptotic behaviors of the solution of \eqref{eq-yfc}-\eqref{eq-divp0tj} is the corresponding planar rarefaction wave as $t$ tends to infinity. In order to show the Theorem \ref{thm-main}, we just need to prove the following two Theorems \ref{thm-V-main} and \ref{thm-vp-main}.
\begin{thm}\label{thm-V-main} Assume that \eqref{eq-f-tj1}, \eqref{eq-f-tj3} and $0\le f'(u_-)<f'(u_+)$ holds. Suppose that $V_0\in H^4(\mathbbm{R}_+)\cap L^1(\mathbbm{R}_+)$ and there exists a small positive constant $\delta_1'$ such that $\| V_0 \|_{H^4(\mathbbm{R}_+)}+\delta\le \delta_1' $, then the problem \eqref{eq-rd-VP-1} has a unique global solution satisfying
\begin{equation*}
  \begin{cases}
    V\in C^0([0,\infty);H^4(\mathbbm{R}_+)) \cap C^0([0,\infty);L^1(\mathbbm{R}_+)),\\[1mm]
    P\in C^0([0,\infty);H^5(\mathbbm{R}_+))\cap L^{2}(0,\infty;H^5(\mathbbm{R}_+)),
  \end{cases}
\end{equation*}
and for sufficiently large $t$
\begin{equation*}
  \begin{cases}
    \|V(t) \|_{L^\infty(\mathbbm{R}_+)} \le C(1+t)^{-\frac{1}{2} }\mathrm{log}^3(2+t),\\[1mm]
    \| \partial_x^kV(t) \|_{L^\infty(\mathbbm{R}_+)} \le C(1+t)^{-\frac{3}{4} }\mathrm{log}^5(2+t),\quad k=1,2,3,\\[1mm]
    \| \partial_x^kP(t) \|_{L^\infty(\mathbbm{R}_+)} \le C(1+t)^{-\frac{3}{4} }\mathrm{log}^5(2+t),\quad k=0,1,2,3,4.
  \end{cases}
\end{equation*}
\end{thm}

\begin{thm}\label{thm-vp-main} Assume that \eqref{eq-f-tj1}, \eqref{eq-f-tj3} and $0\le f'(u_-)<f'(u_+)$ holds. Suppose that $v_0\in H^3(\mathbbm{R}_+^2)\cap L^1(\mathbbm{R}_+^2)$ and there exists a small positive constant $\delta_2'$ such that $\|v_0 \|_{H^3(\mathbbm{R}_+^2)}+\delta\le \delta_2' $, then the problem $\eqref{eq-rdfc-vp-1}$ has a unique global solution satisfying
\begin{equation*}
  \begin{cases}
    v\in C^0([0,\infty);H^3(\mathbbm{R}_+^2)),\quad \nabla{v}\in L^2(0,\infty;H^2(\mathbbm{R}_+^2)),\\[1mm]
    p, {\rm{div}}{p}\in C^0([0,\infty);H^3(\mathbbm{R}_+^2))\cap L^{2}(0,\infty;H^3(\mathbbm{R}_+^2)),
  \end{cases}
\end{equation*}
and
\begin{equation}\label{eq-largebehavior-v}
  \begin{cases}
    &\sup\limits_{(x,y)\in\mathbbm{R}_+^2} |\nabla^k v(x,y,t)|\rightarrow 0,\quad \text{as} \quad t\rightarrow \infty, \quad k=0,1,\\
    &\sup\limits_{(x,y)\in\mathbbm{R}_+^2} |\nabla^k p(x,y,t)|\rightarrow 0,\quad \text{as} \quad t\rightarrow \infty, \quad k=0,1,\\
    &\sup\limits_{(x,y)\in\mathbbm{R}_+^2} |\nabla {\rm{div}}{p}(x,y,t)|\rightarrow 0,\quad \text{as} \quad t\rightarrow \infty.
  \end{cases}
\end{equation}
\end{thm}

In the following, we try to prove the Theorem \ref{thm-V-main} and Theorem \ref{thm-vp-main}. We note that the main difference between $f'(u_-)>0$ and $f'(u_-)=0$ is the boundary value $f'(u_-) \partial_x^k V(0,t)$ and $f'(u_-) \nabla^l v(0,y,t)$, $k=1,2,3,4$, $l=1,2,3$. The proof of the case of $f'(u_-)=0$ is simpler than that of $f'(u_-)>0$. Thus, we only prove the case of $f'(u_-)>0$ and the proof of $f'(u_-)=0$ is omitted.

\subsection{Estimates for the perturbation on one-dimensional half space}\label{sec-xsb1}
In this section, we consider the initial-boundary value problem on one-dimensional half space:
\begin{equation}\label{eq-rd-VP}
\begin{cases}
  V_t+(f(V+\tilde{u})-f(\tilde{u}))_x+P_x=R_1,\\[1mm]
  -P_{xx}+P+V_x=R_2,
\end{cases}
\end{equation}
with initial data
\begin{equation}\label{eq-VP-initial}
   V(x,0)=V_0(x)=U_0(x)-\tilde{u}_0(x),
\end{equation}
and boundary condition
\begin{equation}\label{eq-VP-boundary}
  V(0,t)=0,\quad P_x(0,t)=0.
\end{equation}
Here, $R_1$ and $R_2$ is defined in \eqref{eq-def-R1R2}. From \eqref{eq-rd-VP}, we have
\begin{align}
	&(\partial_x^k\partial_t^l V_x)^2\le 3\left((\partial_x^k P_{xx})^2+(\partial_x^k P)^2+(\partial_x^k R_2)^2\right),\quad k+l=0,1,2,3,\quad k,l\in\mathbbm{N},\label{eq-VxP-rela}\\
    &(\partial_x^k\partial_t^l P_{xx})^2\le 3\left((\partial_x^k V_{x})^2+(\partial_x^k P)^2+(\partial_x^k R_2)^2\right),\quad k+l=0,1,2,3,\quad k,l\in\mathbbm{N},\label{eq-P-Vrela}
\end{align}
which will be often used later and plays an important role in \emph{a priori} estimates.

 The solution of the reformulated problem \eqref{eq-rd-VP}-\eqref{eq-VP-boundary} is sought in the set of functional space $X(0,T)$, where for $0\le T\le +\infty$, we define
\begin{equation*}
X(0,T) = \left\{
\begin{tabular}{c|c}
   \multirow{2}{*}{({\it V, P})}    & $V\in C^0([0,T);H^4(\mathbbm{R}_+)), \ \ V_x\in L^2(0,T;H^3(\mathbbm{R}_+))$   \\[1mm]
             & $P\in C^0([0,T);H^5(\mathbbm{R}_+))\cap L^2(0,T;H^5(\mathbbm{R}_+))$ \ \ \ \
\end{tabular}
       \right\},
\end{equation*}

\begin{prop}\label{prop-V}
Suppose that the boundary condition and far field states satisfy $0\le f'(u_-)<f'(u_+) $, the initial data $V_0\in H^4(\mathbbm{R}_+)$ and the wavelength $\delta=|u_--u_+|$ are sufficiently small.
Then there are the positive constants $\tilde{\delta}_1$ and $C=C(\tilde{\delta}_1)$ such that if $\| V_0 \|_{H^4(\mathbbm{R}_+)}  +\delta\le \tilde{\delta}_1$,
the problem \eqref{eq-rd-VP}-\eqref{eq-VP-boundary} admits a unique solution $(V(x,t),P(x,t))\in X(0,+\infty)$ satisfying
\begin{equation}\label{eq-V-propjg}
\|V(t)\|_{H^4(\mathbbm{R}_+)}^2  +\|P(t)\|_{H^5(\mathbbm{R}_+)}^2+\int_0^t \ (\| V_x(\tau)\|_{H^3(\mathbbm{R}_+)}^2+\| P(\tau)\|_{H^5(\mathbbm{R}_+)}^2 )  \,d{\tau}
\le C(\|V_0\|_{H^4(\mathbbm{R}_+)}^2+\delta^\frac{1}{2} ),\ \ \forall t\in[0,\infty).
\end{equation}
\end{prop}
Since the proof for the local-in-time existence and uniqueness of the solution to \eqref{eq-rd-VP}-\eqref{eq-VP-boundary} is standard, the details will be omitted. By Sobolev inequality and Lemma \ref{lem-jzxsb}, we get
\begin{equation}\label{eq-def-dl}
  \begin{cases}
    \sup\limits_{t\ge 0}\sum\limits_{k=0}^3 \| \partial_x^k V(t) \|_{L^\infty(\mathbbm{R}_+)}\le C \delta',\\
    \sup\limits_{t\ge 0}\sum\limits_{k=1}^3 \| \partial_x^k U(t) \|_{L^\infty(\mathbbm{R}_+)}\le C \delta',
  \end{cases}
\end{equation}
where $\delta'=(\|V_0\|_{H^4(\mathbbm{R}_+)}^2+\delta^\frac{1}{2} )^{\frac{1}{2} }$ is sufficiently small provided that $\tilde{\delta}_1$ is small enough.\\
\indent In other to prove Proposition \ref{prop-V}, we suffices to show the following $\emph{a priori}$ estimates.
\begin{prop}[\emph{A priori} estimates] \label{prop-V-priori}
Let $T$ be a positive constant. Assume that $0\le f'(u_-)<f'(u_+) $. Suppose that the problem \eqref{eq-rd-VP}-\eqref{eq-VP-boundary} has a unique solution $(V,P)\in X(0,T)$.
Then there exist two positive constants $\tilde{\delta}_2(\le \tilde{\delta}_1)$ and $C=C(\tilde{\delta}_2)$ such that if
$\|V_0\|_{H^4(\mathbbm{R}_+)}+\delta\le\tilde{\delta}_2$, then we have the estimate
\begin{equation*}
\|V(t)\|_{H^4(\mathbbm{R}_+)}^2  +\|P(t)\|_{H^5(\mathbbm{R}_+)}^2+\int_0^t \ (\| V_x(\tau)\|_{H^3(\mathbbm{R}_+)}^2+\| P(\tau)\|_{H^5(\mathbbm{R}_+)}^2 )  \,d{\tau}
\le C(\|V_0\|_{H^4(\mathbbm{R}_+)}^2+\delta^\frac{1}{2} ),\ \ \forall t\in[0,T].
\end{equation*}
\end{prop}

Before proving the $\emph{a priori}$ estimates, we first give some basic estimates.
\begin{lem}\label{lem-f-jf} There is a positive constant $C$ such that the following estimates hold:
\begin{flalign*}
&(\mathrm{A1}) \ \ \ \int_{\mathbbm{R}_+} (f(V+\tilde{u})-f(\tilde{u}))_{x}V \,\mathrm{d}x\ge \frac{\alpha}{2}\int_{\mathbbm{R}_+} \tilde{u}_x V^2 \,\mathrm{d}x;&\nonumber
\end{flalign*}
\begin{flalign*}
&(\mathrm{A2}) \ \ \ \int_{\mathbbm{R}_+} (f(V+\tilde{u})-f(\tilde{u}))_{xx}V_x \,\mathrm{d}x&\nonumber\\
&\ \ \ \ \ \ \ \ \ \ge \frac{3 \alpha}{2} \int_{\mathbbm{R}_+} \tilde{u}_x V_x^2 \,\mathrm{d}x- \frac{f'(u_-)}{2}V_x^2(0,t) -C \int_{\mathbbm{R}_+} |V_x|(| \tilde{u}_{xx} ||V|+\tilde{u}_x^2|V|+V_x^2 )  \,\mathrm{d}x;&\nonumber
\end{flalign*}
\begin{flalign*}
&(\mathrm{A3}) \ \ \ \int_{\mathbbm{R}_+} (f(V+\tilde{u})-f(\tilde{u}))_{xxx}V_{xx} \,\mathrm{d}x&\nonumber\\
&\ \ \ \ \ \ \ \ \ \ge \frac{5 \alpha}{2}\int_{\mathbbm{R}_+} \tilde{u}_x V_{xx}^2 \,\mathrm{d}x- \frac{f'(u_-)}{2}V_{xx}^2(0,t)&\nonumber\\
&\ \ \ \ \ \ \ \ \  -C \int_{\mathbbm{R}_+} |V_{xx}|\left\{(| \tilde{u}_{xxx} |+|\tilde{u}_x|^3+\tilde{u}_x|\tilde{u}_{xx}|)|V| +(| \tilde{u}_{xx} |+\tilde{u}_x^2+V_x^2 )|V_x|+|V_x||V_{xx}|\right\} \,\mathrm{d}x;&\nonumber
\end{flalign*}
\begin{flalign*}
&(\mathrm{A4}) \ \ \ \int_{\mathbbm{R}_+} (f(V+\tilde{u})-f(\tilde{u}))_{xxxt}V_{xxt} \,\mathrm{d}x&\nonumber\\
&\ \ \ \ \ \ \ \ \ \ge \frac{5 \alpha}{2}\int_{\mathbbm{R}_+} \tilde{u}_x V_{xxt}^2 \,\mathrm{d}x- \frac{f'(u_-)}{2}V_{xxt}^2(0,t)-C\int_{\mathbbm{R}_+} |V_x||V_{xxt}|^2\,\mathrm{d}x &\nonumber\\
&\ \ \ \ \ \ \ \ \ -C \int_{\mathbbm{R}_+} |V_{xxt}|\left\{(|V_x|^3+\tilde{u}_x^3+\tilde{u}_x|\tilde{u}_{xx}|+|\tilde{u}_{xx}||V_x|+|\tilde{u}_{xxx}|)|V_t|+(|V_x|^2+|\tilde{u}_x|^2+|\tilde{u}_{xx}|)|V_{xt}| \right\} \,\mathrm{d}x&\nonumber\\
&\ \ \ \ \ \ \ \ \ -C \int_{\mathbbm{R}_+} |V_{xxt}|\left\{(|\tilde{u}_t||V_x|^2+|\tilde{u}_{xt}||V_x|+|\tilde{u}_{xxt}|+|\tilde{u}_t||\tilde{u}_{xx}|+\tilde{u}_x^2| \tilde{u}_t|+\tilde{u}_x|\tilde{u}_{xt}|)|V_x|\right\} \,\mathrm{d}x&\nonumber\\
&\ \ \ \ \ \ \ \ \ -C \int_{\mathbbm{R}_+} |V_{xxt}|\left\{(|V_t||V_x|+|\tilde{u}_t||V_x|+|V_t||\tilde{u}_x|+|\tilde{u}_t||\tilde{u}_x|+|V_{xt}|+|\tilde{u}_{xt}|)|V_{xx}|+(|V_t|+|\tilde{u}_t|)||V_{xxx}|\right\} \,\mathrm{d}x&\nonumber\\
&\ \ \ \ \ \ \ \ \ -C \int_{\mathbbm{R}_+} |V_{xxt}|(|\tilde{u}_t|\tilde{u}_x^3+\tilde{u}_x^2|\tilde{u}_{xt}|+|\tilde{u}_t|\tilde{u}_x|\tilde{u}_{xx}|+|\tilde{u}_{xt}||\tilde{u}_{xx}|+\tilde{u}_x|\tilde{u}_{xxt}|+|\tilde{u}_t||\tilde{u}_{xxx}|+|\tilde{u}_{xxxt}|)|V| \,\mathrm{d}x;&\nonumber
\end{flalign*}
\begin{flalign*}
&(\mathrm{A5}) \ \ \ \int_{\mathbbm{R}_+} (f(V+\tilde{u})-f(\tilde{u}))_{xxxtt}V_{xxtt} \,\mathrm{d}x&\nonumber\\
&\ \ \ \ \ \ \ \ \ \ge \frac{ 5\alpha}{2}\int_{\mathbbm{R}_+} \tilde{u}_x V_{xxtt}^2 \,\mathrm{d}x-\frac{f'(u_-)}{2}V_{xxtt}^2(0,t)-C\int_{\mathbbm{R}_+} |V_{xxtt}|\left\{  |V_x||V_{xxtt}| +(|V_t|+|\tilde{u}_t|)|V_{xxxt}|\right\} \,\mathrm{d}x&\nonumber\\
&\ \ \ \ \ \ \ \ \  -C\int_{\mathbbm{R}_+} |V_{xxtt}|\left\{ (\tilde{u}_t^2\!+\! |\tilde{u}_{tt}|)(V_x^2\!+\!\tilde{u}_x^2)\!+\!(|\tilde{u}_t \tilde{u}_{xt}|\!+\!|\tilde{u}_{xtt}|)(|V_x|\!+\!\tilde{u}_x)\!+\!V_{xt}^2\!+\!|\tilde{u}_{xx}|V_t^2\!+\!|\tilde{u}_{xxt}|(|V_t|\!+\!|\tilde{u}_t|)  \right\}|V_x| \,\mathrm{d}x &\nonumber\\
&\ \ \ \ \ \ \ \ \  -C\int_{\mathbbm{R}_+} |V_{xxtt}|\left\{ (|\tilde{u}_{xxtt}|+|\tilde{u}_{xx}|\tilde{u}_t^2+\tilde{u}_{xt}^2+|\tilde{u}_{tt}\tilde{u}_{xx}|)  |V_x|+(|\tilde{u}_{xx}\tilde{u}_{xt}|+\tilde{u}_x|\tilde{u}_{xxt}|+|\tilde{u}_{xxxt}|)|V_t|\right\} \,\mathrm{d}x &\nonumber\\
&\ \ \ \ \ \ \ \ \  -C\int_{\mathbbm{R}_+} |V_{xxtt}|\left\{ (|V_x|^3+\tilde{u}_x^3+\tilde{u}_x |\tilde{u}_{xx}|+|\tilde{u}_{xxx}|)(|V_t|+|\tilde{u}_t|)+|\tilde{u}_{xt}|(V_x^2+\tilde{u}_x^2) \right\}|V_t| \,\mathrm{d}x &\nonumber\\
&\ \ \ \ \ \ \ \ \  -C\int_{\mathbbm{R}_+} |V_{xxtt}|\left\{ (|V_x|\!+\!\tilde{u}_x)(V_{t}^2\!+\!\tilde{u}_t^2)\!+\!(|V_x|\!+\!\tilde{u}_x)(|V_{tt}|\!+\!|\tilde{u}_{tt}|)\!+\!(|V_t|+|\tilde{u}_t|)(|V_{xt}|\!+\!|\tilde{u}_{xt}|)\!+\!|\tilde{u}_{xtt}| \right\}|V_{xx}| \,\mathrm{d}x &\nonumber\\
&\ \ \ \ \ \ \ \ \  -C\int_{\mathbbm{R}_+} |V_{xxtt}|\left\{ (|V_t|+|\tilde{u}_t|)(V_{x}^2+\tilde{u}_x^2)+|\tilde{u}_{xx}|(|V_{t}|+|\tilde{u}_{t}|)+\tilde{u}_x(|V_{xt}|+|\tilde{u}_{xt}|)+|\tilde{u}_{xxt}| \right\}|V_{xt}| \,\mathrm{d}x &\nonumber\\
&\ \ \ \ \ \ \ \ \  -C\int_{\mathbbm{R}_+} |V_{xxtt}|\left\{  (|V_x|^3+\tilde{u}_x^3+|\tilde{u}_{xx}V_x|+|\tilde{u}_{xx}|\tilde{u}_x)|V_{tt}| +(V_x^2+\tilde{u}_x^2+|V_{xx}|+|\tilde{u}_{xx}|)|V_{xtt}|\right\} \,\mathrm{d}x &\nonumber\\
&\ \ \ \ \ \ \ \ \  -C\int_{\mathbbm{R}_+} |V_{xxtt}|\left\{  (V_t^2+\tilde{u}_t^2+|V_{tt}|+|\tilde{u}_{tt}|)|V_{xxx}| +(|V_xV_t|+\tilde{u}_x|V_t|+|\tilde{u}_tV_x|+\tilde{u}_x|\tilde{u}_t|+|V_{xt}|+|\tilde{u}_{xt}|)|V_{xxt}|\right\} \,\mathrm{d}x &\nonumber\\
&\ \ \ \ \ \ \ \ \  -C \int_{\mathbbm{R}_+} |V_{xxtt}|\left\{ \tilde{u}_x|\tilde{u}_{xxtt}|+\tilde{u}_x^3|\tilde{u}_{tt}|+\tilde{u}_x^2|\tilde{u}_{xt}|+\tilde{u}_x\tilde{u}_{xt}^2+\tilde{u}_t^2 \tilde{u}_x|\tilde{u}_{xx}|+\tilde{u}_x^2|\tilde{u}_{xtt}|+|\tilde{u}_{tt}||\tilde{u}_{xxx}|+\tilde{u}_x|\tilde{u}_{tt}||\tilde{u}_{xx}| \right\}|V| \,\mathrm{d}x&\nonumber\\
&\ \ \ \ \ \ \ \ \  -C \int_{\mathbbm{R}_+} |V_{xxtt}|\left\{ |\tilde{u}_{t}\tilde{u}_{xt}\tilde{u}_{xx}|+\tilde{u}_x|\tilde{u}_{t}\tilde{u}_{xxt}|+|\tilde{u}_{xx}\tilde{u}_{xtt}|+|\tilde{u}_{xxt}\tilde{u}_{xt}|+\tilde{u}_t^2|\tilde{u}_{xxx}|+|\tilde{u}_{t}\tilde{u}_{xxxt}|+|\tilde{u}_{xxxtt}| \right\}|V| \,\mathrm{d}x.&\nonumber
\end{flalign*}
\end{lem}
 {\it\bfseries Proof.}
By direct calculation, we can obtain
\begin{equation*}
	\begin{aligned}
	(f(V+\tilde{u})-f(\tilde{u}))_{x}V\!=\!(f(V\!+\!\tilde{u})\!-\!f(\tilde{u})\!-\!f'(\tilde{u})V )\tilde{u}_x\!+\!\left\{ (f(V\!+\!\tilde{u})\!-\!f(\tilde{u}))V\!-\! \int_{\tilde{u}}^{V\!+\!\tilde{u}}f(s) \,\mathrm{d}s\!+\!f(\tilde{u})V\right\}_x,
	\end{aligned}
\end{equation*}
and
\begin{equation*}
  \begin{aligned}[b]
    (f(V+\tilde{u})-f(\tilde{u}))_{xx}V_x
    =&\frac{1}{2}f''(V+\tilde{u})V_x^3+\frac{3}{2}f''(v+\tilde{u})\tilde{u}_x V_x^2+\frac{1}{2}\{ f'(V+\tilde{u})V_x^2\}_x\\
    &+(f''(V+\tilde{u})-f''(\tilde{u})) \tilde{u}_x^2V_x+(f'(V+\tilde{u})-f'(\tilde{u}))\tilde{u}_{xx}V_x.
  \end{aligned}
\end{equation*}
Integrating the above two equations over $\mathbbm{R}_+$, from $V(0,t)=0$, we have
\begin{equation*}
	\int_{\mathbbm{R}_+} (f(V+\tilde{u})-f(\tilde{u}))_{x}V \,\mathrm{d}x \ge \frac{\alpha}{2} \int_{\mathbbm{R}_+} \tilde{u}_xV^2 \,\mathrm{d}x,
\end{equation*}
and
\begin{equation*}
  \begin{aligned}
    &\int_{\mathbbm{R}_+} (f(V+\tilde{u})-f(\tilde{u}))_{xx}V_x \,\mathrm{d}x \\
    &\ge \frac{3 \alpha}{2} \int_{\mathbbm{R}_+} \tilde{u}_x V_x^2 \,\mathrm{d}x- \frac{f'(u_-)}{2}V_x^2(0,t) -C \int_{\mathbbm{R}_+} |V_x|(| \tilde{u}_{xx} ||V|+\tilde{u}_x^2|V|+V_x^2 )  \,\mathrm{d}x,
  \end{aligned}
\end{equation*}
which yields the desired estimates $(\mathrm{A1})$ and $(\mathrm{A2})$. The inequalities $(\mathrm{A3})$-$(\mathrm{A5})$ can be obtained by using the similar way.
 $ \hfill\Box $
\begin{lem}\label{lem-bj-gj} Suppose that $f(u_-)>0$. The solution $V(x,t)$ satisfies the following boundary estimates:
\begin{flalign*}
\hspace{1mm}
&(\mathrm{B1})\ \ \ V_x^2(0,t)\le C \delta \mathrm{e}^{-c(1+t)}; &\nonumber\\
\hspace{1mm}
&(\mathrm{B2})\ \ \ V_{xt}^2(0,t)\le C \delta\mathrm{e}^{-c(1+t)}; &\nonumber\\
&(\mathrm{B3})\ \ \ V_{xx}^2(0,t)\le C \delta\mathrm{e}^{-c(1+t)}+C P_{xx}^2(0,t); &\nonumber\\
&(\mathrm{B4})\ \ \ V_{xtt}^2(0,t)\le C \delta\mathrm{e}^{-c(1+t)}; &\nonumber\\
&(\mathrm{B5})\ \ \ V_{xxt}^2(0,t)\le C \delta \mathrm{e}^{-c(1+t)}+C P_{xxt}^2(0,t); &\nonumber\\
&(\mathrm{B6})\ \ \ V_{xttt}^2(0,t)\le C \delta\mathrm{e}^{-c(1+t)}; &\nonumber\\
&(\mathrm{B7})\ \ \ V_{xxtt}^2(0,t)\le C \delta\mathrm{e}^{-c(1+t)}+C P_{xxtt}^2(0,t). &\nonumber\\
\end{flalign*}
\end{lem}
 {\it\bfseries Proof.}
 From \eqref{eq-rd-VP}$_1$ and \eqref{eq-VP-boundary}, it holds that
\begin{equation}\label{eq-Vx0t}
  f'(u_-)V_x(0,t)=R_1(0,t)-P_x(0,t)=R_1(0,t).
\end{equation}
By utilizing Corollary \ref{cor-R1R2} $\mathrm{(iv)}$, we have
\begin{equation*}
  V_x^2(0,t)\le C \delta\mathrm{e}^{-c(1+t)}.
\end{equation*}
Differentiating \eqref{eq-Vx0t} with respect to $t$, we get
\begin{equation*}
  f'(u_-)V_{xt}(0,t)=R_{1t}(0,t).
\end{equation*}
It follows that
\begin{equation*}
  V_{xt}^2(0,t)\le C R_{1t}^2(0,t)\le C \delta\mathrm{e}^{-c(1+t)}.
\end{equation*}
Similarly, we can get the remaining boundary estimates $(\mathrm{B3})$-$(\mathrm{B7})$ by utilizing $\eqref{eq-rd-VP}_1$.
 $ \hfill\Box $

 \begin{remark} For the case of $f'(u_-)=0$, the boundary terms $(\mathrm{B1})$-$(\mathrm{B7})$ in Lemma \ref{lem-bj-gj} will disappear. Because the coefficients of all boundary terms $(\partial_x^k \partial_t^lV(0,y,t))^2$ are $\frac{f'(u_-)}{2} $ which is given in Lemma \ref{lem-f-jf} $(\mathrm{A2})$-$(\mathrm{A5})$. Similarly, the boundary terms $(\mathrm{D3})$-$(\mathrm{D5})$ in Lemma \ref{lem-bjv-gj} will also disappear. In other words, in the case of $f'(u_-)=0$, the results of the initial boundary value problem \eqref{eq-yfc}-\eqref{eq-divp0tj} is the same as that of the Cauchy problem. For details, we refer to \cite{Gao2008}.
 \end{remark}

\subsubsection{\emph{A priori} estimates}
In this section, we will show Proposition \ref{prop-V-priori} under the $\emph{a priori}$ assumption
\begin{equation*}
  \| V_x(t) \|_{H^4(\mathbbm{R}_+)} \le \varepsilon_0,
\end{equation*}
where $0<\varepsilon_0\ll 1$.
By Sobolev inequality, there exists a positive constant $C$ such that
\begin{equation*}
  \|\partial_x^k V(t) \|_{L^\infty(\mathbbm{R}_+)} \le C \varepsilon_0,\quad k=0,1,2,3.
\end{equation*}
 For simplicity, we divide the proof of the \emph{a priori} estimates into several lemmas.
\begin{lem}\label{lem-V-xsb} Under the same assumption as Proposition \ref{prop-V-priori}, there exists a positive constant $C$ such that
  \begin{equation}\label{eq-V-L2-jg}
    \begin{aligned}[b]
      \| V(t) \|_{L^2(\mathbbm{R}_+)}^2+\int_0^t (\| \sqrt{\tilde{u}_x}V(\tau) \|_{L^2(\mathbbm{R}_+)}^2+\| P_{x}(\tau) \|_{L^2(\mathbbm{R}_+)}^2+\| P(\tau) \|_{L^2(\mathbbm{R}_+)}^2)  \,\mathrm{d}\tau \le C (\| V_0 \|_{H^4(\mathbbm{R}_+)}^2+\delta^\frac{1}{2}  ).
    \end{aligned}
  \end{equation}
\end{lem}
 {\it\bfseries Proof.}
We can get from $\eqref{eq-rd-VP}_1\times V+\eqref{eq-rd-VP}_2\times P$ that
\begin{equation}\label{eq-V-L2}
  \begin{aligned}[b]
    \frac{1}{2}\frac{\mathrm{d}}{\mathrm{d}t} V^2&+(f(V+\tilde{u})-f(\tilde{u}))_xV+P_x^2+P^2+\left\{ PV- P_xP\right\}_x=R_1V+R_2P.
  \end{aligned}
\end{equation}
For $V(0,t)=0$ and $P_x(0,t)=0$, the terms in $\{\cdots\}_x$ disappear after integration in $x\in \mathbbm{R}_+$. Thus, integrating $\eqref{eq-V-L2}$ over $\mathbbm{R}_+$, by Lemma \ref{lem-f-jf} $(\mathrm{A1})$, we have
\begin{equation}\label{eq-V-L2-1}
  \begin{aligned}
    &\frac{1}{2} \frac{\mathrm{d}}{\mathrm{d}t}\| V(t) \|_{L^2(\mathbbm{R}_+)}^2  +\frac{\alpha}{2} \| \sqrt{\tilde{u}_x}V(t) \|_{L^2(\mathbbm{R}_+)}^2 +\|P_x(t)\|_{L^2(\mathbbm{R}_+)}^2+\|P(t)\|_{L^2(\mathbbm{R}_+)}^2\le C \int_{\mathbbm{R}_+} (|R_1V|+|R_2P|) \,\mathrm{d}x.
  \end{aligned}
\end{equation}
We treat the terms on the right-hand side of $\eqref{eq-V-L2-1}$ as follows. From Corollary \ref{cor-R1R2} $\mathrm{(i)}$ and $\mathrm{(iii)}$,
\begin{equation}\label{eq-V-L2-2}
  \begin{aligned}[b]
    &C\int_{\mathbbm{R}_+} |R_1V| \,\mathrm{d}x \le C\| R_1(t) \|_{L^2(\mathbbm{R}_+)} \| V(t) \|_{L^2(\mathbbm{R}_+)}\le C\left(\delta^{\frac{1}{2} }(1+t)^{-\frac{3}{2} }+\delta \mathrm{e}^{-c(1+t)}\right)\left(1+\| V(t) \|_{L^2(\mathbbm{R}_+)}^2\right),\\
    &C\int_{\mathbbm{R}_+} |R_2P| \,\mathrm{d}x \le \frac{1}{4} \| P(t) \|_{L^2(\mathbbm{R}_+)}^2+C \| R_2(t) \|_{L^2(\mathbbm{R}_+)}^2\le \frac{1}{4} \| P(t) \|_{L^2(\mathbbm{R}_+)}^2+C \delta^2(1+t)^{-\frac{5}{2} }+C \delta \mathrm{e}^{-c(1+t)}.
  \end{aligned}
\end{equation}
Substituting \eqref{eq-V-L2-2} into \eqref{eq-V-L2-1}, we can deduce that
\begin{equation}\label{eq-V-L2-3}
  \begin{aligned}[b]
     &\frac{\mathrm{d}}{\mathrm{d}t}\| V(t) \|_{L^2(\mathbbm{R}_+)}^2  +\| \sqrt{\tilde{u}_x}V(t) \|_{L^2(\mathbbm{R}_+)}^2 +\|P_x(t)\|_{L^2(\mathbbm{R}_+)}^2+\|P(t)\|_{L^2(\mathbbm{R}_+)}^2\\
     \le& C\delta^{\frac{1}{2} }\left((1+t)^{-\frac{3}{2} }+ \mathrm{e}^{-c(1+t)}\right)\left(1+\| V(t) \|_{L^2(\mathbbm{R}_+)}^2\right).
  \end{aligned}
\end{equation}
Integrating $\eqref{eq-V-L2-3}$ over $[0,t]$, for some small $\delta$, we can obtain the desired estimate \eqref{eq-V-L2-jg}.
 $\hfill\Box$\\
\begin{lem}\label{lem-Vx-jg} Under the same assumptions as Proposition \ref{prop-V-priori}, there exists a positive constant $C$ such that
  \begin{equation}\label{eq-V-H1-jg}
    \begin{aligned}
      \| V_x(t) \|_{L^2(\mathbbm{R}_+)}^2+\int_0^t (\| \sqrt{\tilde{u}_x}V_x(\tau) \|_{L^2(\mathbbm{R}_+)}^2+\| P_{xx}(\tau) \|_{L^2(\mathbbm{R}_+)}^2+\| P_{x}(\tau) \|_{L^2(\mathbbm{R}_+)}^2  ) \,\mathrm{d}\tau\le C(\| V_0 \|_{H^4(\mathbbm{R}_+)}^2+\delta^\frac{1}{2}  ).
    \end{aligned}
  \end{equation}
\end{lem}
 {\it\bfseries Proof.}
We can obtain from $\partial_x\eqref{eq-rd-VP}_1\times V_x- \eqref{eq-rd-VP}_2\times P_{xx}$ that
\begin{equation}\label{eq-V-H1-1}
  \begin{aligned}[b]
    \frac{1}{2} \frac{\mathrm{d}}{\mathrm{d}t}V_x^2+(f(V+\tilde{u})-f(\tilde{u}))_{xx}V_x+P_{xx}^2+P_x^2-\{P_xP\}_x
    =R_{1x}V_x-R_{2}P_{xx}.
  \end{aligned}
\end{equation}
From Lemma \ref{lem-f-jf} $\mathrm{(A2)}$, integrating $\eqref{eq-V-H1-1}$ over $\mathbbm{R}_+$, we have
\begin{equation}\label{eq-V-H1-2}
\begin{aligned}[b]
  &\frac{1}{2} \frac{\mathrm{d}}{\mathrm{d}t}\| V_x(t) \|_{L^2(\mathbbm{R}_+)}^2+\frac{3 \alpha}{2}\| \sqrt{\tilde{u}_x}V_x(t) \|_{L^2(\mathbbm{R}_+)}^2   +\| P_{xx}(t) \|_{L^2(\mathbbm{R}_+)}^2+\| P_x(t) \|_{L^2(\mathbbm{R}_+)}^2 \\
  \le& \frac{f'(u_-)}{2}V_x^2(0,t) +C \int_{\mathbbm{R}_+} |V_x|(| \tilde{u}_{xx} ||V|+\tilde{u}_x^2|V|+V_x^2 )  \,\mathrm{d}x+\int_{\mathbbm{R}_+} (|R_{1x}V_x|+|R_{2}P_{xx}|) \,\mathrm{d}x.
\end{aligned}
\end{equation}
Now we estimate the terms on the right-hand side of $\eqref{eq-V-H1-2}$, successively. By applying $\eqref{eq-f-parx}$ and Lemma \ref{lem-jzxsb} $\mathrm{(iv)}$ with $k=2$, the next term can be estimated as
\begin{equation}\label{eq-V-H1-4}
  \begin{aligned}[b]
    C \int_{\mathbbm{R}_+} |V_x|(| \tilde{u}_{xx} ||V|+ \tilde{u}_{x}^2 |V|+V_x^2 )  \,\mathrm{d}x
    \le & (\| V_x(t) \|_{L^\infty(\mathbbm{R}_+)}+\mu ) \| V_x(t) \|_{L^2(\mathbbm{R}_+)}^2+C \mu^{-1}\| \tilde{u}_{xx}(t) \|_{L^\infty(\mathbbm{R}_+)}^2 \| V(t) \|_{L^2(\mathbbm{R}_+)}^2 \\
    \le & (\varepsilon_0+\mu ) \| V_x(t) \|_{L^2(\mathbbm{R}_+)}^2+C \mu^{-1} \delta^2(1+t)^{-2 }.
  \end{aligned}
\end{equation}
The last two terms of \eqref{eq-V-H1-2} are bounded by
\begin{equation}\label{eq-V-H1-5}
  \begin{aligned}[b]
     &\int_{\mathbbm{R}_+} (|R_{1x}V_x|+|R_2P_{xx}|) \,\mathrm{d}x\\
     \le& \mu\left(\| V_x(t) \|_{L^2(\mathbbm{R}_+)}^2 +\| P_{xx}(t) \|_{L^2(\mathbbm{R}_+)}^2 \right)+ C \mu^{-1}\left(\| R_{1x}(t) \|_{L^2(\mathbbm{R}_+)}^2+\| R_{2}(t) \|_{L^2(\mathbbm{R}_+)}^2 \right)\\
     \le &\mu\left(\| V_x(t) \|_{L^2(\mathbbm{R}_+)}^2 +\| P_{xx}(t) \|_{L^2(\mathbbm{R}_+)}^2 \right)+ C \mu^{-1} \delta^2(1+t)^{-\frac{5}{2} }.
  \end{aligned}
\end{equation}
Substituting \eqref{eq-V-H1-4}-\eqref{eq-V-H1-5} into \eqref{eq-V-H1-2}, by employing Lemma \ref{lem-bj-gj} $\mathrm{(B1)}$, for some small but fixed $\mu$, we have
\begin{equation}\label{eq-H1-djf0}
  \begin{aligned}[b]
    &\frac{\mathrm{d}}{\mathrm{d}t} \| V_x(t) \|_{L^2(\mathbbm{R}_+)}^2+\| \sqrt{\tilde{u}_x}V_x(t) \|_{L^2(\mathbbm{R}_+)}^2 +\| P_{xx}(t) \|_{L^2(\mathbbm{R}_+)}^2+\| P_x(t) \|_{L^2(\mathbbm{R}_+)}^2 \\
    \le& C(\varepsilon_0+\mu ) \| V_x(t) \|_{L^2(\mathbbm{R}_+)}^2+C \delta^2(1+t)^{-2 }+ C \delta\mathrm{e}^{-c(1+t)}.
  \end{aligned}
\end{equation}
The inequality \eqref{eq-VxP-rela} with $k=l=0$ gives
\begin{equation}\label{eq-Vx-P}
  \| V_x(t) \|_{L^2(\mathbbm{R}_+)}^2\le 3(\| P_{xx} (t)\|_{L^2(\mathbbm{R}_+)}^2+\| P(t) \|_{L^2(\mathbbm{R}_+)}^2+\| R_2(t) \|_{L^2(\mathbbm{R}_+)}^2 ) ,
\end{equation}
it follows from \eqref{eq-H1-djf0} that
\begin{equation}\label{eq-H1-djf1}
  \begin{aligned}[b]
    &\frac{\mathrm{d}}{\mathrm{d}t} \| V_x(t) \|_{L^2(\mathbbm{R}_+)}^2+\| \sqrt{\tilde{u}_x}V_x(t) \|_{L^2(\mathbbm{R}_+)}^2 +\| P_{xx}(t) \|_{L^2(\mathbbm{R}_+)}^2+\| P_x(t) \|_{L^2(\mathbbm{R}_+)}^2 \\
    \le&   C\| P(t) \|_{L^2(\mathbbm{R}_+)}^2+C \delta^2(1+t)^{-2 }+ C \delta\mathrm{e}^{-c(1+t)}.
  \end{aligned}
\end{equation}
Integrating \eqref{eq-H1-djf1} over $[0,t]$, combining \eqref{eq-V-L2-jg}, we get the desired estimates \eqref{eq-V-H1-jg}.
  $\hfill\Box$\\
\indent Combining \eqref{eq-V-L2-jg}, \eqref{eq-V-H1-jg} and \eqref{eq-Vx-P}, we can easily obtain the following Corollary \ref{cor-Vx}.
\begin{cor}\label{cor-Vx} Under the same assumptions as Proposition \ref{prop-V-priori}, there exists a positive constant $C$ such that
  \begin{equation}\label{eq-Vx-cjf-jg}
    \int_0^t \| V_x(t) \|_{L^2(\mathbbm{R}_+)}^2  \,\mathrm{d}\tau\le C(\| V_0 \|_{H^4(\mathbbm{R}_+)}^2+\delta^\frac{1}{2}).
  \end{equation}
\end{cor}
\begin{lem}\label{lem-PH2-jg} Under the same assumptions as Proposition \ref{prop-V-priori}, there exists a positive constant $C$ such that
\begin{equation}
	\begin{aligned}[b]
		\| P(t) \|_{H^2(\mathbbm{R}_+)}^2\le  C(\| V_0 \|_{H^4(\mathbbm{R}_+)}^2+\delta^\frac{1}{2}).
	\end{aligned}
\end{equation}
\end{lem}
{\it\bfseries Proof.}
Rewriting $\eqref{eq-rd-VP}_2$ in the form $P_{xx}-P=V_x-R_1$ and squaring this equation, we get
\begin{equation}\label{eq-P-xsb}
  P_{xx}^2+2P_x^2+P^2-2\{P_xP\}_x=V_x^2+R_1^2-2V_xR_1.
\end{equation}
 Integrating $\eqref{eq-P-xsb}$ over $\mathbbm{R}_+$, from $P_x(0,t)=0$, combining $\eqref{eq-Vx-cjf-jg}$ and Corollary \ref{cor-R1R2} $\mathrm{(i)}$, we obtain
 \begin{equation}\label{eq-P-xsb-1}
   \| P(t) \|_{H^2(\mathbbm{R}_+)}^2\le 2 (\| V_x(t) \|_{L^2(\mathbbm{R}_+)}^2+\| R_1(t) \|_{L^2(\mathbbm{R}_+)}^2 )\le C(\| V_0 \|_{H^4(\mathbbm{R}_+)}^2 +\delta^{\frac{1}{2} }),
 \end{equation}
 which completes the proof of Lemma \ref{lem-PH2-jg}.
$ \hfill\Box$

\begin{lem}\label{lem-Vxx-jg} Under the same assumptions as Proposition \ref{prop-V-priori}, there exists a positive constant $C$ such that
\begin{equation}\label{eq-V-H2-jg}
  \begin{aligned}[b]
    \| V_{xx}(t) \|_{L^2(\mathbbm{R}_+)}^2&+\| P_{xxx}(t) \|_{L^2(\mathbbm{R}_+)}^2+\int_0^t (\| \sqrt{\tilde{u}_x}V_{xx}(\tau) \|_{L^2(\mathbbm{R}_+)}^2+\| P_{xx}(\tau) \|_{H^1(\mathbbm{R}_+)}^2+\| V_{xx}(\tau) \|_{L^2(\mathbbm{R}_+)}^2)   \,\mathrm{d}\tau\\
    &\le C(\| V_0 \|_{H^4(\mathbbm{R}_+)}^2+\delta^\frac{1}{2}  ).
  \end{aligned}
\end{equation}
\end{lem}
 {\it\bfseries Proof.}
 We can get from $\partial_x^2\eqref{eq-rd-VP}_1\times V_{xx}-\partial_x \eqref{eq-rd-VP}_2\times P_{xxx}$ that
 \begin{equation}\label{eq-V-H2-1}
   \begin{aligned}[b]
     \frac{1}{2} \frac{\mathrm{d}}{\mathrm{d}t}V_{xx}^2+(f(V+\tilde{u})-f(\tilde{u}))_{xxx}V_{xx}+P_{xxx}^2+P_{xx}^2-\{P_{xx}P_x\}_x=R_{1xx}V_{xx}-R_{2x}P_{xxx}.
   \end{aligned}
 \end{equation}
 Integrating \eqref{eq-V-H2-1} over $\mathbbm{R}_+$, by employing Lemma \ref{lem-f-jf} $\mathrm{(A3)}$, we obtain
 \begin{equation}\label{eq-V-H2-2}
   \begin{aligned}[b]
     &\frac{1}{2} \frac{\mathrm{d}}{\mathrm{d}t}\| V_{xx}(t) \|_{L^2(\mathbbm{R}_+)}^2+\frac{5 \alpha}{2}\| \sqrt{\tilde{u}_x}V_{xx}(t) \|_{L^2(\mathbbm{R}_+)}^2 +\| P_{xxx}(t) \|_{L^2(\mathbbm{R}_+)}^2+\| P_{xx}(t) \|_{L^2(\mathbbm{R}_+)}^2\\
     \le& \frac{f'(u_-)}{2}  V_{xx}^2(0,t)+C \int_{\mathbbm{R}_+} |V_{xx}|\left\{(| \tilde{u}_{xxx} |+|\tilde{u}_x|^3+\tilde{u}_x|\tilde{u}_{xx}|)|V| +(| \tilde{u}_{xx} |+\tilde{u}_{x}^2+V_x^2 )|V_x|+|V_x||V_{xx}|\right\} \,\mathrm{d}x\\
     &+\int_{\mathbbm{R}_+} (|R_{1xx}V_{xx}|+|R_{2x}P_{xxx}|) \,\mathrm{d}x.
   \end{aligned}
 \end{equation}
The terms on the right-hand side of \eqref{eq-V-H2-2} can be estimated as follows. By using Lemma \ref{lem-bj-gj} $\mathrm{(B3)}$, the first term can be estimated as
\begin{equation}\label{eq-V-H2-3}
  \begin{aligned}
    \frac{f'(u_-)}{2} V_{xx}^2(0,t)\le C \delta\mathrm{e}^{-c(1+t)}+\frac{1}{4} \| P_{xxx}(t) \|_{L^2(\mathbbm{R}_+)}^2+C \| P_{xx}(t) \|_{L^2(\mathbbm{R}_+)}^2.
  \end{aligned}
\end{equation}
Using Cauchy inequality and \eqref{eq-f-parx}, we have
\begin{equation}\label{eq-V-H2-4}
  \begin{aligned}[b]
   &C \int_{\mathbbm{R}_+} |V_{xx}|\left\{(| \tilde{u}_{xxx} |+|\tilde{u}_x|^3+\tilde{u}_x|\tilde{u}_{xx}|)|V| +(| \tilde{u}_{xx} |+\tilde{u}_{x}^2+V_x^2 )|V_x|+|V_x||V_{xx}|\right\} \,\mathrm{d}x\\
   \le& (\| V_x \|_{L^\infty(\mathbbm{R}_+)} +\mu )\| V_{xx}(t) \|_{L^2(\mathbbm{R}_+)}^2+C \mu^{-1}\| V_x(t) \|_{L^2(\mathbbm{R}_+)}^2+C  \mu^{-1}\| \tilde{u}_{xxx}(t) \|_{L^\infty(\mathbbm{R}_+)}^2 \| V(t) \|_{L^2(\mathbbm{R}_+)}^2\\
   \le& C(\varepsilon_0+\mu)\| V_{xx}(t) \|_{L^2(\mathbbm{R}_+)}^2+C \mu^{-1}\| V_x(t) \|_{L^2(\mathbbm{R}_+)}^2+C \delta^2(1+t)^{-3},
  \end{aligned}
\end{equation}
and
\begin{equation}\label{eq-V-H2-5}
  \begin{aligned}[b]
    &\int_{\mathbbm{R}_+} |R_{1xx}V_{xx}|+|R_{2x}P_{xxx}| \,\mathrm{d}x\\
    \le& \mu(\| V_{xx}(t) \|_{L^2(\mathbbm{R}_+)}^2+\| P_{xxx}(t) \|_{L^2(\mathbbm{R}_+)}^2 ) +C \mu^{-1}(\| R_{1xx}(t) \|_{L^2(\mathbbm{R}_+)}^2+\| R_{2x}(t) \|_{L^2(\mathbbm{R}_+)}^2 )\\
    \le& \mu(\| V_{xx}(t) \|_{L^2(\mathbbm{R}_+)}^2+\| P_{xxx}(t) \|_{L^2(\mathbbm{R}_+)}^2 ) +C \mu^{-1} \delta^2(1+t)^{-\frac{7}{2} }+C \mu^{-1} \delta \mathrm{e}^{-c(1+t)}.
  \end{aligned}
\end{equation}
Substituting \eqref{eq-V-H2-3}-\eqref{eq-V-H2-5} into \eqref{eq-V-H2-2}, choosing small but fixed $\mu$, we have
\begin{equation}\label{eq-V-H2-5-1}
  \begin{aligned}[b]
     &\frac{\mathrm{d}}{\mathrm{d}t}\| V_{xx}(t) \|_{L^2(\mathbbm{R}_+)}^2+\| \sqrt{\tilde{u}_x}V_{xx}(t) \|_{L^2(\mathbbm{R}_+)}^2+\| P_{xxx}(t) \|_{L^2(\mathbbm{R}_+)}^2+\| P_{xx}(t) \|_{L^2(\mathbbm{R}_+)}^2\\
     \le&(\varepsilon_0+\mu)\| V_{xx}(t) \|_{L^2(\mathbbm{R}_+)}^2+C \left(\| V_x(t) \|_{L^2(\mathbbm{R}_+)}^2+\| P_{xx}(t) \|_{L^2(\mathbbm{R}_+)}^2 \right)+C \delta^2(1+t)^{-3 }+C \delta\mathrm{e}^{-c(1+t)}.\\
  \end{aligned}
\end{equation}
From \eqref{eq-VxP-rela} with $k=1$ and $l=0$, it follows that
\begin{equation}\label{eq-rela-Vxx-P}
  \| V_{xx}(t) \|_{L^2(\mathbbm{R}_+)}^2\le 3(\| P_{xxx}(t) \|_{L^2(\mathbbm{R}_+)}^2+\| P_{x}(t) \|_{L^2(\mathbbm{R}_+)}^2+\| R_{2x}(t) \|_{L^2(\mathbbm{R}_+)}^2 ).
\end{equation}
Substitute \eqref{eq-rela-Vxx-P} into \eqref{eq-V-H2-5-1}, and then integrate the resulting inequality over $[0,t]$. Consequently, combining \eqref{eq-V-H1-jg} and \eqref{eq-Vx-cjf-jg}, for some small $\delta$, $\mu$ and $\varepsilon_0$, we obtain
\begin{equation*}
  \begin{aligned}[b]
    &\| V_{xx}(t) \|_{L^2(\mathbbm{R}_+)}^2+\int_0^t (\| \sqrt{\tilde{u}_x}V_{xx}(\tau) \|_{L^2(\mathbbm{R}_+)}^2+\| P_{xxx}(\tau) \|_{L^2(\mathbbm{R}_+)}^2+\| P_{xx}(\tau) \|_{L^2(\mathbbm{R}_+)}^2)   \,\mathrm{d}\tau\\
    \le& C(\| V_0 \|_{H^4(\mathbbm{R}_+)}^2+\delta^\frac{1}{2}  ).
  \end{aligned}
\end{equation*}
Using the relation \eqref{eq-rela-Vxx-P} again, we can deduce that
\begin{equation*}
  \begin{aligned}[b]
    \int_0^t \| V_{xx}(\tau) \|_{L^2(\mathbbm{R}_+)}^2   \,\mathrm{d}\tau\le C(\| V_0 \|_{H^4(\mathbbm{R}_+)}^2+\delta^\frac{1}{2}  ).
  \end{aligned}
\end{equation*}
On the other hand, \eqref{eq-P-Vrela} with $k=1$ and $l=0$ gives
$$\| P_{xxx} \|_{L^2(\mathbbm{R}_+)}^2\le C(\| P_{xx} \|_{L^2(\mathbbm{R}_+)}^2+\| V_{xx} \|_{L^2(\mathbbm{R}_+)}^2+\| R_{2x} \|_{L^2(\mathbbm{R}_+)}^2), $$
which completes the proof of \eqref{eq-V-H2-jg}.
 $ \hfill\Box $

\begin{lem}\label{lem-VH2-jg} Under the same assumptions as Proposition \ref{prop-V-priori}, there exists a positive constant $C$ such that
  \begin{equation}\label{eq-Vt-cjf-jg}
    \| V_t(t) \|_{H^1(\mathbbm{R}_+)}^2+\int_0^t (\| V_t(\tau) \|_{H^1(\mathbbm{R}_+)}^2+\| P_t(\tau) \|_{H^2(\mathbbm{R}_+)}^2)  \,\mathrm{d}\tau\le C(\| V_0 \|_{H^4(\mathbbm{R}_+)}^2+\delta^\frac{1}{2}).
  \end{equation}
\end{lem}
 {\it\bfseries Proof.}
The equation $\eqref{eq-rd-VP}_1$ gives
\begin{equation}\label{eq-Vt-cjf}
\begin{aligned}[b]
  &\int_{\mathbbm{R}_+}V_t^2  \,\mathrm{d}x \le C\int_{\mathbbm{R}_+}( V_x^2+\tilde{u}_x^2V^2+P_x^2+R_1^2) \,\mathrm{d}x,\\
  &\int_{\mathbbm{R}_+}V_{xt}^2  \,\mathrm{d}x \le C\int_{\mathbbm{R}_+}( V_{xx}^2+V_x^4+\tilde{u}_x^2V_x^2+\tilde{u}_{xx}^2V^2+\tilde{u}_x^4V^2+P_{xx}^2+R_{1x}^2) \,\mathrm{d}x.
\end{aligned}
\end{equation}
It follows that
\begin{equation*}
\begin{aligned}
	&\| V_t(t) \|_{L^2(\mathbbm{R}_+)}^2\le \| V(t) \|_{H^1(\mathbbm{R}_+)}^2+\| P_x(t) \|_{L^2(\mathbbm{R}_+)}^2+C \delta\le  C(\| V_0 \|_{H^4(\mathbbm{R}_+)}^2+\delta^\frac{1}{2}  ),\\
	&\| V_{xt}(t) \|_{L^2(\mathbbm{R}_+)}^2\le \| V(t) \|_{H^2(\mathbbm{R}_+)}^2+\| P_{xx}(t) \|_{L^2(\mathbbm{R}_+)}^2+C \delta\le  C(\| V_0 \|_{H^4(\mathbbm{R}_+)}^2+\delta^\frac{1}{2}  ),	
\end{aligned}
\end{equation*}
and
\begin{equation*}
\begin{aligned}
  \int_0^t \| V_t(\tau) \|_{H^1(\mathbbm{R}_+)}^2  \,\mathrm{d}\tau
  \le& C \int_0^t (\| V_x(\tau) \|_{H^1(\mathbbm{R}_+)}^2+\| \sqrt{\tilde{u}_x}V(\tau) \|_{L^2(\mathbbm{R}_+)}^2+\| P_{x}(\tau) \|_{H^1(\mathbbm{R}_+)}^2)  \,\mathrm{d}\tau+ C \delta\\
  \le&C(\| V_0 \|_{H^4(\mathbbm{R}_+)}^2+\delta^\frac{1}{2}  ).
\end{aligned}
\end{equation*}
From \eqref{eq-rd-VP}$_2$ and $P_{xt}(0,t)=0$, it is easy to obtain
\begin{equation}\label{eq-PtH2}
	 \| P_t(t) \|_{H^2(\mathbbm{R}_+)}^2 \le C \| V_{xt}(t) \|_{L^2(\mathbbm{R}_+)}^2  +C \| R_{2t}(t) \|_{L^2(\mathbbm{R}_+)}^2 \le C  \| V_{xt}(t) \|_{L^2(\mathbbm{R}_+)}^2  +C \delta(1+t)^{-\frac{7}{2} }.
\end{equation}
Integrating \eqref{eq-PtH2} over $[0,t]$, we can obtain \eqref{eq-Vt-cjf-jg}. Therefore, we complete the proof of Lemma \ref{lem-VH2-jg}.
 $ \hfill\Box$

\begin{lem} Under the same assumptions as Proposition \ref{prop-V-priori}, there exists a positive constant $C$ such that
  \begin{equation}\label{eq-Vxxt-H1-jg}
  \begin{aligned}
       \| V_{xxt}(t) \|_{L^2(\mathbbm{R}_+)}^2&+\int_0^t (\| P_{xxxt}(\tau) \|_{L^2(\mathbbm{R}_+)}^2+\| P_{xxt}(\tau) \|_{L^2(\mathbbm{R}_+)}^2+\| V_{xxt}(\tau) \|_{L^2(\mathbbm{R}_+)}^2)   \,\mathrm{d}\tau\\
      &\le C(\| V_0 \|_{H^4(\mathbbm{R}_+)}^2+\delta^\frac{1}{2}  )+C(\| V_0 \|_{H^4(\mathbbm{R}_+)}^2+\delta^\frac{1}{2}  )\int_0^t \| V_{xxx} (\tau)\|_{L^2(\mathbbm{R}_+)}^2  \,\mathrm{d}\tau.
  \end{aligned}
  \end{equation}
\end{lem}
 {\it\bfseries Proof.}
We can get from $\partial_{xxt}\eqref{eq-rd-VP}_1\times V_{xxt}-\partial_{xt}\eqref{eq-rd-VP}_2\times P_{xxxt}$ that
\begin{equation}\label{eq-Vxxt-H1-1}
  \begin{aligned}
    \frac{1}{2} \frac{\mathrm{d}}{\mathrm{d}t}V_{xxt}^2+(f(V+\tilde{u})-f(\tilde{u}))_{xxxt}V_{xxt}-\{P_{xxt}P_{xt}\}_x+P_{xxxt}^2+P_{xxt}^2=R_{1xxt}V_{xxt}-R_{2xt}P_{xxxt}.
  \end{aligned}
\end{equation}
Integrating \eqref{eq-Vxxt-H1-1} over $\mathbbm{R}_+$, using Lemma \ref{lem-f-jf} $(\mathrm{A4})$, $\eqref{eq-f-parx}$ and Cauchy inequality, we can get
\begin{equation}\label{eq-Vxxt-H1-2}
\begin{aligned}[b]
   &\frac{1}{2} \frac{\mathrm{d}}{\mathrm{d}t}\| V_{xxt}(t) \|_{L^2(\mathbbm{R}_+)}^2+\frac{5 \alpha}{2} \| \sqrt{\tilde{u}_x}V_{xxt}(t) \|_{L^2(\mathbbm{R}_+)}^2+\| P_{xxxt}(t) \|_{L^2(\mathbbm{R}_+)}^2+\| P_{xxt}(t) \|_{L^2(\mathbbm{R}_+)}^2\\
  \le& CV_{xxt}^2(0,t)  \!+\!(\| V_x(t) \|_{L^\infty(\mathbbm{R}_+)} \!+\!\mu )\| V_{xxt}(t) \|_{L^2(\mathbbm{R}_+)}^2\!+\!C \mu^{-1}(\| V_t \|_{L^\infty(\mathbbm{R}_+)}^2\!+\!\| \tilde{u}_t \|_{L^\infty(\mathbbm{R}_+)}^2  )\| V_{xxx}(t) \|_{L^2(\mathbbm{R}_+)}^2\\
  &+C \mu^{-1}(\| V_x \|_{L^\infty(\mathbbm{R}_+)}^6+\| \tilde{u}_{xx}(t) \|_{L^\infty(\mathbbm{R}_+)}^2 \| V_x \|_{L^\infty(\mathbbm{R}_+)}^2+\| \tilde{u}_{xxx} \|_{L^\infty(\mathbbm{R}_+)}^2    )\| V_t(t) \|_{L^2(\mathbbm{R}_+)}^2\\
  &+C \mu^{-1}(\| V_x \|_{L^\infty(\mathbbm{R}_+)}^4 \| \tilde{u}_t \|_{L^\infty(\mathbbm{R}_+)}^2+\| \tilde{u}_{xt} \|_{L^\infty(\mathbbm{R}_+)}^2\| V_x(t) \|_{L^\infty(\mathbbm{R}_+)}^2+\| \tilde{u}_{xxt} \|_{L^\infty(\mathbbm{R}_+)}^2   )\| V_{x}(t) \|_{L^2(\mathbbm{R}_+)}^2 \\
  &+ C \mu^{-1}(\| V_x \|_{L^\infty(\mathbbm{R}_+)}^4+ \| \tilde{u}_{xx} \|_{L^\infty(\mathbbm{R}_+)}^2)\| V_{xt}(t) \|_{L^2(\mathbbm{R}_+)}^2+\| \tilde{u}_{xxxxt} \|_{L^\infty(\mathbbm{R}_+)}^2 \| V(t) \|_{L^2(\mathbbm{R}_+)}^2\\
   &+C \mu^{-1}(\| V_{xt} \|_{L^\infty(\mathbbm{R}_+)}^2\!+\!\| \tilde{u}_{xt} \|_{L^\infty(\mathbbm{R}_+)}^2+\| V_{t}V_{x} \|_{L^\infty(\mathbbm{R}_+)}+\| \tilde{u}_tV_{x} \|_{L^\infty(\mathbbm{R}_+)}+\| V_{t}\tilde{u}_{x} \|_{L^\infty(\mathbbm{R}_+)} )\| V_{xx}(t) \|_{L^2(\mathbbm{R}_+)}^2\\
   & + \mu \| P_{xxxt}(t) \|_{L^2(\mathbbm{R}_+)}^2+C \mu^{-1} (\| R_{1xxt}(t) \|_{L^2(\mathbbm{R}_+)}^2+\| R_{2xt}(t) \|_{L^2(\mathbbm{R}_+)}^2) .
\end{aligned}
\end{equation}
By Lemma \ref{lem-bj-gj} $(\mathrm{B5})$, the first term on the right-hand side of \eqref{eq-Vxxt-H1-2} is bounded by
\begin{equation}\label{eq-Vxxt-4}
  V_{xxt}^2(0,t)\le C\delta\mathrm{e}^{-c(1+t)}+\frac{1}{4} \|P_{xxxt}(t) \|_{L^2(\mathbbm{R}_+)}^2+C\|P_{xxt}(t) \|_{L^2(\mathbbm{R}_+)}^2 .
\end{equation}
From \eqref{eq-VxP-rela} with $k=1$, $l=1$, the second term on the right-hand side is bounded by
\begin{equation}\label{eq-Vxxt-5}
   \| V_{xxt} \|_{L^2(\mathbbm{R}_+)}^2 \le 3(\| P_{xxxt} \|_{L^2(\mathbbm{R}_+)}^2 +\| P_{xt} \|_{L^2(\mathbbm{R}_+)}^2+\| R_{2xt} \|_{L^2(\mathbbm{R}_+)}^2).
\end{equation}
Substituting \eqref{eq-Vxxt-4}-\eqref{eq-Vxxt-5} into \eqref{eq-Vxxt-H1-2}, and then integrating the resulting inequality over $[0,t]$, choosing small but fixed $\mu$, we get
\begin{equation*}
  \begin{aligned}[b]
    \| V_{xxt}(t) \|_{L^2(\mathbbm{R}_+)}^2&+\int_0^t (\| \sqrt{\tilde{u}_x}V_{xxt}(\tau) \|_{L^2(\mathbbm{R}_+)}^2+\| P_{xxxt}(\tau) \|_{L^2(\mathbbm{R}_+)}^2+\| P_{xxt}(\tau) \|_{L^2(\mathbbm{R}_+)}^2)   \,\mathrm{d}\tau\\
      &\le C\| V_0 \|_{H^4(\mathbbm{R}_+)}^2+C \int_0^t (\| V_t(\tau) \|_{H^1(\mathbbm{R}_+)}^2+\| V_x(\tau) \|_{H^1(\mathbbm{R}_+)}^2+\|P_{xt}(\tau) \|_{H^1(\mathbbm{R}_+)}^2 ) \,\mathrm{d}\tau\\
      &+C (\| V_0 \|_{H^4(\mathbbm{R}_+)}^2+\delta^\frac{1}{2} )\int_0^t   \| V_{xxx}(\tau) \|_{L^2(\mathbbm{R}_+)}^2\,\mathrm{d}\tau+C \delta.
  \end{aligned}
\end{equation*}
By using \eqref{eq-Vxxt-5} again, and combining \eqref{eq-Vx-cjf-jg}, \eqref{eq-V-H2-jg} and \eqref{eq-Vt-cjf-jg}, the desired estimate \eqref{eq-Vxxt-H1-jg} can be proved.
$  \hfill\Box $

\begin{lem}\label{lem-Vxxx-jg} Under the same assumptions as Proposition \ref{prop-V-priori}, there exists a positive constant $C$ such that
\begin{equation}\label{eq-V-H3-jg}
  \begin{aligned}[b]
    \| V_{xxx}(t) \|_{L^2(\mathbbm{R}_+)}^2+\int_0^t \| V_{xxx}(\tau) \|_{L^2(\mathbbm{R}_+)}^2   \,\mathrm{d}\tau
   \le C(\| V_0 \|_{H^4(\mathbbm{R}_+)}^2+\delta^\frac{1}{2}  ).
  \end{aligned}
\end{equation}
\end{lem}
 {\it\bfseries Proof.}
By $f'(u_+)>0$ and using Taylor's expansion, we can get from $\partial_x^2\eqref{eq-rd-VP}_1$ that
\begin{equation}\label{eq-V-H3-1}
  \begin{aligned}
    f'(u_+)V_{xxx}=&-f''(\xi)(V+\tilde{u}-u_+)V_{xxx}-f'''(V+\tilde{u})(V_x+\tilde{u}_x)^3-3f''(V+\tilde{u})(V_x+\tilde{u}_x)(V_{xx}+\tilde{u}_{xx})\\
    &-f'(V+\tilde{u})\tilde{u}_{xxx}+f'''(\tilde{u})\tilde{u}_x^3+3f''(\tilde{u})\tilde{u}_x\tilde{u}_{xx}+f'(\tilde{u})\tilde{u}_{xxx}+R_{1xx}-P_{xxx}-V_{xxt},
  \end{aligned}
\end{equation}
where $\xi$ is between $V+\tilde{u}$ and $u_+$.

Squaring $\eqref{eq-V-H3-1}$ and then integrating the resulting equation over $\mathbbm{R}_+$, we have
\begin{equation}\label{eq-V-H3-2}
  \begin{aligned}[b]
    \int_{\mathbbm{R}_+} (f'(u_+)V_{xxx})^2 \,\mathrm{d}x  \le &C \int_{\mathbbm{R}_+} (V+\tilde{u}-u_+)^2V_{xxx}^2 \,\mathrm{d}x +C \int_{\mathbbm{R}_+} (\tilde{u}_x^6V^2+\tilde{u}_x^2 \tilde{u}_{xx}^2V^2+\tilde{u}_{xxx}^2V^2) \,\mathrm{d}x \\
    &+C\int_{\mathbbm{R}_+} (V_x^6+\tilde{u}_x^2V_x^4+\tilde{u}_x^4V_x^2+V_x^2V_{xx}^2+\tilde{u}_{xx}^2V_x^2+\tilde{u}_{x}^2V_{xx}^2+R_{1xx}^2+P_{xxx}^2+V_{xxt}^2) \,\mathrm{d}x\\
    \le& C(\varepsilon_0+\delta)\int_{\mathbbm{R}_+} V_{xxx}^2 \,\mathrm{d}x+C(\| \tilde{u}_{x} \|_{L^\infty(\mathbbm{R}_+)}^6+\| \tilde{u}_{x}\tilde{u}_{xx} \|_{L^\infty(\mathbbm{R}_+)}^2+\| \tilde{u}_{xxx} \|_{L^\infty(\mathbbm{R}_+)}^2) \| V(t) \|_{L^2(\mathbbm{R}_+)}^2\\
    &+C(\| V_x \|_{L^\infty(\mathbbm{R}_+)}^4+\| \tilde{u}_x \|_{L^\infty(\mathbbm{R}_+)}^4+\| \tilde{u}_{xx} \|_{L^\infty(\mathbbm{R}_+)}^2   )\| V_x(t) \|_{L^2(\mathbbm{R}_+)}^2+C\| V_{xx}(t) \|_{L^2(\mathbbm{R}_+)}^2\\
    &+C\| R_{1xx}(t) \|_{L^2(\mathbbm{R}_+)}^2+C\| P_{xxx}(t) \|_{L^2(\mathbbm{R}_+)}^2 +C\| V_{xxt}(t) \|_{L^2(\mathbbm{R}_+)}^2    .
  \end{aligned}
\end{equation}
Combining \eqref{eq-V-H2-jg} and \eqref{eq-Vxxt-H1-jg}, we can obtain the desired inequality \eqref{eq-V-H3-jg} for some small $\delta$ and $\varepsilon_0$.
 $ \hfill\Box $

 The equation \eqref{eq-P-Vrela} gives $\| P_{xxxx} \|_{L^2(\mathbbm{R}_+)}^2\le C(\| P_{xx} \|_{L^2(\mathbbm{R}_+)}^2+\| V_{xxx} \|_{L^2(\mathbbm{R}_+)}^2+\| R_{2xx} \|_{L^2(\mathbbm{R}_+)}^2   ) $. Consequently, we have the following Corollary \ref{cor-Pxxxx}.
\begin{cor}\label{cor-Pxxxx} Under the same assumptions as Proposition \ref{prop-V-priori}, there exists a positive constant $C$ such that
\begin{equation}\label{eq-Pxxxx-jg}
  \begin{aligned}[b]
    \| P_{xxxx}(t) \|_{L^2(\mathbbm{R}_+)}^2+\int_0^t \| P_{xxxx}(\tau) \|_{L^2(\mathbbm{R}_+)}^2   \,\mathrm{d}\tau
    \le C(\| V_0 \|_{H^4(\mathbbm{R}_+)}^2+\delta^\frac{1}{2}  ).
  \end{aligned}
\end{equation}
\end{cor}

\begin{lem} Under the same assumptions as Proposition \ref{prop-V-priori}, there exists a positive constant $C$ such that
  \begin{equation}\label{eq-Vtt-L2-jg}
    \begin{aligned}
      \int_0^t (\| V_{tt}(\tau) \|_{H^1(\mathbbm{R}_+)}^2+\| P_{tt}(\tau) \|_{H^2(\mathbbm{R}_+)}^2)  \,\mathrm{d}\tau\le C(\| V_0 \|_{H^4(\mathbbm{R}_+)}^2 +\delta^{\frac{1}{2} }).
    \end{aligned}
  \end{equation}
\end{lem}
 {\it\bfseries Proof.}
We can get from $\partial_t\eqref{eq-rd-VP}$ and $\partial_t\partial_x\eqref{eq-rd-VP}$ that
\begin{equation*}
  \begin{aligned}[b]
    \int_{\mathbbm{R}_+} V_{tt}^2 \,\mathrm{d}x
    \le& C \int_{\mathbbm{R}} (V_x^2V_t^2+\tilde{u}_t^2 V_x^2+\tilde{u}_x^2V_t^2+\tilde{u}_x^2 \tilde{u}_t^2V^2+V_{xt}^2+\tilde{u}_{xt}^2V^2+P_{xt}^2+R_{1t}^2) \,\mathrm{d}x \\
    \le& C(\| V_t \|_{L^\infty(\mathbbm{R}_+)}^2+\| \tilde{u}_t \|_{L^\infty(\mathbbm{R}_+)}^2 )\| V_x \|_{L^2(\mathbbm{R}_+)}^2+C \| \tilde{u}_x \|_{L^\infty}^2 \| V_t \|_{L^2(\mathbbm{R}_+)}^2+C \| V_{xt} \|_{L^2(\mathbbm{R}_+)}^2\\
    &+C  (\| \tilde{u}_{x}\tilde{u}_{t} \|_{L^\infty(\mathbbm{R}_+)}^2+\| \tilde{u}_{xt} \|_{L^\infty(\mathbbm{R}_+)}^2)\| V \|_{L^2(\mathbbm{R}_+)}^2+C\| P_{xt} \|_{L^2(\mathbbm{R}_+)}^2+C \| R_{1t} \|_{L^2(\mathbbm{R}_+)}^2,
\end{aligned}
\end{equation*}
and
\begin{equation}
	\begin{aligned}[b]
		\int_{\mathbbm{R}_+} V_{xtt}^2 \,\mathrm{d}x
        \le& C(\| V_x \|_{L^\infty(\mathbbm{R}_+)}^2+\| \tilde{u}_xV_x \|_{L^\infty(\mathbbm{R}_+)}^2+\| \tilde{u}_x \|_{L^\infty(\mathbbm{R}_+)}^4+\| \tilde{u}_{xx} \|_{L^\infty(\mathbbm{R}_+)}^2    )\| V_t(t) \|_{L^2(\mathbbm{R}_+)}^2\\
        &+C(\| \tilde{u}_tV_x \|_{L^\infty(\mathbbm{R}_+)}^2+\| \tilde{u}_x \tilde{u}_t \|_{L^\infty(\mathbbm{R}_+)}^2  )\| V_x(t) \|_{L^2(\mathbbm{R}_+)}^2+C\| V_{xt}(t) \|_{H^1(\mathbbm{R}_+)}^2+C\| V_{xx}(t) \|_{L^2(\mathbbm{R}_+)}^2  \\
        &+C(\| \tilde{u}_{t}\tilde{u}_{xx} \|_{L^\infty(\mathbbm{R}_+)}^2+\| \tilde{u}_{xxt} \|_{L^\infty(\mathbbm{R}_+)}^2+\| \tilde{u}_x \tilde{u}_{xt} \|_{L^\infty(\mathbbm{R}_+)}^2+\| \tilde{u}_t \tilde{u}_{x}^2 \|_{L^\infty(\mathbbm{R}_+)}^2  ) \| V(t) \|_{L^2(\mathbbm{R}_+)}^2\\
        &+  C\| P_{xxt}(t) \|_{L^2(\mathbbm{R}_+)}^2+C\| R_{1xt}(t) \|_{L^2(\mathbbm{R}_+)}^2.
	\end{aligned}
\end{equation}
On the other hand, we can get from \eqref{eq-rd-VP}$_2$ and $P_{xtt}(0,t)=0$ that
\begin{equation*}
	\| P_{tt}(t) \|_{H^2(\mathbbm{R}_+)}^2\le C(\| V_{xtt}(t) \|_{L^2(\mathbbm{R}_+)}^2+\| V_{tt}(t) \|_{L^2(\mathbbm{R}_+)}^2 ).
\end{equation*}
 Integrating the above three inequalities over $[0,t]$, and using Corollary \ref{cor-R1R2}, the desired estimate \eqref{eq-Vtt-L2-jg} can be obtained.
$  \hfill\Box $

\begin{lem}\label{lem-Vxxtt-jg} Under the same assumptions as Proposition \ref{prop-V-priori}, there exists a positive constant $C$ such that
  \begin{equation}\label{eq-Vxxtt-jg}
    \begin{aligned}[b]
      &\| V_{xxtt}(t) \|_{L^2(\mathbbm{R}_+)}^2+\int_0^t (\| V_{xxtt}(\tau) \|_{L^2(\mathbbm{R}_+)}^2+\| P_{xxxtt}(\tau) \|_{L^2(\mathbbm{R}_+)}^2+\| P_{xxtt}(\tau) \|_{L^2(\mathbbm{R}_+)}^2)  \,\mathrm{d}\tau\\
      \le& C(\| V_0 \|_{H^4(\mathbbm{R}_+)}^2 +\delta^{\frac{1}{2} })+C(\| V_0 \|_{H^4(\mathbbm{R}_+)}^2 +\delta^2)\int_0^t \| V_{xxxt}(\tau) \|_{L^2(\mathbbm{R}_+)}^2  \,\mathrm{d}\tau .
    \end{aligned}
  \end{equation}
\end{lem}
 {\it\bfseries Proof.}
We can get from $\partial_x^2\partial_t^2\eqref{eq-rd-VP}_1\times V_{xxtt}-\partial_x\partial_t^2\eqref{eq-rd-VP}_2\times P_{xxxtt}$ that
\begin{equation}\label{eq-Vxxtt-1}
  \frac{1}{2} \frac{\mathrm{d}}{\mathrm{d}t}V_{xxtt}^2+(f(V+\tilde{u})-f(\tilde{u}))_{xxxtt}V_{xxtt}+P_{xxxtt}^2+P_{xxtt}^2-\{P_{xxtt}P_{xtt}\}_x=R_{1xxtt}V_{xxtt}-P_{xxxtt}R_{2xtt}.
\end{equation}
Integrating \eqref{eq-Vxxtt-1} over $\mathbbm{R}_+$, using Lemma \ref{lem-f-jf} $(\mathrm{A5})$ and Lemma \ref{lem-bj-gj} $(\mathrm{B7})$, we have
\begin{equation*}
  \begin{aligned}[b]
    &\frac{1}{2} \frac{\mathrm{d}}{\mathrm{d}t}\| V_{xxtt}(t) \|_{L^2(\mathbbm{R}_+)}^2+\frac{5 \alpha}{2} \|\sqrt{\tilde{u}_x} V_{xxtt}(t) \|_{L^2(\mathbbm{R}_+)}^2 +\|P_{xxxtt}(t)  \|_{L^2(\mathbbm{R}_+)}^2+ \|P_{xxtt}(t)  \|_{L^2(\mathbbm{R}_+)}^2\\
    \le&C V_{xxtt}^2(0,t)+(\| V_x \|_{L^\infty(\mathbbm{R}_+)}+\mu )\| V_{xxtt}(t) \|_{L^2(\mathbbm{R}_+)}^2+\mu \| P_{xxxtt}(t) \|_{L^2(\mathbbm{R}_+)}^2\\
    &+C\mu^{-1} \| \tilde{u}_{xxxtt} \|_{L^\infty(\mathbbm{R}_+)}^2 \| V(t) \|_{L^2(\mathbbm{R}_+)}^2+C\mu^{-1}(\| V_x(t) \|_{H^2(\mathbbm{R}_+)}^2+\| V_{tt} \|_{H^1(\mathbbm{R}_+)}^2+\| V_{t}(t) \|_{H^2(\mathbbm{R}_+)}^2)\\
    &+C\mu^{-1}(\| V_t \|_{L^\infty(\mathbbm{R}_+)}^2+\| \tilde{u}_t \|_{L^\infty(\mathbbm{R}_+)}^2 )\| V_{xxxt} \|_{L^2(\mathbbm{R}_+)}^2 +C\mu^{-1} (\| R_{1xxtt}(t) \|_{L^2(\mathbbm{R}_+)}^2+\| R_{2xtt}(t) \|_{L^2(\mathbbm{R}_+)}^2)\\
    \le& C(\varepsilon_0+\mu )\| V_{xxtt}(t) \|_{L^2(\mathbbm{R}_+)}^2+2\mu\| P_{xxxtt}(t) \|_{L^2(\mathbbm{R}_+)}^2+C \delta^2(1+t)^{-5}+C \delta \mathrm{e}^{-c(1+t)}+C\| V_x(t) \|_{H^2(\mathbbm{R}_+)}^2\\
    &+C\| V_{tt}(t) \|_{H^1(\mathbbm{R}_+)}^2+C\| V_{t}(t) \|_{H^2(\mathbbm{R}_+)}^2+C(\| V_t \|_{L^\infty(\mathbbm{R}_+)}^2+\| \tilde{u}_t \|_{L^\infty(\mathbbm{R}_+)}^2 )\| V_{xxxt}(t) \|_{L^2(\mathbbm{R}_+)}^2+C\| P_{xxtt}(t) \|_{L^2(\mathbbm{R}_+)}^2.
  \end{aligned}
\end{equation*}
From \eqref{eq-VxP-rela} with $k=1$ and $l=2$, the first term on the right-hand side can be treated as
\begin{equation}\label{eq-VxxttP}
  \| V_{xxtt}(t) \|_{L^2(\mathbbm{R}_+)}^2\le 3(\| P_{xxxtt}(t) \|_{L^2(\mathbbm{R}_+)}^2+\| P_{xtt}(t) \|_{L^2(\mathbbm{R}_+)}^2+\| R_{2xtt} \|_{L^2(\mathbbm{R}_+)}^2),
\end{equation}
it follows that
\begin{equation}\label{eq-Vxxtt-cjf-jg}
  \begin{aligned}[b]
    &\frac{\mathrm{d}}{\mathrm{d}t}\| V_{xxtt}(t) \|_{L^2(\mathbbm{R}_+)}^2+ \|\sqrt{\tilde{u}_x} V_{xxtt}(t) \|_{L^2(\mathbbm{R}_+)}^2 +\|P_{xxxtt}(t)  \|_{L^2(\mathbbm{R}_+)}^2+ \|P_{xxtt}(t)  \|_{L^2(\mathbbm{R}_+)}^2\\
    \le& C(\| P_{xtt}(t) \|_{H^1(\mathbbm{R}_+)}^2+\| V_x(t) \|_{H^2(\mathbbm{R}_+)}^2+\| V_{tt}(t) \|_{H^1(\mathbbm{R}_+)}^2+\| V_{t}(t) \|_{H^2(\mathbbm{R}_+)}^2)\\
    &+C(\| V_t \|_{L^\infty(\mathbbm{R}_+)}^2+\| \tilde{u}_t \|_{L^\infty(\mathbbm{R}_+)}^2 )\| V_{xxxt} \|_{L^2(\mathbbm{R}_+)}^2+C \delta^2(1+t)^{-5}+C \delta \mathrm{e}^{-c(1+t)}.
  \end{aligned}
\end{equation}
Integrating $\eqref{eq-Vxxtt-cjf-jg}$ over $[0,t]$ and using \eqref{eq-VxxttP} again, the estimate $\eqref{eq-Vxxtt-jg}$ can be proved.
 $ \hfill\Box $
\begin{lem}\label{lem-Vxxxt-jg} Under the same assumptions as Proposition \ref{prop-V-priori}, there exists a positive constant $C$ such that
  \begin{equation}\label{eq-Vxxxt-jg}
    \begin{aligned}[b]
      \| V_{xxxt}(t) \|_{L^2(\mathbbm{R}_+)}^2+\int_0^t \| V_{xxxt}(\tau) \|_{L^2(\mathbbm{R}_+)}^2  \,\mathrm{d}\tau
      \le C(\| V_0 \|_{H^4(\mathbbm{R}_+)}^2 +\delta^{\frac{1}{2} }).
    \end{aligned}
  \end{equation}
\end{lem}
 {\it\bfseries Proof.}
By applying Taylor's expansion, we can get from $\partial_x^2\partial_t\eqref{eq-rd-VP}_1$ that
\begin{equation}\label{eq-Vxxxt-1}
  \begin{aligned}[b]
    f'(u_+)V_{xxxt}=&-f''(\xi)(V+\tilde{u}-u_+)V_{xxxt}-f^{(4)}(V+\tilde{u})(V_t+\tilde{u}_t)(V_x+\tilde{u}_x)^3\\
    &-3f^{(3)}(V+\tilde{u})(V_x+\tilde{u}_x)^2(V_{xt}+\tilde{u}_{xt})-3f^{(3)}(V+\tilde{u})(V_t+\tilde{u}_t)(V_x+\tilde{u}_x)(V_{xx}+\tilde{u}_{xx})\\
    &-3f''(V+\tilde{u})(V_{xt}+\tilde{u}_{xt})(V_{xx}+\tilde{u}_{xx})-3f''(V+\tilde{u})(V_x+\tilde{u}_x)(V_{xxt}+\tilde{u}_{xxt})\\
    &-f''(V+\tilde{u})(V_t+\tilde{u}_t)(V_{xxx}+\tilde{u}_{xxx})-f'(V+\tilde{u})\tilde{u}_{xxxt}-f^{(4)}(\tilde{u})\tilde{u}_t\tilde{u}_x^3-3f^{(3)}(\tilde{u})\tilde{u}_x^2\tilde{u}_{xt}\\
    &-3f^{(3)}\tilde{u}_t\tilde{u}_x\tilde{u}_{xx}-3f''(\tilde{u})\tilde{u}_{xt}\tilde{u}_{xx}-3f''(\tilde{u})\tilde{u}_x\tilde{u}_{xxt}-f''(\tilde{u})\tilde{u}_t\tilde{u}_{xxx}-f'(\tilde{u})\tilde{u}_{xxxt}\\
    &-V_{xxtt}-P_{xxxt}+R_{1xxt},
  \end{aligned}
\end{equation}
where $\xi$ is between $V+\tilde{u}$ and $u_+$.

Square $\eqref{eq-Vxxxt-1}$ and then integrate the resulting equation over $\mathbbm{R}_+$. Consequently, we choose small $\varepsilon_0$ and $\delta$ such that the first term $\int_{\mathbbm{R}_+} (V+\tilde{u}-u_+)^2V_{xxxt}^2 \,\mathrm{d}x \le C(\varepsilon_0+\delta)^2\int_{\mathbbm{R}_+} V_{xxxt}^2 \,\mathrm{d}x  $ on the right-hand side of \eqref{eq-Vxxxt-1}  can be absorbed in the left-hand side of \eqref{eq-Vxxxt-1}. Then, we get
\begin{equation}
	\begin{aligned}[b]
		\int_{\mathbbm{R}_+} V_{xxxt}^2 \,\mathrm{d}x \le& C (\| V_x \|_{L^\infty(\mathbbm{R}_+)}^4+\| \tilde{u}_x \|_{L^\infty(\mathbbm{R}_+)}^4+\| \tilde{u}_x \tilde{u}_{xx} \|_{L^\infty(\mathbbm{R}_+)}^2+\| \tilde{u}_{xxx} \|_{L^\infty(\mathbbm{R}_+)}^2   ) \| V_t(t) \|_{L^2(\mathbbm{R}_+)}^2\\
        &+C(\| \tilde{u}_tV_x^2 \|_{L^\infty(\mathbbm{R}_+)}^2+\| \tilde{u}_{xt}V_x \|_{L^\infty(\mathbbm{R}_+)}^2+\| \tilde{u}_{xx}V_t \|_{L^\infty(\mathbbm{R}_+)}^2+\| \tilde{u}_{xxt} \|_{L^\infty(\mathbbm{R}_+)}^2)\| V_x(t) \|_{L^2(\mathbbm{R}_+)}^2\\[1mm]
        &+C (\| \tilde{u}_{t}\tilde{u}_{x}^3 \|_{L^\infty(\mathbbm{R}_+)}^2+\| \tilde{u}_{xt}\tilde{u}_{x}^2 \|_{L^\infty(\mathbbm{R}_+)}^2+\| \tilde{u}_{t}\tilde{u}_{x}\tilde{u}_{xx} \|_{L^\infty(\mathbbm{R}_+)}^2+\| \tilde{u}_{xxxt} \|_{L^\infty(\mathbbm{R}_+)}^2)\| V(t) \|_{L^2(\mathbbm{R}_+)}^2\\[1mm]
        &+C(\| \tilde{u}_{xx}\tilde{u}_{xt} \|_{L^\infty(\mathbbm{R}_+)}^2\!+\!\| \tilde{u}_{xxt}\tilde{u}_{x} \|_{L^\infty(\mathbbm{R}_+)}^2\!+\!\| \tilde{u}_{xxx}\tilde{u}_{t} \|_{L^\infty(\mathbbm{R}_+)}^2)\| V(t) \|_{L^2(\mathbbm{R}_+)}^2\!+\!C \| V_{xt}(t) \|_{H^1(\mathbbm{R}_+)}^2\\[1mm]
        &+C (\| V_{xx}(t) \|_{H^1(\mathbbm{R}_+)}^2+\| V_{xxtt}(t) \|_{L^2(\mathbbm{R}_+)}^2+\| P_{xxtt}(t) \|_{L^2(\mathbbm{R}_+)}^2+\| R_{1xxt}(t) \|_{L^2(\mathbbm{R}_+)}^2 ).
	\end{aligned}
\end{equation}
Combining Lemmas \ref{lem-V-xsb}-\ref{lem-Vxxtt-jg}, we can easily obtain the desired estimate \eqref{eq-Vxxxt-jg}.
  $\hfill\Box $

\begin{lem}\label{lem-Vxxxx-jg}  Under the same assumptions as Proposition \ref{prop-V-priori}, there exists a positive constant $C$ such that
  \begin{equation}\label{eq-Vxxxx-jg}
    \begin{aligned}[b]
      \| \partial_x^4V(t) \|_{L^2(\mathbbm{R}_+)}^2+\| \partial_x^5 P(t) \|_{L^2(\mathbbm{R}_+)}^2+\int_0^t (\| \partial_x^4 V(\tau) \|_{L^2(\mathbbm{R}_+)}^2+\| \partial_x^4 P(\tau) \|_{H^1(\mathbbm{R}_+)}^2)  \,\mathrm{d}\tau
      \le C(\| V_0 \|_{H^4(\mathbbm{R}_+)}^2 +\delta^{\frac{1}{2} }).
    \end{aligned}
  \end{equation}
\end{lem}
 {\it\bfseries Proof.}
Similar to \eqref{eq-Vxxxt-1}, we can get from $\partial_x^3\eqref{eq-rd-VP}_1$ that
\begin{equation}\label{eq-Vxxxx-1}
  \begin{aligned}[b]
    f'(u_+)\partial_x^4V=&-f''(\xi)(V+\tilde{u}-u_+)\partial_x^4 V-f'(V+\tilde{u})\partial_x^4\tilde{u}-f^{(4)}(V+\tilde{u})(V_x+\tilde{u}_x)^4\\
    &-6f^{(3)}(V+\tilde{u})(V_x+\tilde{u}_x)^2(V_{xx}+\tilde{u}_{xx})-3f''(V+\tilde{u})(V_{xx}+\tilde{u}_{xx})^2\\
    &-4f''(V+\tilde{u})(V_x+\tilde{u}_x)(V_{xxx}+\tilde{u}_{xxx})-f^{(4)}(\tilde{u})\tilde{u}_x^4-6f^{(3)}(\tilde{u})\tilde{u}_x^2\tilde{u}_{xx}\\
    &-3f''(\tilde{u})\tilde{u}_{xx}^2-4f''(\tilde{u})\tilde{u}_x\tilde{u}_{xxx}-f'(\tilde{u})\tilde{u}_{xxxx}-V_{xxxt}-\partial_x^4P+R_{1xxx},
  \end{aligned}
\end{equation}
where $\xi$ is between $V+\tilde{u}$ and $u_+$.

Squaring $\eqref{eq-Vxxxx-1}$ and then integrating the resulting equation over $\mathbbm{R}_+$, we obtain
\begin{equation}
	\begin{aligned}[b]
		\int_{\mathbbm{R}_+} (\partial_x^4V)^2 \,\mathrm{d}x
		\le& C(\| V_x \|_{L^\infty(\mathbbm{R}_+)}^6\!+\!\| \tilde{u}_x \|_{L^\infty(\mathbbm{R}_+)}^6\!+\!\| \tilde{u}_{xxx} \|_{L^\infty(\mathbbm{R}_+)}^2\!+\!\| \tilde{u}_x \tilde{u}_{xx} \|_{L^\infty(\mathbbm{R}_+)}^2\!+\!\| \tilde{u}_{xx}V_x \|_{L^\infty(\mathbbm{R}_+)}^2 )\| V_x(t) \|_{L^2(\mathbbm{R}_+)}^2\\
		&+C (\| \partial_x^4 \tilde{u} \|_{L^\infty(\mathbbm{R}_+)}^2\!+\!\| \tilde{u}_x \|_{L^\infty(\mathbbm{R}_+)}^8\!+\!\| \tilde{u}_x \tilde{u}_{xxx} \|_{L^\infty(\mathbbm{R}_+)}^2\!+\!\| \tilde{u}_{xx} \|_{L^\infty(\mathbbm{R}_+)}^4\!+\!\| \tilde{u}_x^2 \tilde{u}_{xx} \|_{L^\infty(\mathbbm{R}_+)}^2 )\| V(t) \|_{L^2(\mathbbm{R}_+)}^2\\
        &+C(\| V_{xx}(t) \|_{H^1(\mathbbm{R}_+)}^2+\| \partial_x^4P(t) \|_{L^2(\mathbbm{R}_+)}^2+\| R_{1xxx}(t) \|_{L^2(\mathbbm{R}_+)}^2 ).
	\end{aligned}
\end{equation}
On the other hand, utilizing $\eqref{eq-P-Vrela}$ with $k=3$ and $l=0$, the fifth derivative of $P$ can be estimated as
\begin{equation*}
 	\| \partial_x^5 P(t) \|_{L^2(\mathbbm{R}_+)}^2\le C(\| P_{xxx}(t) \|_{L^2(\mathbbm{R}_+)}^2+\| \partial_x^4 V(t) \|_{L^2(\mathbbm{R}_+)}^2+\| R_{2xxx}(t) \|_{L^2(\mathbbm{R}_+)}^2).
 \end{equation*}
By combining Lemmas \ref{lem-Vx-jg}, \ref{lem-Vxx-jg}, \ref{lem-Vxxx-jg}, \ref{lem-Vxxxt-jg} and Corollaries \ref{cor-Vx}, \ref{cor-Pxxxx}, the desired estimate \eqref{eq-Vxxxx-jg} can be obtained.
  $\hfill\Box $

Combining Lemmas \ref{lem-V-xsb}-\ref{lem-Vxxxx-jg} , we finish the proof of the \emph{a priori} estimates (Proposition \ref{prop-V-priori}). Next, we denvote ourselves to obtain the decay rate of $V$ by employing the $L^1$-estimate.

\subsubsection{Decay estimates}
To give the decay estimates for the perturbation $V$, we further assume that $V_0\in L^1(\mathbbm{R}_+)$. We define $\phi_\mu(x)$ and $\Phi_\mu(x)$ as follows:
\begin{equation}\label{eq-phi-def}
  \begin{aligned}[b]
    &\phi_\mu(x):=(\rho_\mu\ast\mathrm{sgn})(x)=\int_{-\infty}^{+\infty} \rho_\mu(x-y)\mathrm{sgn}(y) \,\mathrm{d}y,\\
     &\Phi_\mu(x):=\int_0^{x} \phi_\mu(y) \,\mathrm{d}{y},
  \end{aligned}
\end{equation}
where ``$\mathrm{sgn}$'' is a unsual signature function defined as
\begin{equation*}
  \mathrm{sgn}(y):=
  \begin{cases}
    -1, &y<0,\\
    0, &y=0,\\
    1, &y>0.\\
  \end{cases}
\end{equation*}
The symbol $\rho_\mu$ denotes the Friedrichs mollifier defined as
\begin{equation*}
  \rho_\mu(x):=\frac{1}{\mu}\rho \left( \frac{x}{\mu}  \right),
\end{equation*}
where $\rho$ is a smooth function which has a compact support and satisfies $\int_{-\infty}^{+\infty} \rho(x) \,\mathrm{d}x=1 $. We here recall the following properties of $\phi_\mu(x)$ and $\Phi_\mu(x)$. The details can see \cite{Hashimoto2009,Ito1996}.

\begin{lem}\label{lem-phi} Suppose that $\phi_\mu(x)$ and $\Phi_\mu(x)$ is defined in \eqref{eq-phi-def}. $\phi_\mu(x)$ and $\Phi_\mu(x)$ satisfy
	\begin{enumerate}
		\item[(1)] $\lim\limits_{\mu\rightarrow 0} \phi_\mu(x)={\rm sgn}(x),\quad x\in\mathbbm{R}$,\\[-3mm]
		\item[(2)] $\lim\limits_{\mu\rightarrow 0} \Phi_\mu(x)=|x|,\quad x\in\mathbbm{R}$,\\[-3mm]
		\item[(3)] $\phi_\mu(0)=0$,\\[-3mm]
		\item[(4)] $\frac{\mathrm{d}}{\mathrm{d}x} \phi_\mu(x)=2\rho_\mu(x)\ge 0,\quad x\in\mathbbm{R}$.
	\end{enumerate}
\end{lem}

By utilizing Lemma \ref{lem-phi}, we can obtain the following $L^1$-estimate.

\begin{lem}\label{lem-VL1-xsb} ($L^1$-estimate) Suppose that $V_0\in L^1(\mathbbm{R}_+)\cap H^3(\mathbbm{R}_+)$, then the solution $V$ of the problem \eqref{eq-rd-VP}-\eqref{eq-VP-boundary} satisfies
\begin{equation}\label{eq-VL1-xsb-jg}
 \| V(t) \|_{L^1(\mathbbm{R}_+)}\le C(\| V_0 \|_{L^1(\mathbbm{R}_+)}+\delta \mathrm{log}(2+t)).
\end{equation}
\end{lem}

 {\it\bfseries Proof.}
We denote $F(x,t)$ by $F(x,t)=-V_x(x,t)+R_2(x,t)$ and then extend the function $F(x,t)$ such that
\begin{equation*}
	\tilde{F}(x,t):=
	\begin{cases}
		F(x,t),~& x\ge0,\\
		F(-x,t),~&x<0.
	\end{cases}
\end{equation*}
Then $P$ in \eqref{eq-rd-VP}$_2$ can be solve as
\begin{equation}\label{eq-P-bds}
	P(x,t)=\frac{1}{2} \int_{\mathbbm{R}} \mathrm{e}^{-|x-y|}F(y,t) \,\mathrm{d}y= \frac{1}{2} \int_{\mathbbm{R}_+} (\mathrm{e}^{-|x-y|}+\mathrm{e}^{-|x+y|})(-V_x(y,t)+R_2(y,t)) \,\mathrm{d}y, \quad x\in\mathbbm{R}_+.
\end{equation}
 Differentiating \eqref{eq-P-bds} with respect to $x$, we get
\begin{equation}\label{eq-Px-bds}
\begin{aligned}[b]
	P_x(x,t)&=\frac{1}{2} \int_{\mathbbm{R}_+} \left(\frac{\partial \mathrm{e}^{-|x-y|}}{\partial x}+\frac{\partial \mathrm{e}^{-|x+y|}}{\partial x}\right)(-V_x(y,t)+R_2(y,t))  \,\mathrm{d}y\\
	&=\frac{1}{2} \int_{\mathbbm{R}_+} \left(\frac{-\partial \mathrm{e}^{-|x-y|}}{\partial y}+\frac{\partial \mathrm{e}^{-|x+y|}}{\partial y}\right)(-V_x(y,t)+R_2(y,t))  \,\mathrm{d}y\\
	&=\frac{1}{2} \int_{\mathbbm{R}_+} \left(\mathrm{e}^{-|x-y|}- \mathrm{e}^{-|x+y|}\right)(-V_{xx}(y,t)+R_{2x}(y,t))  \,\mathrm{d}y\\
	&=V-\frac{1}{2} \int_{\mathbbm{R}_+} \left(\mathrm{e}^{-|x-y|}- \mathrm{e}^{-|x+y|}\right)(V(y,t)-R_{2x}(y,t))  \,\mathrm{d}y.
\end{aligned}
\end{equation}
It is easy to verify that $P_x(0,t)=0$. In the deriving the last equality of \eqref{eq-Px-bds}, we have used the fact that
\begin{equation*}
	\begin{cases}
		V(x,t)=\int_{\mathbbm{R}_+} (\mathrm{e}^{-|x-y|}- \mathrm{e}^{-|x+y|})(-V_{xx}(y,t)+V(y,t)) \,\mathrm{d}y,\qquad x\in\mathbbm{R_+},\\
		V(0,t)=0.
	\end{cases}
\end{equation*}
We define the operator $K$ as
\begin{equation}\label{eq-def-K}
	K(f)(x)=\frac{1}{2} \int_{\mathbbm{R}_+} \left(\mathrm{e}^{-|x-y|}- \mathrm{e}^{-|x+y|}\right)f(y)  \,\mathrm{d}y.
\end{equation}
Then the equation \eqref{eq-Px-bds} can be rewritten as
\begin{equation}\label{eq-Px-K-bds}
	P_x=V-KV+KR_{2x}.
\end{equation}
Substituting \eqref{eq-Px-K-bds} into \eqref{eq-rd-VP}$_1$, we obtain
\begin{equation}\label{eq-VKV}
	V_t+(f(V+\tilde{u})-f(\tilde{u}))_x+V-KV=R_1-KR_{2x}.
\end{equation}
Multipying \eqref{eq-VKV} by $\phi_\mu(V)$ and then integrating the resulting equation over $\mathbbm{R}_+\times[0,t]$, we have
\begin{equation*}
	\begin{aligned}[b]
		\int_{\mathbbm{R}_+} \Phi_\mu(V) \,\mathrm{d}x+\int_0^t\int_{\mathbbm{R}_+} \phi_\mu(V)(f(V+\tilde{u})-f(\tilde{u}))_x \,\mathrm{d}x \mathrm{d}\tau  +\int_0^t\int_{\mathbbm{R}_+} \phi_\mu(V)(V-KV) \,\mathrm{d}x \mathrm{d}\tau\\
		=\int_{\mathbbm{R}_+} \Phi_\mu(V_0) \,\mathrm{d}x+\int_0^t\int_{\mathbbm{R}_+} \phi_\mu(V)(R_1-KR_{2x}) \,\mathrm{d}x \mathrm{d}\tau.
	\end{aligned}
\end{equation*}
Letting $\mu\rightarrow0$, we can obtain that
\begin{equation*}
	\begin{aligned}
		\int_0^t\int_{\mathbbm{R}_+} \phi_\mu(V)(f(V+\tilde{u})-f(\tilde{u}))_x \,\mathrm{d}x \mathrm{d}\tau
		=\int_0^t\int_{\mathbbm{R}_+}\int_{0}^{V} 2\rho_\mu(\eta)(f'(\eta+\tilde{u})-f'(\tilde{u}))\tilde{u}_x \,\mathrm{d}{\eta} \mathrm{d}x \mathrm{d}\tau\ge 0.
	\end{aligned}
\end{equation*}
By using Young's inequality
$$\| f*g \|_{L^r(\mathbbm{R}_+)}\le \| f \|_{L^p(\mathbbm{R}_+)}\| g \|_{L^q(\mathbbm{R}_+)}, $$
where $1\le r,p,q\le\infty$ with $1/r=1/p+1/q-1$ and ``$*$'' denotes the convolution with respect to the space variable $x$, we can deduce that
$$\| Kf \|_{L^1(\mathbbm{R}_+)}\le \| f \|_{L^1(\mathbbm{R}_+)}. $$
It follows that
\begin{equation*}
	\int_{\mathbbm{R}_+} \phi_\mu(V)(V-KV) \,\mathrm{d}x \ge0,
\end{equation*}
and
\begin{equation*}
	\int_{\mathbbm{R}_+} \phi_\mu(V)(R_1-KR_{2x}) \,\mathrm{d}x \le \| R_1 \|_{L^1(\mathbbm{R}_+)}+\| R_{2x} \|_{L^1(\mathbbm{R}_+)}\le C \delta(1+t)^{-1}+C \delta\mathrm{e}^{-c(1+t)},
\end{equation*}
for $\mu \rightarrow 0$. Thus, letting $\mu\rightarrow 0$, we can obtain the desired estimate \eqref{eq-VL1-xsb-jg}.
 $ \hfill\Box$

\begin{prop}\label{prop-Vrate-xsb} (Decay estimates of $V$) Suppose that $f'(u_-)>0$ and $(V(x,t),P(x,t))$ is a solution of the problem \eqref{eq-rd-VP}-\eqref{eq-VP-boundary}, for $\varepsilon\in(0,\frac{1}{2} )$ and sufficiently large $t$, $(V(x,t),P(x,t))$ satisfies
  \begin{align}
    (1+t)^{\frac{1}{2} +\varepsilon}&\int_{\mathbbm{R}_+} (| V|^2+| V_x|^2) \,\mathrm{d}x +\int_0^t (1+\tau)^{\frac{1}{2} +\varepsilon}\int_{\mathbbm{R}_+} \tilde{u}_x(| V|^2+|V_x|^2) \,\mathrm{d}x  \,\mathrm{d}\tau\nonumber\\
    & +\int_0^t (1+\tau)^{\frac{1}{2} +\varepsilon}\int_{\mathbbm{R}_+} |V_x|^2 \,\mathrm{d}x  \,\mathrm{d}\tau\le C (1+t)^\varepsilon \mathrm{log}^2(2+t),\label{eq-rate-V-jg-1}\\
    (1+t)^{\frac{3}{2} +\varepsilon}&\int_{\mathbbm{R}_+} (| V_x|^2+|V_{xx}|^2) \,\mathrm{d}x +\int_0^t (1+\tau)^{\frac{1}{2} +\varepsilon}\int_{\mathbbm{R}_+} \tilde{u}_x(| V_x|^2+|V_{xx}|^2) \,\mathrm{d}x  \,\mathrm{d}\tau\nonumber\\
    & +\int_0^t (1+\tau)^{\frac{1}{2} +\varepsilon}\int_{\mathbbm{R}_+} |V_{xx}|^2 \,\mathrm{d}x  \,\mathrm{d}\tau\le C (1+t)^\varepsilon \mathrm{log}^{10}(2+t),\label{eq-rate-V-jg-2}\\
    (1+t)^{\frac{3}{2} +\varepsilon}&\int_{\mathbbm{R}_+} (| V_{xxx}|^2\!+\!\sum\limits_{j=0}^{2}|\partial_x^j V_t|^2) \,\mathrm{d}x \!+\!\int_0^t (1\!+\!\tau)^{\frac{3}{2} +\varepsilon}\int_{\mathbbm{R}_+} |V_{xxx}|^2 \,\mathrm{d}x  \,\mathrm{d}\tau\le C  (1\!+\!t)^\varepsilon \mathrm{log}^{10}(2\!+\!t),\label{eq-rate-V-jg-3}\\
    (1+t)^{\frac{3}{2} +\varepsilon}&\int_{\mathbbm{R}_+} (|\partial_x^4 V|^2\!+\!\sum\limits_{j=0}^{2}|\partial_x^j V_{tt}|^2\!+\!| V_{xxxt}|^2) \,\mathrm{d}x +\int_0^t (1\!+\!\tau)^{\frac{3}{2} \!+\!\varepsilon}\int_{\mathbbm{R}_+} |\partial_x^4 V|^2 \,\mathrm{d}x  \,\mathrm{d}\tau\le C  (1\!+\!t)^\varepsilon \mathrm{log}^{10}(2\!+\!t),\label{eq-rate-V-jg-4}
  \end{align}
  and
  \begin{equation}\label{eq-P-rate}
  	(1+t)^{\frac{3}{2} +\varepsilon}\int_{\mathbbm{R}_+} |\partial_x^j P|^2 \,\mathrm{d}x\le C  (1+t)^\varepsilon \mathrm{log}^{10}(2+t),\quad j=0,1,2,3,4,5.
  \end{equation}
\end{prop}

 {\it\bfseries Proof.}
Adding $\eqref{eq-V-L2-3}$ and $\eqref{eq-H1-djf0}$, we have
\begin{equation}\label{eq-rate-V-1}
  \begin{aligned}[b]
    &\frac{\mathrm{d}}{\mathrm{d}t}(\| V(t) \|_{L^2(\mathbbm{R}_+)}^2+\| V_x(t) \|_{L^2(\mathbbm{R}_+)}^2 )+\| \sqrt{\tilde{u}_x}V(t) \|_{L^2(\mathbbm{R}_+)}^2+ \| \sqrt{\tilde{u}_x}V_x(t) \|_{L^2(\mathbbm{R}_+)}^2\\
    &\ \ \ \ \ \ \ +\| P_{xx}(t) \|_{L^2(\mathbbm{R}_+)}^2+ 2\| P_{x}(t) \|_{L^2(\mathbbm{R}_+)}^2+\| P(t) \|_{L^2(\mathbbm{R}_+)}^2\\
    \le& C (\varepsilon_0+\mu)\| V_x(t) \|_{L^2(\mathbbm{R}_+)}^2+C (1+t)^{-\frac{3}{2} } +C \mathrm{e}^{-c(1+t)}.
  \end{aligned}
\end{equation}
By $\eqref{eq-rd-VP}_2$, we get
\begin{equation*}
  \int_{\mathbbm{R}_+} (P_{xx}^2+2P_x+P^2)  \,\mathrm{d}x=\int_{\mathbbm{R}_+} (V_x^2+2R_2V_x+R_2^2) \,\mathrm{d}x,
\end{equation*}
which implies that
\begin{equation}\label{eq-rate-V-2}
\begin{aligned}[b]
    \int_{\mathbbm{R}_+} V_x^2 \,\mathrm{d}x
  &\le \int_{\mathbbm{R}_+} (P_{xx}^2+2P_x+P^2)  \,\mathrm{d}x+ C \int_{\mathbbm{R}_+} R_2^2 \,\mathrm{d}x\\
  &\le  \int_{\mathbbm{R}_+} (P_{xx}^2+2P_x+P^2)  \,\mathrm{d}x+ C(1+t)^{-\frac{7}{2} }+C \mathrm{e}^{-c(1+t)}.
\end{aligned}
\end{equation}
Therefore, adding $\eqref{eq-rate-V-1}$ and $\eqref{eq-rate-V-2}$, and then multipying the resulting inequality by $(1+t)^{\frac{1}{2} +\varepsilon}$ and integrating it over $[0,t]$, we have
\begin{equation}\label{eq-rate-V-3}
  \begin{aligned}[b]
    &(1+t)^{\frac{1}{2} +\varepsilon}(\| V(t) \|_{L^2(\mathbbm{R}_+)}^2+\| V_x(t) \|_{L^2(\mathbbm{R}_+)}^2)+\int_0^t (1+\tau)^{\frac{1}{2} +\varepsilon}\| \sqrt{\tilde{u}_x}V(\tau) \|_{L^2(\mathbbm{R}_+)}^2 \,\mathrm{d}\tau \\
    &\ \ \ \ +\int_0^t (1+\tau)^{\frac{1}{2} +\varepsilon}\| \sqrt{\tilde{u}_x}V_x(\tau) \|_{L^2(\mathbbm{R}_+)}^2 \,\mathrm{d}\tau +\int_0^t(1+\tau)^{\frac{1}{2} +\varepsilon} \| V_x(\tau) \|_{L^2(\mathbbm{R}_+)}^2 \,\mathrm{d}\tau \\
    \le& C \| V_0 \|_{H^4(\mathbbm{R}_+)}^2(1+t)^\varepsilon+(\frac{1}{2} +\varepsilon) \int_0^t (1+\tau)^{-\frac{1}{2}+\varepsilon }\| V(\tau) \|_{L^2(\mathbbm{R}_+)}^2  \,\mathrm{d}\tau +(\frac{1}{2} +\varepsilon)\int_0^t (1+\tau)^{-\frac{1}{2} +\varepsilon}\| V_x(\tau) \|_{L^2(\mathbbm{R}_+)}^2  \,\mathrm{d}\tau.
  \end{aligned}
\end{equation}
By employing Gagliardo-Nirenberg inequality
\begin{equation}\label{eq-f-L2}
  \|f\|_{L^2(\mathbbm{R}_+)}^2\le C\|f\|_{L^1(\mathbbm{R}_+)}^{\frac{4}{3}}\| f_x \|_{L^2(\mathbbm{R}_+)}^{\frac{2}{3}},
\end{equation}
the second term on the right-hand side of $\eqref{eq-rate-V-3}$ can be estimated as
\begin{equation}\label{eq-rate-V-4}
  \begin{aligned}[b]
    &(\frac{1}{2} +\varepsilon)\int_0^t (1+\tau)^{-\frac{1}{2}+\varepsilon }\| V(\tau) \|_{L^2(\mathbbm{R}_+)}^2  \,\mathrm{d}\tau\\
    \le& C \int_0^t (1+\tau)^{-1+\varepsilon }\| V(\tau) \|_{L^1(\mathbbm{R}_+)}^2 \,\mathrm{d}\tau+\frac{1}{4} \int_0^t (1+\tau)^{\frac{1}{2} +\varepsilon }\| V_x(\tau) \|_{L^2(\mathbbm{R}_+)}^2 \,\mathrm{d}\tau\\
    \le& C(1+t)^\varepsilon \mathrm{log}^2(2+t)+\frac{1}{4}
    \int_0^t (1+\tau)^{\frac{1}{2} +\varepsilon }\| V_x(\tau) \|_{L^2(\mathbbm{R}_+)}^2 \,\mathrm{d}\tau.
  \end{aligned}
\end{equation}
 We treat the last term on the right-hand side of $\eqref{eq-rate-V-3}$. For any $\varepsilon\in(0,\frac{1}{2} )$, it holds that
 \begin{equation}\label{eq-rate-V-5}
   \begin{aligned}[b]
     (\frac{1}{2} +\varepsilon)\int_0^t (1+\tau)^{-\frac{1}{2} +\varepsilon}\| V_x(\tau) \|_{L^2(\mathbbm{R}_+)}^2  \,\mathrm{d}\tau
     \le C\int_0^t \| V_x(\tau) \|_{L^2(\mathbbm{R}_+)}^2  \,\mathrm{d}\tau
     \le C (\| V_0 \|_{H^4(\mathbbm{R}_+)}^2+\delta^{\frac{1}{2}})  .
   \end{aligned}
 \end{equation}
 Substituting \eqref{eq-rate-V-4}-\eqref{eq-rate-V-5} into \eqref{eq-rate-V-3}, we obtain \eqref{eq-rate-V-jg-1}.\\
\indent Next, we show \eqref{eq-rate-V-jg-2}. From \eqref{eq-V-H1-2} and Lemma \ref{lem-bj-gj} $\mathrm{(B1)}$, we can get
\begin{equation}\label{eq-rete-Vx-1}
  \begin{aligned}[b]
    & \frac{\mathrm{d}}{\mathrm{d}t}\| V_x(t) \|_{L^2(\mathbbm{R}_+)}^2+\| \sqrt{\tilde{u}_x}V_x(t) \|_{L^2(\mathbbm{R}_+)}^2   +\| P_{xx}(t) \|_{L^2(\mathbbm{R}_+)}^2+\| P_x(t) \|_{L^2(\mathbbm{R}_+)}^2 \\
  \le&C \mathrm{e}^{-c(1+t)}+(1+t)^{-1}\int_{\mathbbm{R}_+} V_x^2 \,\mathrm{d}x+C(1+t)\int_{\mathbbm{R}} R_{1x}^2 \,\mathrm{d}x+ \mu \int_{\mathbbm{R}_+} P_{xx}^2 \,\mathrm{d}x +C\mu^{-1} \int_{\mathbbm{R}_+} R_{2}^2  \,\mathrm{d}x\\
     &+C \int_{\mathbbm{R}_+} |V_x|(| \tilde{u}_{xx} ||V|+\tilde{u}_x^2|V|+V_x^2 )  \,\mathrm{d}x.
  \end{aligned}
\end{equation}
The last term can be estimated as follows. From $\eqref{eq-f-parx}$ and \eqref{eq-rate-V-jg-1},
\begin{equation}\label{eq-rete-Vx-2}
  \begin{aligned}[b]
    \int_{\mathbbm{R}_+} |V_x|(| \tilde{u}_{xx} |+\tilde{u}_x^2 )|V|  \,\mathrm{d}x
    \le&\| \tilde{u}_{xx}(t) \|_{L^\infty(\mathbbm{R}_+)}\| V_x(t) \|_{L^2(\mathbbm{R}_+)}\| V(t) \|_{L^2(\mathbbm{R}_+)}\\[-1mm]
    \le& \| \tilde{u}_{xx} \|_{L^\infty}^{\frac{2}{3} }\| V_x(t) \|_{L^2(\mathbbm{R}_+)}^2+ \| \tilde{u}_{xx} \|_{L^\infty}^{\frac{4}{3} }\| V(t) \|_{L^2(\mathbbm{R}_+)}^2 \\[1mm]
    \le& (1+t)^{-1}\| V_x(t) \|_{L^2(\mathbbm{R}_+)}^2+C(1+t)^{-2}\| V(t) \|_{L^2(\mathbbm{R}_+)}^2\\[1mm]
    \le& (1+t)^{-1}\| V_x(t) \|_{L^2(\mathbbm{R}_+)}^2+C(1+t)^{-\frac{5}{2} }\mathrm{log}^{2}(2+t).
  \end{aligned}
\end{equation}
Using Gagliardo-Nirenberg inequality
$$\| f_x \|_{L^3(\mathbbm{R}_+)}^3\le C \| f \|_{L^2(\mathbbm{R}_+)}^{\frac{5}{4}}\| f_{xx} \|_{L^2(\mathbbm{R}_+)}^{\frac{7}{4}},   $$
we have
\begin{equation}\label{eq-rete-Vx-3}
  \begin{aligned}
  \int_{\mathbbm{R}_+} |V_x|^3 \,\mathrm{d}x
  &\le  \mu \| V_{xx}(t) \|_{L^2(\mathbbm{R}_+)}^2  +C \mu^{-1}\| V(t) \|_{L^2(\mathbbm{R}_+)}^{10}\\
  &\le \mu \| V_{xx}(t) \|_{L^2(\mathbbm{R}_+)}^2  +C\mu^{-1}(1+t)^{-\frac{5}{2} }\mathrm{log}^{10}(2+t).
  \end{aligned}
\end{equation}
Substituting \eqref{eq-rete-Vx-2}-\eqref{eq-rete-Vx-3} into \eqref{eq-rete-Vx-1}, for some small but fixed $\mu$, we can deduce that
\begin{equation}\label{eq-rete-Vx-3-1}
  \begin{aligned}[b]
    & \frac{\mathrm{d}}{\mathrm{d}t}\| V_x(t) \|_{L^2(\mathbbm{R}_+)}^2+\| \sqrt{\tilde{u}_x}V_x(t) \|_{L^2(\mathbbm{R}_+)}^2   +\| P_{xx}(t) \|_{L^2(\mathbbm{R}_+)}^2+\| P_x(t) \|_{L^2(\mathbbm{R}_+)}^2 \\
  \le& \mu \| V_{xx}(t) \|_{L^2(\mathbbm{R}_+)}^2+C (1+t)^{-1}\| V_x(t) \|_{L^2(\mathbbm{R}_+)}^2+C(1+t)^{-\frac{5}{2} }\mathrm{log}^{10}(2+t)+ C \mathrm{e}^{-c(1+t)}.
  \end{aligned}
\end{equation}
By \eqref{eq-V-H2-3}, we can get from \eqref{eq-V-H2-2} that
\begin{equation}\label{eq-rete-Vx-4}
  \begin{aligned}[b]
    &\frac{\mathrm{d}}{\mathrm{d}t}\| V_{xx}(t) \|_{L^2(\mathbbm{R}_+)}^2+\| \sqrt{\tilde{u}_x}V_{xx}(t) \|_{L^2(\mathbbm{R}_+)}^2 +\| P_{xxx}(t) \|_{L^2(\mathbbm{R}_+)}^2+\| P_{xx}(t) \|_{L^2(\mathbbm{R}_+)}^2\\
     \le& C (\varepsilon_0+\mu) \| V_{xx}(t) \|_{L^2(\mathbbm{R}_+)}^2+CV_{xx}^2(0,t)+\mu \| P_{xxx}(t) \|_{L^2(\mathbbm{R}_+)}^2+C \mu^{-1}(\| R_{1xx} \|_{L^2(\mathbbm{R}_+)}^2+\| R_{2x} \|_{L^2(\mathbbm{R}_+)}^2 ) \\
     &+C \mu^{-1} \left(\| \tilde{u}_{xxx} \|_{L^\infty(\mathbbm{R}_+)}^2 \| V(t) \|_{L^2(\mathbbm{R}_+)}^2+\| \tilde{u}_{xx} \|_{L^\infty(\mathbbm{R}_+)}^2 \| V_x(t) \|_{L^2(\mathbbm{R}_+)}^2+\| V_x(t) \|_{L^6(\mathbbm{R}_+)}^6 \right)\\
     \le&C (\varepsilon_0+\mu) \| V_{xx}(t) \|_{L^2(\mathbbm{R}_+)}^2+\mu \| P_{xxx}(t) \|_{L^2(\mathbbm{R}_+)}^2+C \mathrm{e}^{-c(1+t)}+C \| P_{xx}(t) \|_{L^2(\mathbbm{R}_+)}^2\\
     &+C (1+t)^{-4 }\| V(t) \|_{L^2(\mathbbm{R}_+)}^2+C (1+t)^{-3 } \| V_x(t) \|_{L^2(\mathbbm{R}_+)}^2+C(1+t)^{-\frac{9}{2} }.
  \end{aligned}
\end{equation}
In the deriving of the last inequality of \eqref{eq-rete-Vx-4}, we have used the fact that
$$\| V_x(t) \|_{L^6(\mathbbm{R}_+)}^6\le C \| V(t) \|_{L^2(\mathbbm{R}_+)}^4 \| V_{xx}(t) \|_{L^2(\mathbbm{R}_+)}^2\le C \varepsilon_0^4\| V_{xx}(t) \|_{L^2(\mathbbm{R}_+)}^2.  $$
Multipy $\eqref{eq-rete-Vx-3-1}$ by a sufficiently large number $\lambda$ such that the fourth term on the right-hand side of $\eqref{eq-rete-Vx-4}$ can be absorbed in the third on the left-hand side of $\eqref{eq-rete-Vx-3-1}$, then add the resulting inequality to $\eqref{eq-rete-Vx-4}$. Consequently, choosing small $\varepsilon_0$ and $\mu$ and using $\eqref{eq-rela-Vxx-P}$, we have
\begin{equation}\label{eq-rete-Vx-5}
  \begin{aligned}[b]
    &\frac{\mathrm{d}}{\mathrm{d}t}(\| V_{x}(t) \|_{L^2(\mathbbm{R}_+)}^2\!+\!\| V_{xx}(t) \|_{L^2(\mathbbm{R}_+)}^2)\!+\!\| \sqrt{\tilde{u}_x}V_{x}(t) \|_{L^2(\mathbbm{R}_+)}^2\!+\!\| \sqrt{\tilde{u}_x}V_{xx}(t) \|_{L^2(\mathbbm{R}_+)}^2 \!+\!\| V_{xx}(t) \|_{L^2(\mathbbm{R}_+)}^2\!+\!\| P_x(t) \|_{H^2(\mathbbm{R}_+)}^2 \\
    &\le C (1+t)^{-1}\| V_x(t) \|_{L^2(\mathbbm{R}_+)}^2+C(1+t)^{-4}\| V(t) \|_{L^2(\mathbbm{R}_+)}^2  +C(1+t)^{-\frac{5}{2} }\mathrm{log}^{10}(2+t)+ C \mathrm{e}^{-c(1+t)}.
  \end{aligned}
\end{equation}
Multipying $\eqref{eq-rete-Vx-5}$ by $(1+t)^{\frac{3}{2}+\varepsilon }$ and then integrating it over $[0,t]$, by $\eqref{eq-V-propjg}$, we can deduce that
\begin{equation*}
  \begin{aligned}[b]
    &(1+t)^{\frac{3}{2}+\varepsilon }\| V_{x}(t) \|_{H^1(\mathbbm{R}_+)}^2+\int_0^t (1+\tau)^{\frac{3}{2}+\varepsilon}(\| \sqrt{\tilde{u}_x}V_{x}(\tau) \|_{L^2(\mathbbm{R}_+)}^2+\| \sqrt{\tilde{u}_x}V_{xx}(\tau) \|_{L^2(\mathbbm{R}_+)}^2) \,\mathrm{d}\tau\\
     &\ \ \ \ +\int_0^t (1+\tau)^{\frac{3}{2}+\varepsilon}(\| V_{xx}(\tau) \|_{L^2(\mathbbm{R}_+)}^2+\| P_{x}(\tau) \|_{H^2(\mathbbm{R}_+)}^2 ) \,\mathrm{d}\tau\\
    \le& C (1+t)^\varepsilon\mathrm{log}^{10}(2+t)\\
    &+(\frac{3}{2}+\varepsilon )\int_0^t (1+\tau)^{\frac{1}{2} +\varepsilon} \| V_{x}(\tau) \|_{L^2(\mathbbm{R}_+)}^2\,\mathrm{d}\tau +(\frac{3}{2}+\varepsilon ) \int_0^t (1+\tau)^{\frac{1}{2} +\varepsilon} \| V_{xx}(\tau) \|_{L^2(\mathbbm{R}_+)}^2\,\mathrm{d}\tau.
  \end{aligned}
\end{equation*}
For sufficiently large $T$ such that $(\frac{3}{2}+\varepsilon )\frac{1}{1+T}\le \frac{1}{2} $, the last term can be estimated as
\begin{equation}\label{eq-Vxx-eps}
  \begin{aligned}[b]
    &(\frac{3}{2}+\varepsilon )\int_0^t (1+\tau)^{\frac{1}{2} +\varepsilon} \| V_{xx}(\tau) \|_{L^2(\mathbbm{R}_+)}^2\,\mathrm{d}\tau\\
    =&(\frac{3}{2}+\varepsilon )\int_0^T (1+\tau)^{\frac{1}{2} +\varepsilon} \| V_{xx}(\tau) \|_{L^2(\mathbbm{R}_+)}^2\,\mathrm{d}\tau+(\frac{3}{2}+\varepsilon )\frac{1}{1+T}\int_T^t (1+\tau)^{\frac{3}{2} +\varepsilon} \| V_{xx}(\tau) \|_{L^2(\mathbbm{R}_+)}^2\,\mathrm{d}\tau\\
    \le& C \int_0^t \| V_{xx}(\tau) \|_{L^2(\mathbbm{R}_+)}^2  \,\mathrm{d}\tau+\frac{1}{2} \int_0^t  (1+\tau)^{\frac{3}{2}+\varepsilon}\| V_{xx}(\tau) \|_{L^2(\mathbbm{R}_+)}^2 \,\mathrm{d}\tau.
  \end{aligned}
\end{equation}
Combining \eqref{eq-V-H2-jg} and \eqref{eq-rate-V-jg-1}, we have
\begin{equation}\label{eq-PxH2-sjl}
  \begin{aligned}[b]
    &(1+t)^{\frac{3}{2}+\varepsilon }\| V_{x}(t) \|_{H^1(\mathbbm{R}_+)}^2+\int_0^t (1+\tau)^{\frac{3}{2}+\varepsilon}(\| \sqrt{\tilde{u}_x}V_{x}(\tau) \|_{L^2(\mathbbm{R}_+)}^2+\| \sqrt{\tilde{u}_x}V_{xx}(\tau) \|_{L^2(\mathbbm{R}_+)}^2\,\mathrm{d}\tau\\
    &\ \ \ \ +\int_0^t (1+\tau)^{\frac{3}{2}+\varepsilon}(\| V_{xx}(\tau) \|_{L^2(\mathbbm{R}_+)}^2+\| P_{x}(\tau) \|_{H^2(\mathbbm{R}_+)}^2 ) \,\mathrm{d}\tau\le C (1+t)^\varepsilon\mathrm{log}^{10}(2+t),
  \end{aligned}
\end{equation}
which completes the proof of $\eqref{eq-rate-V-jg-2}$. Moreover, it follows that
\begin{equation}\label{eq-Vx-infrate}
  \| V_x(t) \|_{L^\infty(\mathbbm{R}_+)}^2\le C(1+t)^{-\frac{3}{2} }\mathrm{log}^{10}(2+t).
\end{equation}
On the other hand, using \eqref{eq-P-xsb-1} and \eqref{eq-P-Vrela} with $k=1$, we also have 
\begin{equation}\label{eq-P-sjljs}
	(1+t)^{\frac{3}{2}+\varepsilon }\| P(t) \|_{H^3(\mathbbm{R}_+)}^2\le C (1+t)^\varepsilon\mathrm{log}^{10}(2+t).
\end{equation}
\indent We now show $\eqref{eq-rate-V-jg-2}$. We can get from $\eqref{eq-Vt-cjf}$ that
\begin{equation*}
\begin{aligned}[b]
    &(1+t)^{\frac{3}{2} +\varepsilon}\| V_t(t) \|_{L^2(\mathbbm{R}_+)}^2+\int_0^t (1+\tau)^{\frac{1}{2} +\varepsilon}\| V_t(\tau) \|_{L^2(\mathbbm{R}_+)}^2 \,\mathrm{d}\tau\le C(1+t)^\varepsilon\mathrm{log}^{10}(2+t),\\
    &(1+t)^{\frac{3}{2} +\varepsilon}\| V_{xt}(t) \|_{L^2(\mathbbm{R}_+)}^2+\int_0^t (1+\tau)^{\frac{3}{2} +\varepsilon}\| V_{xt}(\tau) \|_{L^2(\mathbbm{R}_+)}^2 \,\mathrm{d}\tau\le C(1+t)^\varepsilon\mathrm{log}^{10}(2+t).
\end{aligned}
\end{equation*}
It follows from \eqref{eq-PtH2} that 
\begin{equation}\label{eq-Pt-sjl}
	\int_0^t (1+\tau)^{\frac{3}{2} +\varepsilon}\| P_{t}(\tau) \|_{H^2(\mathbbm{R}_+)}^2 \,\mathrm{d}\tau\le C(1+t)^\varepsilon\mathrm{log}^{10}(2+t).
\end{equation}
By utilizing  \eqref{eq-Vxxt-4}-\eqref{eq-Vxxt-5} and \eqref{eq-Vx-infrate}, we can get from \eqref{eq-Vxxt-H1-2} that
\begin{equation}\label{eq-ratek2-5}
  \begin{aligned}[b]
    & \frac{\mathrm{d}}{\mathrm{d}t}\| V_{xxt}(t) \|_{L^2(\mathbbm{R}_+)}^2+ \| \sqrt{\tilde{u}_x}V_{xxt}(t) \|_{L^2(\mathbbm{R}_+)}^2+\| P_{xxxt}(t) \|_{L^2(\mathbbm{R}_+)}^2+\| P_{xxt}(t) \|_{L^2(\mathbbm{R}_+)}^2+\| V_{xxt}(t) \|_{L^2(\mathbbm{R}_+)}^2\\
    \le& C \mathrm{e}^{-c(1+t)}+C (\varepsilon_0+\mu) \| V_{xxt}(t) \|_{L^2(\mathbbm{R}_+)}^2+C(1+t)^{-4  }\| V_x(t) \|_{L^2(\mathbbm{R}_+)}^2\\
    &+C(1+t)^{-4  }\| V_t(t) \|_{L^2(\mathbbm{R}_+)}^2 +C\| V_{xt}(t) \|_{L^2(\mathbbm{R}_+)}^2+C\| V_{xx}(t) \|_{L^2(\mathbbm{R}_+)}^2+C \| P_{xt}(t) \|_{H^1(\mathbbm{R}_+)}^2 \\
    &+C (\varepsilon_0+\delta)\| V_{xxx}(t) \|_{L^2(\mathbbm{R}_+)}^2+C(1+t)^{-6}\| V(t) \|_{L^2(\mathbbm{R}_+)}^2+C(1+t)^{-\frac{11}{2} }.
  \end{aligned}
\end{equation}
Multipying $\eqref{eq-ratek2-5}$ by $(1+t)^{\frac{3}{2}+\varepsilon }$ and then integrating \eqref{eq-ratek2-5} over $[0,t]$, for some small $\varepsilon_0$ and $\mu$, we obtain
\begin{equation}\label{eq-ratek2-6}
  \begin{aligned}[b]
    & (1+t)^{\frac{3}{2}+\varepsilon }\| V_{xxt}(t) \|_{L^2(\mathbbm{R}_+)}^2+ \int_0^t (1+\tau)^{\frac{3}{2} +\varepsilon}(\| \sqrt{\tilde{u}_x}V_{xxt}(\tau) \|_{L^2(\mathbbm{R}_+)}^2 +\| P_{xxt}(\tau) \|_{H^1(\mathbbm{R}_+)}^2+\| V_{xxt}(\tau) \|_{L^2(\mathbbm{R}_+)}^2) \,\mathrm{d}\tau\\
    \le& C (1+t)^\varepsilon\mathrm{log}^{10}(2+t)+(\frac{3}{2} +\varepsilon)\int_0^t (1+\tau)^{\frac{1}{2} +\varepsilon}\| V_{xxt}(\tau) \|_{L^2(\mathbbm{R}_+)}^2  \,\mathrm{d}\tau\\
    &+C (\varepsilon_0+\delta)\int_0^t (1+\tau)^{\frac{3}{2} +\varepsilon}\| V_{xxx}(\tau) \|_{L^2(\mathbbm{R}_+)}^2  \,\mathrm{d}\tau.
  \end{aligned}
\end{equation}
The second term on the right-hand side of \eqref{eq-ratek2-6} can be estimated by using the similar method in \eqref{eq-Vxx-eps}. To deal with the last term of $\eqref{eq-ratek2-6}$, we can get from $\eqref{eq-V-H3-2}$ that
\begin{equation}\label{eq-ratek2-7}
  \begin{aligned}[b]
     \| V_{xxx}(t) \|_{L^2(\mathbbm{R}_+)}^2
    \le &C (1+t)^{-6}\| V(t) \|_{L^2(\mathbbm{R}_+)}^2+ C (1+t)^{-3  }\mathrm{log}^{20}(2+t)\| V_x(t) \|_{L^2(\mathbbm{R}_+)}^2\\
    &+C \| V_{xx}(t) \|_{L^2(\mathbbm{R}_+)}^2+C \| P_{xxx}(t) \|_{L^2(\mathbbm{R}_+)}^2+C \| V_{xxt} \|_{H^1(\mathbbm{R}_+)}^2+C(1+t)^{-\frac{9}{2} }.
  \end{aligned}
\end{equation}
By combining \eqref{eq-rate-V-jg-2}, \eqref{eq-P-sjljs} and \eqref{eq-Pt-sjl}, it follows from \eqref{eq-ratek2-7} that
\begin{equation}\label{eq-ratek2-8}
  \begin{aligned}[b]
     (1+t)^{\frac{3}{2} +\varepsilon}\| V_{xxx}(t) \|_{L^2(\mathbbm{R}_+)}^2+\int_0^t (1+\tau)^{\frac{3}{2} +\varepsilon} \| V_{xxx}(\tau) \|_{L^2(\mathbbm{R}_+)}^2\,\mathrm{d}\tau
    \le C(1+t)^\varepsilon\mathrm{log}^{10}(2+t).
  \end{aligned}
\end{equation}
Therefore, we complete the proof of \eqref{eq-rate-V-jg-3}. By using the similar method, we can prove \eqref{eq-rate-V-jg-4}. By utilizing\eqref{eq-P-Vrela} with $k=2,3$, combining \eqref{eq-rate-V-jg-1}-\eqref{eq-rate-V-jg-4} and \eqref{eq-PxH2-sjl}, we can prove \eqref{eq-P-rate}. Thus, we complete the proof of Proposition \ref{prop-Vrate-xsb}.
 $ \hfill\Box $

\begin{remark}We can not multiply the inequality \eqref{eq-rete-Vx-4} by $(1+t)^{\frac{5}{2}+\varepsilon }$ as in \cite{Gao2008}. Because the decay rate of the boundary term $f'(u_-)V_{xx}^2(0,t)$ is $O(t^{-\frac{3}{2} })$ if $f'(u_-)$ is not large enough by using $f'(u_-)V_{xx}^2(0,t)\le C\frac{1}{f'(u_-)}  P_{xx}^2(0,t)+C \mathrm{e}^{-c(1+t)}\le C \| P_{xx} \|_{H^1(\mathbbm{R}_+)}^2+C \mathrm{e}^{-c(1+t)}$. Although our decay rate is not as good as that of the Cauchy problem in \cite{Gao2008}, we focus on the stability of the initial-boundary value problem of \eqref{eq-yfc} rather than calculating the optimal decay rate.
\end{remark}
\begin{remark}
	For the case of $f'(u_-)=0$, we can get the better decay estimate like in \cite{Gao2008} owing to the boundary term $f'(u_-) \partial_x^kV(0,t)=0 \ \ (k=1,2,3,4)$. The details can refer to \cite{Gao2008}.
\end{remark}

By Sobolev inequality and Proposition \ref{prop-Vrate-xsb}, we obtain the following Corollary \ref{cor-Vinf-rate}.
\begin{cor}\label{cor-Vinf-rate}
  The solution $(V,P)$ of \eqref{eq-rd-VP}-\eqref{eq-VP-boundary} satisfies
\begin{equation}
  \begin{cases}
    \|  V(t) \|_{L^\infty(\mathbbm{R}_+)} \le C(1+t)^{-\frac{1}{2} }\mathrm{log}^3(2+t),\\
    \| \partial_x^kV(t) \|_{L^\infty(\mathbbm{R}_+)} \le C(1+t)^{-\frac{3}{4} }\mathrm{log}^5(2+t),\quad k=1,2,3,\\
    \| \partial_x^kV_t(t) \|_{L^\infty(\mathbbm{R}_+)} \le C(1+t)^{-\frac{3}{4} }\mathrm{log}^5(2+t),\quad k=0,1,2,\\
    \| \partial_x^kP(t) \|_{L^\infty(\mathbbm{R}_+)} \le C(1+t)^{-\frac{3}{4} }\mathrm{log}^5(2+t),\quad k=0,1,2,3,4.
  \end{cases}
\end{equation}
\end{cor}
Combining the results of Proposition \ref{prop-V} and Corollary \ref{cor-Vinf-rate}, we can complete the proof of Theorem \ref{thm-V-main}.\\
\indent Furthermore, by Corollaries \ref{cor-Vinf-rate}, \ref{prop-Vrate-xsb} and Lemma \ref{lem-jzxsb}, we get the following Corollary.
\begin{cor}\label{cor-Uinf-rate}
  The solution $(U,Q)$ of \eqref{eq-UQ-fc}-\eqref{eq-UQ-boundary} satisfies
\begin{equation}\label{eq-U-sjgj}
  \begin{cases}
    \|\partial_x^k U(t) \|_{L^2(\mathbbm{R}_+)}\le C \min\{\delta',(1+t)^{-\frac{3}{4} }\mathrm{log}^{5}(2+t)\},\quad k=2,3,4,\\
    \|\partial_x^k U_t(t) \|_{L^2(\mathbbm{R}_+)}\le C \min\{\delta',(1+t)^{-\frac{3}{4} }\mathrm{log}^{5}(2+t)\},\quad k=1,2,3,\\
    \| \partial_x^kU(t) \|_{L^\infty(\mathbbm{R}_+)} \le C \min\{\delta',(1+t)^{-\frac{3}{4} }\mathrm{log}^5(2+t)\},\quad k=2,3,\\
    \| \partial_x^kU_t(t) \|_{L^\infty(\mathbbm{R}_+)} \le C \min\{\delta',(1+t)^{-\frac{3}{4} }\mathrm{log}^5(2+t)\},\quad k=1,2,
  \end{cases}
\end{equation}
where $\delta'$ is defined in $\eqref{eq-def-dl}$.
\end{cor}

\subsection{Estimates for the perturbation on two-dimensional half space}\label{sec-xsb2}
In this section, we consider the initial-boundary value problem on two-dimensional half space:
\begin{equation}\label{eq-rdfc-vp}
  \begin{cases}
    v_t+(f(v+U)-f(U))_x+g(v+U)_y+{\rm{div}}{p}=0,\\[1mm]
    -\nabla{\rm{div}}{p}+p+\nabla{v}=0,
  \end{cases}
\end{equation}
with initial data
\begin{equation}
  v(x,y,0)=v_0(x,y)=u_0(x,y)-U_0(x),
\end{equation}
 and boundary condition
\begin{equation}\label{eq-boundary-vp}
  v(0,y,t)=0,\quad {\rm{div}}{p}(0,y,t)=0.
\end{equation}
By calculating the curl of \eqref{eq-rdfc-vp}$_2$, we can deduce that
\begin{equation}\label{eq-p1p2}
	p_{1y}=p_{2x},
\end{equation}
which will be used in estimating the perturbation $p$. Note that the equation \eqref{eq-rdfc-vp}$_2$ is equivalent to
\begin{equation}\label{eq-rd2-dengjia}
	\begin{cases}
		-{\rm{div}}{p}_x+p_1+{v}_x=0,\\
		-{\rm{div}}{p}_y+p_2+{v}_y=0.
	\end{cases}
\end{equation}

 The global existence follows from the combination of the local estimate Proposition \ref{prop-local-vp} and the a \emph{priori} estimates Proposition \ref{prop-priori-vp}. In this section, we will denote ourselves to establish a \emph{priori} estimates under the \emph{a priori} assumption
\begin{equation}\label{eq-pri-assump-H3}
  \sup\limits_{t\ge 0}\| v(t)\|_{H^3(\mathbbm{R}_+^2)}^2\le  \varepsilon_1^2,
\end{equation}
where $0<\varepsilon_1\ll 1$.
Then by the Sobolev inequality
\begin{equation}\label{eq-Soblev-f}
	\| f \|_{L^\infty(\mathbbm{R}_+^2)}^2\le C \| f \|_{H^2(\mathbbm{R}_+^2)}^2,
\end{equation}
we obtain
\begin{equation}\label{eq-pri-assump-Linf}
  \|\nabla v(t)\|_{L^\infty(\mathbbm{R}_+^2)}\le C \varepsilon_1.
\end{equation}
For any $0<T\le \infty$, we seek the solution of the initial boundary value problem \eqref{eq-rdfc-vp}-\eqref{eq-boundary-vp} in the set of functions $\tilde{X}(0,T)$ defined by
\begin{equation*}
\tilde{X}(0,T) = \left\{
\begin{tabular}{c|c}
   \multirow{2}{*}{({\it v, p})}    & $v\in C^0([0,T);H^3), \ \ \nabla{v}\in L^2(0,T;H^2)$   \\[1mm]
             & $p\in C^0([0,T);H^3)\cap L^2(0,T;H^3), \ \ {\rm{div}}{p}\in C^0([0,T);H^3)\cap L^2(0,T;H^3)$
\end{tabular}
       \right\}.
\end{equation*}
In order to state the results on the \emph{a priori} estimates, we define $M_0$ by
$$M_0^2:=\| v_0 \|_{H^3(\mathbbm{R}_+^2)}^2. $$
\begin{prop}[Local existence]\label{prop-local-vp}
Suppose the boundary condition satisfies $0\le f'(u_-)<f'(u_+)$ and the initial data satisfies $v_0 \in H^3(\mathbbm{R}_+^2)$. Also suppose that the initial data $\| v_0 \|_{H^3(\mathbbm{R}_+^2)} $ and $\delta$ are both small enough. Then there are two positive constants $C$ and $T_0$ such that the problem \eqref{eq-rdfc-vp}-\eqref{eq-boundary-vp} has a unique solution $(v,p)\in \tilde{X}(0,T_0)$, which satisfies
\begin{equation*}
\begin{aligned}[b]
\|v(t)\|_{H^3(\mathbbm{R}_+^2)}^2&+\|p(t)\|_{H^3(\mathbbm{R}_+^2)}^2+\| {\rm{div}}{p}(t) \|_{H^3(\mathbbm{R}_+^2)}^2 \\
&+\int_0^t \  ( \| \nabla{v}(\tau) \|_{H^2(\mathbbm{R}_+^2)}^2+\| p(\tau) \|_{H^3(\mathbbm{R}_+^2)}^2+\| {\rm{div}}{p}(\tau) \|_{H^3(\mathbbm{R}_+^2)}^2)  \,d{\tau} \le C\| v_0 \|_{H^3(\mathbbm{R}_+^2)}^2, \quad \forall t\in [0,T_0].	
\end{aligned}
\end{equation*}
\end{prop}

\begin{prop}[A priori estimates]\label{prop-priori-vp}
Let $T$ be a positive constant. Suppose that $0\le f'(u_-)<f'(u_+)$ and the problem \eqref{eq-rdfc-vp}-\eqref{eq-boundary-vp} has a unique solution $(v,p)\in \tilde{X}(0,T)$. Then there exists positive constants $C$ and $\delta_3'$ such that if $\| v_0 \|_{H^3(\mathbbm{R}_+^2)}+\delta\le \delta_3'$, $(0<\delta_3'\ll 1)$, then we get the estimate
\begin{equation*}
\begin{aligned}[b]
\|v(t)\|_{H^3(\mathbbm{R}_+^2)}^2&+\|p(t)\|_{H^3(\mathbbm{R}_+^2)}^2+\| {\rm{div}}{p}(t) \|_{H^3(\mathbbm{R}_+^2)}^2\\
&+\int_0^t \ ( \| \nabla v(\tau)\|_{H^2(\mathbbm{R}_+^2)}^2+\| p(\tau) \|_{H^3(\mathbbm{R}_+^2)}^2+\| {\rm{div}}{p}(\tau) \|_{H^3(\mathbbm{R}_+^2)}^2)  \,d{\tau}
\le C \| v_0 \|_{H^3(\mathbbm{R}_+^2)}^2, \quad \forall t\in [0,T].	
\end{aligned}
\end{equation*}
\end{prop}

\begin{lem}\label{lem-fv} Under the same assumptions of Proposition \ref{prop-priori-vp}, there exists a positive constant $C$ such that the following estimates hold:
\begin{flalign*}
\hspace{1mm}
&(C1) \ \ \ \iint_{\mathbbm{R}_+^2} (f(v+U)-f(U))_x {v} \,\mathrm{d}x\mathrm{d}y \ge \frac{\alpha}{2} \iint_{\mathbbm{R}_+^2} U_xv^2 \,\mathrm{d}x\mathrm{d}y;&\nonumber
\end{flalign*}
\begin{flalign*}
&(C2) \ \ \ \iint_{\mathbbm{R}_+^2} \nabla(f(v+U)-f(U))_x\cdot\nabla{v} \,\mathrm{d}x &\nonumber\\[1mm]
&\ \ \ \ \ \ \ \ \ \ \ge \frac{\alpha}{2} \iint_{\mathbbm{R}_+^2} U_x|\nabla{v}|^2 \,\mathrm{d}x\mathrm{d}y\mathrm{d}y+ \alpha\iint_{\mathbbm{R}_+^2} U_x|{v}_x|^2 \,\mathrm{d}x\mathrm{d}y- \frac{f'(u_-)}{2}|\nabla{v}(0,y,t)|^2 &\nonumber\\[1mm]
&\ \ \ \ \ \ \ \ \ \ -C \iint_{\mathbbm{R}_+^2} |\nabla{v}|(|\nabla{v}|^2+U_x^2|v|+| U_{xx} | |v|) \,\mathrm{d}x\mathrm{d}y; &\nonumber
\end{flalign*}
\begin{flalign*}
&(C3) \ \ \  \iint_{\mathbbm{R}_+^2} \Delta(f(v+U)-f(U))_x\Delta{v} \,\mathrm{d}x\mathrm{d}y &\nonumber\\[1mm]
&\ \ \ \ \ \ \ \ \ \ \ge \frac{\alpha}{2}\iint_{\mathbbm{R}_+^2} U_x|\Delta{v}|^2 \,\mathrm{d}x\mathrm{d}y- \frac{f'(u_-)}{2}(\Delta{v} (0,y,t))^2  \nonumber\\[1mm]
&\ \ \ \ \ \ \ \ \ \  -C\iint_{\mathbbm{R}_+^2} \left\{|\nabla{v}||\nabla^2{v}|^2+(|\nabla{v}|+|U_{xx}|+U_x)|\nabla{v}|\Delta v \right\} \,\mathrm{d}x\mathrm{d}y   &\nonumber\\[1mm]
&\ \ \ \ \ \ \ \ \ \  -C\iint_{\mathbbm{R}_+^2} (|U_{xxx}|+U_x|U_{xx}|+U_x^3)|v|(\Delta{v}) \,\mathrm{d}x\mathrm{d}y;   &\nonumber
\end{flalign*}
\begin{flalign*}
&(C4) \ \ \ \iint_{\mathbbm{R}_+^2} (f(v+U)-f(U))_{xyy}{v}_{yy} \,\mathrm{d}x\mathrm{d}y &\nonumber\\[1mm]
&\ \ \ \ \ \ \ \ \ \ \ge \frac{\alpha}{2}\iint_{\mathbbm{R}_+^2} U_x v_{yy}^2 \,\mathrm{d}x\mathrm{d}y-C \iint_{\mathbbm{R}_+^2}  \{|\nabla{v}| |\nabla^2{v}|^2+(|\nabla{v}|+U_x)|\nabla{v}||\nabla^2{v}|\}\,\mathrm{d}x\mathrm{d}y; &\nonumber
\end{flalign*}
\begin{flalign*}
&(C5) \ \ \ \iint_{\mathbbm{R}_+^2} \nabla(f(v+U)-f(U))_{xyy}\cdot\nabla{v}_{yy} \,\mathrm{d}x\mathrm{d}y &\nonumber\\[1mm]
&\ \ \ \ \ \ \ \ \ \  \ge \frac{\alpha}{2} \iint_{\mathbbm{R}_+^2} U_x|\nabla{v}_{yy}|^2 \,\mathrm{d}x\mathrm{d}y+ \alpha\iint_{\mathbbm{R}_+^2} U_x|{v}_{xyy}|^2 \,\mathrm{d}x\mathrm{d}y- \frac{f'(u_-)}{2}\int_{\mathbbm{R}}  |\nabla{v}_{yy}(0,y,t)|^2\,\mathrm{d}y &\nonumber\\[1mm]
&\ \ \ \ \ \ \ \ \ \  -C \iint_{\mathbbm{R}_+^2} \{(|\nabla{v}|+U_x)|\nabla{v}||\nabla{v}_{yy}|+(|\nabla{v}|+U_x+|U_{xx}|)|\nabla^2{v}||\nabla{v}_{yy}|\} \,\mathrm{d}x \mathrm{d}y &\nonumber\\[1mm]
&\ \ \ \ \ \ \ \ \ \  -C \iint_{\mathbbm{R}_+^2} \{ |\nabla^2{v}|^2|\nabla{v}_{yy}|+(|\nabla{v}|+U_x)|\nabla^3{v}|^2\} \,\mathrm{d}x \mathrm{d}y;  &\nonumber
\end{flalign*}
\begin{flalign*}
&(C6) \ \ \ \iint_{\mathbbm{R}_+^2} \nabla\Delta(f(v+U)-f(U))_x\cdot\nabla\Delta{v} \,\mathrm{d}x\mathrm{d}y &\nonumber\\[1mm]
&\ \ \ \ \ \ \ \ \ \  \ge \frac{\alpha}{2}\iint_{\mathbbm{R}_+^2} U_x|\nabla \Delta{v}|^2 \,\mathrm{d}x\mathrm{d}y+ \alpha\iint_{\mathbbm{R}_+^2} U_x| \Delta{v}_x|^2 \,\mathrm{d}x\mathrm{d}y- \frac{f'(u_-)}{2}|\nabla\Delta{v}(0,y,t)|^2  &\nonumber\\[1mm]
&\ \ \ \ \ \ \ \ \ \  -C (\| \nabla{v} \|_{L^\infty(\mathbbm{R}_+)}+\sum_{k=1}^3\| \partial_x^k U \|_{L^\infty(\mathrm{R}_+^2)} )\iint_{\mathbbm{R}_+^2} (|\nabla{v}|^2+|\nabla^2{v}|^2+|\nabla^2{v}|^4+|\nabla^3{v}|^2 ) \,\mathrm{d}x\mathrm{d}y, & \nonumber\\[1mm]
&\ \ \ \ \ \ \ \ \ \  -C \iint_{\mathbbm{R}_+^2} (| U_{xxxx}|+|U_x|^4+|U_{xx}|^2+|U_{x}||U_{xxx}|+|U_{x}|^2|U_{xx}|) |v||\Delta{v}_x| \,\mathrm{d}x\mathrm{d}y; & \nonumber
\end{flalign*}
\begin{flalign*}
&(C7) \ \ \ \iint_{\mathbbm{R}_+^2} (g(v+U)_y) v \,\mathrm{d}x \mathrm{d}y=0; & \nonumber\\[2mm]
&(C8) \ \ \ \iint_{\mathbbm{R}_+^2} \nabla (g(v+U)_y)\cdot\nabla v \,\mathrm{d}x \mathrm{d}y \le C(\|\nabla{v}\|_{L^\infty(\mathrm{R}_+^2)}+\| U_x \|_{L^\infty(\mathrm{R}_+^2)})\iint_{\mathbbm{R}_+^2} |\nabla{v}|^2 \,\mathrm{d}x \mathrm{d}y; & \nonumber\\[2mm]
&(C9) \ \ \ \iint_{\mathbbm{R}_+^2} \Delta (g(v+U)_y)\Delta v \,\mathrm{d}x \mathrm{d}y \le C(\|\nabla{v}\|_{L^\infty(\mathrm{R}_+^2)}+\sum\limits_{k=1}^2\| \partial_x^k U \|_{L^\infty(\mathrm{R}_+^2)})\iint_{\mathbbm{R}_+^2} (|\nabla{v}|^2+|\nabla^2{v}|^2) \,\mathrm{d}x \mathrm{d}y; & \nonumber\\[2mm]
&(C10) \ \ \ \iint_{\mathbbm{R}_+^2} (g(v+U)_{yyy} v_{yy} \,\mathrm{d}x \mathrm{d}y \le C(\|\nabla{v}\|_{L^\infty(\mathrm{R}_+^2)}+\| U_x \|_{L^\infty(\mathrm{R}_+^2)})\iint_{\mathbbm{R}_+^2} (|\nabla{v}|^2+|\nabla^2{v}|^2) \,\mathrm{d}x \mathrm{d}y; & \nonumber
\end{flalign*}
\begin{flalign*}
&(C11) \ \ \ \iint_{\mathbbm{R}_+^2} \nabla(g(v+U)_{yyy})\cdot \nabla v_{yy} \,\mathrm{d}x \mathrm{d}y  & \nonumber\\[1mm]
&\ \ \ \ \ \ \ \ \ \  \le C\iint_{\mathbbm{R}_+^2} \{(| \nabla{v} |+|U_x| )|\nabla{v}|+(| \nabla{v} |+| U_x | )|\nabla^2{v}|+|\nabla^2{v}|^2+|\nabla{v}||\nabla{v}_{yy}|\}|\nabla{v}_{yy}| \,\mathrm{d}x \mathrm{d}y; & \nonumber
\end{flalign*}
\begin{flalign*}
&(C12) \ \ \ \iint_{\mathbbm{R}_+^2} \nabla\Delta(g(v+U)_{y}) \nabla\Delta v \,\mathrm{d}x \mathrm{d}y  & \nonumber\\[1mm]
&\ \ \ \ \ \ \ \ \ \  \le C(\| \nabla{v} \|_{L^\infty(\mathbbm{R}_+^2)}+\| U_x \|_{L^\infty(\mathbbm{R}_+)})\iint_{\mathbbm{R}_+^2} (|\nabla{v}|^2+|\nabla^2{v}|^2+|\nabla^3{v}|^2) \,\mathrm{d}x \mathrm{d}y+C \iint_{\mathbbm{R}_+^2} |\nabla^2{v}|^4 \,\mathrm{d}x \mathrm{d}y. & \nonumber
\end{flalign*}
\end{lem}
{\it\bfseries Proof.}
We only give the proof of $(C3)$ and $(C6)$, and the rest can be calculated in the similar way. Notice that $\Delta (fg)=(\Delta f)g+ f(\Delta g)+2\nabla f\cdot\nabla g$. By direct calculation, we have
\begin{equation}\label{eq-delfv}
	\begin{aligned}[b]
		\Delta(f(v+U)-f(U))_x\Delta{v}
		&= \Delta f'(v+U)(v_x+U_x)\Delta{v}+f'(v+U)(\Delta v_x+U_{xxx})\Delta{v}\\
		& \ \ +2f''(v+U)(\nabla v\cdot\nabla v_x+v_xU_{xx}+U_xv_{xx}+U_xU_{xx})\Delta{v}\\
		& \ \ -\Delta f'(U)U_x\Delta{v}-f'(U)U_{xxx}\Delta{v}-2f''(U)U_xU_{xx}\Delta{v}\\
		&=f'''(v+U)(v_x^3+U_xv_x^2+U_x^2v_x)\Delta{v}+f''(v+U)(v_x+U_x)(\Delta{v})^2\\
		& \ \ +f''(v+u)U_{xx}v_x\Delta{v}+f'''(v+U)(v_x+U_x)v_y^2(\Delta{v})+f(v+U)\Delta{v}_x\Delta{v}\\
		& \ \ +2f''(v\!+\!U)(\nabla v\cdot\nabla v_x\!+\!v_xU_{xx}\!+\!U_xv_{xx})\Delta{v}\!+\!(f'''(v+U)\!-\!f'''(U))U_x^3 \Delta{v}\\
		& \ \ +(f'(v+U)-f'(U))U_{xxx} \Delta{v}+3(f''(v+U)-f''(U))U_xU_{xx} \Delta{v},
	\end{aligned}
\end{equation}
where we have used the fact that
$$\Delta f'(v+U)=f'''(v+U)(v_x+U_x)^2+f''(v+U)(\Delta{v}+U_{xx})+f'''(v+U)v_y^2.$$
Integrating \eqref{eq-delfv} over $\mathbbm{R}_+^2$, we can obtain $(C3)$. Based on the equations, we have
\begin{equation}\label{eq-nabdelfv}
	\begin{aligned}[b]
		\nabla\Delta&(f(v+U)-f(U))_x\cdot\nabla\Delta{v}\\
		=& \nabla\Delta f'(v+U)\cdot\nabla\Delta{v}(v_x+U_x)+\Delta f'(v+U)(\nabla v_x\cdot \nabla\Delta v+U_{xx}\Delta v)\\
		& +f''(v+U)(\nabla{v}\cdot \nabla\Delta{v}+U_x\Delta v_x)(\Delta v_x+U_{xxx})+f'(v+U)(\nabla\Delta v_x\cdot \nabla\Delta{v}+U_{xxxx}\Delta v_x)\\
		& +2f'''(v+U)(\nabla{v}\cdot \nabla\Delta{v}+U_x\Delta v_x)(\nabla v\cdot\nabla v_x+v_xU_{xx}+U_xv_{xx}+U_xU_{xx})\\
		& +2f''(v+U)(\nabla v_x\cdot \nabla\Delta{v}v_{xx}+v_x\nabla v_{xx}\cdot \nabla\Delta{v}+\nabla v_y\cdot \nabla\Delta{v}v_{xy}+v_y\nabla v_{xy}\cdot \nabla\Delta{v}+\nabla v_x\cdot \nabla\Delta{v}U_{xx})\\
		& +2f''(v+U)(v_xU_{xxx}\Delta v_x+v_{xx}U_{xx}\Delta v_x+U_{x}\nabla v_{xx}\cdot\nabla\Delta v+U_{xx}^2\Delta v_x+U_{x}U_{xxx}\Delta v_x)\\
		& -\left(f^{(4)}(U)U_x^4+6f'''(U)U_x^2U_{xx}+3f'(U)U_{xx}^2+4f'(U)U_xU_{xxx}+f'(U)U_{xxxx}\right)\Delta v_x.
	\end{aligned}
\end{equation}
The first term on the right-hand side is equal to
\begin{equation}\label{eq-nabdelfv1}
	\begin{aligned}[b]
		\nabla\Delta& f'(v+U)\cdot\nabla\Delta{v}(v_x+U_x)\\
		&=f^{(4)}(v+U)(\nabla v\cdot\nabla\Delta v+U_x\Delta v_x)(v_x+U_x)^3+2f'''(v+U)(v_x+U_x)^2(\nabla{v}\cdot \nabla\Delta{v}+U_x\Delta v_x)\\
		&  +f'''(v+U)(\nabla v\cdot\nabla\Delta v+U_x\Delta v_x)(\Delta v+U_{xx})(v_x+U_x)+f''(v+U)(|\nabla\Delta{v}|^2+U_{xxx}\Delta v_x)(v_x+U_x)\\
		&  +f^{(4)}(v+U)(\nabla{v}\cdot \nabla\Delta{v}+U_x\Delta v_x)v_y^2(v_x+U_x)+2f'''(v+U)v_y\nabla{v}_y\cdot\nabla\Delta{v}(v_x+U_x).
	\end{aligned}
\end{equation}
Substituting \eqref{eq-nabdelfv1} into \eqref{eq-nabdelfv} and then integrating it over $\mathbbm{R}_+^2$, by using Cauchy inequality, the estimate $(C6)$ can be obtained.
 $ \hfill\Box$
\begin{lem}\label{lem-bjv-gj} The solution $v(x,y,t)$ of \eqref{eq-rdfc-vp}-\eqref{eq-boundary-vp} satisfies the following boundary estimates:
\begin{flalign*}
\hspace{3mm}
      &(D1) \ \  \partial_y^k \partial_t^l v(0,t)=0,\quad k,l=0,1,2, \cdots; &\nonumber\\[1mm]
\hspace{3mm}
      &(D2) \ \  v_x(0,y,t)=0; &\nonumber\\[1mm]
\hspace{3mm}
      &(D3) \ \  v_{xx}^2(0,y,t)\le C ({\rm{div}}{p}_{x}(0,y,t))^2; &\nonumber\\[1mm]
\hspace{3mm}
      &(D4) \ \  v_{xxy}^2(0,y,t)\le C ({\rm{div}}{p}_{xy}(0,y,t))^2; &\nonumber\\[1mm]
\hspace{3mm}
      &(D5) \ \  v_{xxx}^2(0,y,t)\le C ({\rm{div}}{p}_{x}(0,t))^2+C ({\rm{div}}{p}_{xt}(0,y,t))^2+C ({\rm{div}}{p}_{xy}(0,y,t))^2+C ({\rm{div}}{p}_{xx}(0,y,t))^2. &\nonumber
\end{flalign*}
\end{lem}

\subsubsection{\emph{A priori} estimates}
\begin{lem}\label{lem-nabv} Under the same assumptions of Proposition \ref{prop-priori-vp}, there exists a positive constant $C$ such that for $0\le t\le T$,
\begin{equation}\label{eq-nabv-jg}
  \begin{aligned}[b]
    \| v(t) \|_{H^1(\mathbbm{R}_+^2)}^2&+\int_0^t (\| \sqrt{U_x}v(t) \|_{L^2(\mathbbm{R}_+^2)}^2+\| \sqrt{U_x}\nabla v(t) \|_{L^2(\mathbbm{R}_+^2)}^2 ) \,\mathrm{d}\tau\\
    &+\int_0^t (\| {\rm{div}}{p}(t) \|_{H^1(\mathbbm{R}_+^2)}^2 +\| p(t) \|_{L^2(\mathbbm{R}_+^2)}^2+\| \nabla{v}(t) \|_{L^2(\mathbbm{R}_+^2)}^2) \,\mathrm{d}\tau\le C M_0^2.
  \end{aligned}
\end{equation}
\end{lem}
 {\it\bfseries Proof.}
We can get from $\eqref{eq-rdfc-vp}_1\times v+\eqref{eq-rdfc-vp}_2\cdot p$ that
\begin{equation}\label{eq-v-1}
  \begin{aligned}[b]
    &\frac{1}{2} \frac{\mathrm{d}}{\mathrm{d}t}v^2+(f(v+U)-f(U))_xv+(g(v+U)_y)v+({\rm{div}}{p})^2+|p|^2+\mathrm{div}\{ pv- {\rm{div}}{p} p \}=0.
  \end{aligned}
\end{equation}
By Lemma \ref{lem-fv} $\mathrm{(C1)},$ $\mathrm{(C7)}$ and \eqref{eq-boundary-vp}, integrating \eqref{eq-v-1} over $\mathbbm{R}_+^2\times[0,t]$, we have
\begin{equation}\label{eq-v-gj}
  \begin{aligned}[b]
    \| v(t) \|_{L^2(\mathbbm{R}_+^2)}^2+\int_0^t \| \sqrt{U_x}v(\tau) \|_{L^2(\mathbbm{R}_+)}^2  \,\mathrm{d}\tau+\int_0^t (\| {\rm{div}}{p}(\tau) \|_{L^2(\mathbbm{R}_+^2)}^2+\| p(\tau) \|_{L^2(\mathbbm{R}_+^2)}^2)  \,\mathrm{d}\tau \le C M_0^2.
  \end{aligned}
\end{equation}
On the other hand, we can get from $\nabla\eqref{eq-rdfc-vp}_1\cdot\nabla v- \eqref{eq-rdfc-vp}_2\cdot\nabla{\rm{div}}{p}$ that
\begin{equation}\label{eq-vx-2}
  \begin{aligned}[b]
    \frac{1}{2} \frac{\mathrm{d}}{\mathrm{d}t}|\nabla{v}|^2+\nabla(f(v+U)-f(U))_x\cdot\nabla{v}+\nabla(g(v+U)_y)\cdot\nabla{v}+|\nabla{\rm{div}}{p}|^2+({\rm{div}}{p})^2- \mathrm{divp}\{{\rm{div}}{p}p\}=0.
  \end{aligned}
\end{equation}
Integrating $\eqref{eq-vx-2}$ over $\mathbbm{R}_+^2$, using Lemma \ref{lem-fv} $\mathrm{(C2)},$ $\mathrm{(C8)}$ and Lemma \ref{lem-bjv-gj} $\mathrm{(D1)}$-$\mathrm{(D2)}$, combining \eqref{eq-v-gj} and \eqref{eq-U-sjgj}, we have
\begin{equation}\label{eq-vx-3}
  \begin{aligned}[b]
    &\frac{\mathrm{d}}{\mathrm{d}t}\|\nabla v(t) \|_{L^2(\mathbbm{R}_+^2)}^2+ \| \sqrt{U_x}\nabla v(t) \|_{L^2(\mathbbm{R}_+^2)}^2  + \| \nabla {\rm{div}}{p}(t) \|_{L^2(\mathbbm{R}_+^2)}^2+\| {\rm{div}}{p}(t) \|_{L^2(\mathbbm{R}_+^2)}^2 \\
     \le& C\| \nabla{v} \|_{L^\infty(\mathbbm{R}_+^2)}\| \nabla{v}(t) \|_{L^2(\mathbbm{R}_+^2)}^2+ C \| v(t) \|_{L^2(\mathbbm{R}_+^2)}^2  \| U_{xx} \|_{L^\infty(\mathbbm{R}_+)}^2+C \| \sqrt{U_x}v(t) \|_{L^2(\mathbbm{R}_+^2)}^2 \\
     \le& C\varepsilon_1\| \nabla{v}(t) \|_{L^2(\mathbbm{R}_+^2)}^2+ C M_0^2(1+t)^{-\frac{3}{2}}\mathrm{log}^{10}(2+t)+C \| \sqrt{U_x}v(t) \|_{L^2(\mathbbm{R}_+^2)}^2.
  \end{aligned}
\end{equation}
From $\eqref{eq-rdfc-vp}_2$, the first term on the right-hand side can be estimated as
\begin{equation}\label{eq-nabv-p}
  \| \nabla{v} \|_{L^2(\mathbbm{R}_+^2)}^2\le 2(\| \nabla{\rm{div}}{p} \|_{L^2(\mathbbm{R}_+^2)}^2+\| p \|_{L^2(\mathbbm{R}_+^2)}^2  ),
\end{equation}
it follows from $\eqref{eq-vx-3}$ after integration over $[0,t]$ that
\begin{equation*}
  \begin{aligned}[b]
    \| \nabla v(t) \|_{L^2(\mathbbm{R}_+^2)}^2+\int_0^t \| \sqrt{U_x}\nabla v(\tau) \|_{L^2(\mathbbm{R}_+)}^2  \,\mathrm{d}\tau+\int_0^t (\| \nabla{\rm{div}}{p}(\tau) \|_{L^2(\mathbbm{R}_+^2)}^2+\| {\rm{div}}{p}(\tau) \|_{L^2(\mathbbm{R}_+^2)}^2)  \,\mathrm{d}\tau
    \le C M_0^2,
  \end{aligned}
\end{equation*}
for some small $\varepsilon_1$. Using $\eqref{eq-nabv-p}$ again, we get
\begin{equation*}
  \int_0^t \| \nabla{v}(\tau) \|_{L^2(\mathbbm{R}_+^2)}^2  \,\mathrm{d}\tau\le  C M_0^2,
\end{equation*}
which completes the proof of Lemma \ref{lem-nabv}.
 $ \hfill\Box $

\begin{lem}\label{lem-nab2v} Under the same assumptions of Proposition \ref{prop-priori-vp}, there exists a positive constant $C$ such that for $0\le t\le T$,
\begin{equation}\label{eq-nab2divp-jg}
\begin{aligned}
  \| \nabla^2 v(t) \|_{L^2(\mathbbm{R}_+^2)}^2+&\int_0^t (\| \sqrt{U_x}\Delta v(\tau) \|_{L^2(\mathbbm{R}_+^2)}^2 +\| \nabla{\rm{div}}{p}(\tau) \|_{H^1(\mathbbm{R}_+^2)}^2+\| \nabla^2{v}(\tau) \|_{L^2(\mathbbm{R}_+^2)}^2) \,\mathrm{d}\tau\le CM_0^2.
\end{aligned}
\end{equation}
\end{lem}
 {\it\bfseries Proof.}
We can get from $\Delta\eqref{eq-rdfc-vp}_1\times\Delta v- {\rm{div}}\eqref{eq-rdfc-vp}_2\times\Delta {\rm{div}}{p}$ that
\begin{equation}\label{eq-delv-1}
  \begin{aligned}[b]
    \frac{1}{2} \frac{\mathrm{d}}{\mathrm{d}t}(\Delta{v})^2+\Delta((f(v+U)-f(U))_x)\Delta v+\Delta(g(v+U)_y)\Delta v+(\Delta{\rm{div}}{p})^2+|\nabla{\rm{div}}{p}|^2-{\rm{div}}\{{\rm{div}}{p} \nabla{\rm{div}}{p}\}=0.
  \end{aligned}
\end{equation}
Notice that
$$(\Delta {\rm{div}}{p})^2=({\rm{div}}{p}_{xx})^2+({\rm{div}}{p}_{yy})^2+2({\rm{div}}{p}_{xy})^2+2\{{\rm{div}}{p}_{xx}{\rm{div}}{p}_y\}_y-2\{{\rm{div}}{p}_{xy}{\rm{div}}{p}_y\}_x,$$
$$(\Delta v)^2=(v_{xx})^2+(v_{yy})^2+2(v_{xy})^2+2\{v_{xx}v_y\}_y-2\{v_{xy}v_y\}_x,$$
which implies that
\begin{align}
  &\iint_{\mathbbm{R}_+^2} (\Delta {\rm{div}}{p})^2 \,\mathrm{d}x \mathrm{d}y =\iint_{\mathbbm{R}_+^2} |\nabla^2 {\rm{div}}{p}|^2 \,\mathrm{d}x \mathrm{d}y,\label{eq-deldivp-zh}\\
  &\iint_{\mathbbm{R}_+^2} (\Delta v)^2 \,\mathrm{d}x \mathrm{d}y =\iint_{\mathbbm{R}_+^2} |\nabla^2 v|^2 \,\mathrm{d}x \mathrm{d}y.\label{eq-delv-zh}
\end{align}
Integrating $\eqref{eq-delv-1}$ over $\mathbbm{R}_+^2$, by utilizing Lemma \ref{lem-fv} $\mathrm{(C3)}$, $\mathrm{(C9)}$, $\eqref{eq-deldivp-zh}$ and Cauchy inequality, we obtain
\begin{equation}\label{eq-delv-2}
  \begin{aligned}[b]
    &\frac{1}{2} \frac{\mathrm{d}}{\mathrm{d}t}\| \nabla^2{v}(t) \|_{L^2(\mathbbm{R}_+^2)}^2+\frac{\alpha}{2}\| \sqrt{U_x}\Delta v(t) \|_{L^2(\mathbbm{R}_+^2)}^2  +\| \nabla^2{\rm{div}}{p}(t) \|_{L^2(\mathbbm{R}_+^2)}^2+\| \nabla{\rm{div}}{p}(t) \|_{L^2(\mathbbm{R}_+^2)}^2  \\
    \le& Cv_{xx}^2(0,y,t)+C(\| \nabla{v} \|_{L^\infty(\mathbbm{R}_+^2)}+\mu  )\| \nabla^2{v}(t) \|_{L^2(\mathbbm{R}_+^2)}^2\\
    & +C\mu^{-1}\| \nabla{v}(t) \|_{L^2(\mathbbm{R}_+^2)}^2+C\mu^{-1} \| v(t) \|_{L^2(\mathbbm{R}_+^2)}^2 \| U_{xxx}(t) \|_{L^\infty(\mathbbm{R}_+)}^2+C\mu^{-1}  \| \sqrt{U_{x}}v(t) \|_{L^2(\mathbbm{R}_+)}^2.
  \end{aligned}
\end{equation}
By Lemma \ref{lem-bjv-gj} $(D3)$, the first term on the right-hand side of \eqref{eq-delv-2}, we have
\begin{equation}\label{eq-delv-3}
  \begin{aligned}[b]
    v_{xx}^2(0,y,t)\le \frac{1}{4} \| {\rm{div}}{p}_{xx}(t) \|_{L^2(\mathbbm{R}_+^2)}^2+ C \| {\rm{div}}{p}_{x}(t) \|_{L^2(\mathbbm{R}_+^2)}^2.
  \end{aligned}
\end{equation}
From $\eqref{eq-delv-zh}$ and $\eqref{eq-rdfc-vp}_2$, we can treat the second term as
\begin{equation}\label{eq-delv-4}
  \begin{aligned}[b]
    \| \nabla^2{v}(t) \|_{L^2(\mathbbm{R}_+^2)}^2
    =\| \Delta{v}(t) \|_{L^2(\mathbbm{R}_+^2)}^2
    \le 2(\| \Delta {\rm{div}}{p}(t) \|_{L^2(\mathbbm{R}_+^2)}^2+\| {\rm{div}}{p}(t) \|_{L^2(\mathbbm{R}_+^2)}^2 ).
  \end{aligned}
\end{equation}
The former term of $\eqref{eq-delv-4}$ can absorbed in the third term on the left-hand side of $\eqref{eq-delv-2}$ if $\varepsilon_1$ and $\mu$ are small enough. Substituting $\eqref{eq-delv-3}$ and $\eqref{eq-delv-4}$ into $\eqref{eq-delv-2}$, we can deduce that
\begin{equation*}
  \begin{aligned}[b]
    &\frac{\mathrm{d}}{\mathrm{d}t}\| \nabla^2{v}(t) \|_{L^2(\mathbbm{R}_+^2)}^2+\| \sqrt{U_x}\Delta v(t) \|_{L^2(\mathbbm{R}_+^2)}^2  +\| \nabla^2{\rm{div}}{p}(t) \|_{L^2(\mathbbm{R}_+^2)}^2+\| \nabla{\rm{div}}{p}(t) \|_{L^2(\mathbbm{R}_+^2)}^2  \\
    &\le C\| {\rm{div}}{p}(t) \|_{H^1(\mathbbm{R}_+^2)}^2  +C \|\nabla{v}(t)\|_{L^2(\mathbbm{R}_+^2)}^2+C\| \sqrt{U_{x}}v(t) \|_{L^\infty(\mathbbm{R}_+)}^2+C M_0^2(1+t)^{-\frac{3}{2} }\mathrm{log}^{10}(2+t).
  \end{aligned}
\end{equation*}
Integrating the above inequality over $[0,t]$, we get
\begin{equation*}
\begin{aligned}
  \| \nabla^2 v(t) \|_{L^2(\mathbbm{R}_+^2)}^2+&\int_0^t (\| \sqrt{U_x}\Delta v(\tau) \|_{L^2(\mathbbm{R}_+^2)}^2 +\| \nabla{\rm{div}}{p}(\tau) \|_{H^1(\mathbbm{R}_+^2)}^2) \,\mathrm{d}\tau\le CM_0^2.
\end{aligned}
\end{equation*}
 Using \eqref{eq-delv-4} again, we complete the proof of Lemma \ref{lem-nab2v}.
 $ \hfill\Box $

\begin{lem}\label{lem-vyy-H1-jg} Under the same assumptions of Proposition \ref{prop-priori-vp}, there exists a positive constant $C$ such that for $0\le t\le T$,
\begin{equation}\label{eq-vyy-H1-jg}
  \begin{aligned}[b]
    &\| v_{yy}(t) \|_{H^1(\mathbbm{R}_+^2)}^2+\int_0^t ( \| \sqrt{U_x}v_{yy}(\tau) \|_{L^2(\mathbbm{R}_+^2)}^2 +\frac{\alpha}{2} \| \sqrt{U_x}\nabla v_{yy}(\tau) \|_{L^2(\mathbbm{R}_+^2)}^2) \,\mathrm{d}\tau\\
    &\ \ \ \ \ \ \ \ \ \ \ +\int_0^t (\| \nabla{v}_{yy}(\tau) \|_{L^2(\mathbbm{R}_+^2)}^2+\|  {\rm{div}}{p}_{yy}(\tau) \|_{H^1(\mathbbm{R}_+^2)}^2+\| {p}_{yy}(\tau) \|_{L^2(\mathbbm{R}_+^2)}^2) \,\mathrm{d}\tau\\
    \le& CM_0^2+C(M_0+\delta')\int_0^t \| \nabla^3{v}(\tau) \|_{L^2(\mathbbm{R}_+^2)}^2 \,\mathrm{d}\tau.
  \end{aligned}
\end{equation}
\end{lem}
 {\it\bfseries Proof.}
We can get from $\partial_y^2\eqref{eq-rdfc-vp}_1\times v_{yy}+\eqref{eq-rdfc-vp}_2\cdot p_{yy}+\nabla\partial_y^2\eqref{eq-rdfc-vp}_1\cdot \nabla v_{yy}-\partial_y^2\eqref{eq-rdfc-vp}_2\cdot\nabla{\rm{div}}{p}_{yy}$ that
\begin{equation}\label{eq-v-H1-1}
  \begin{aligned}[b]
    \frac{1}{2} \frac{\mathrm{d}}{\mathrm{d}t}(v_{yy}^2+|\nabla{v}_{yy}|^2)&+(f(v+U)-f(U))_{xyy}v_{yy}+\nabla(f(v+U)-f(U))_{xyy}\cdot\nabla v_{yy}+g(v+U)_{yyy}v_{yy}\\
    &+\nabla(g(v+U)_{yyy})\cdot\nabla v_{yy}+|\nabla{\rm{div}}{p}_{yy}|^2-2 \nabla{\rm{div}}{p}_{yy}\cdot p_{yy}+|p_{yy}|^2+{\rm{div}}\{v_{yy}p_{yy}\}=0.
  \end{aligned}
\end{equation}
Integrating $\eqref{eq-v-H1-1}$ over $\mathbbm{R}_+^2$, by utilizing Lemma \ref{lem-fv} $\mathrm{(C4)}$-$\mathrm{(C5)}$ and $\mathrm{(C10)}$-$\mathrm{(C11)}$, we have
\begin{equation}\label{eq-v-H1-2}
  \begin{aligned}[b]
    &\frac{1}{2} \frac{\mathrm{d}}{\mathrm{d}t}\| v_{yy}(t) \|_{H^1(\mathbbm{R}_+^2)}^2+\frac{\alpha}{2} \| \sqrt{U_x}v_{yy}(t) \|_{L^2(\mathbbm{R}_+^2)}^2 +\frac{\alpha}{2} \| \sqrt{U_x}\nabla v_{yy}(t) \|_{L^2(\mathbbm{R}_+^2)}^2+\| \nabla{v}_{yy}(t) \|_{L^2(\mathbbm{R}_+^2)}^2\\
    \le& C(\|\nabla{v}  \|_{L^\infty(\mathbbm{R}_+^2)}+\| U_x \|_{L^\infty(\mathbbm{R}_+)}  )\| \nabla^3{v}(t) \|_{L^2(\mathbbm{R}_+^2)}^2+C\| \nabla{v}(t) \|_{H^1(\mathbbm{R}_+^2)}^2+ C\iint_{\mathbbm{R}_+^2} |\nabla^2{v}|^4 \,\mathrm{d}x \mathrm{d}y.
  \end{aligned}
\end{equation}
Here we have used the fact from $\eqref{eq-rdfc-vp}_2$ that
$$|\nabla{v}_{yy}|^2=|\nabla{\rm{div}}{p}_{yy}|^2-2 \nabla{\rm{div}}{p}_{yy}\cdot p_{yy}+|p_{yy}|^2.$$
The last term of \eqref{eq-v-H1-2} can be estimated as
\begin{equation}\label{eq-v-H1-3}
  \begin{aligned}[b]
    C\iint_{\mathbbm{R}_+^2} |\nabla^2{v}|^4 \,\mathrm{d}x \mathrm{d}y
    \le& C\int_{\mathbbm{R}_+} \sup\limits_{y\in\mathbbm{R}}|\nabla^2{v}|^2  \,\mathrm{d}x \int_{\mathbbm{R}} \sup\limits_{x\in\mathbbm{R}_+}|\nabla^2{v}|^2 \,\mathrm{d}y\\
    \le&C\int_{\mathbbm{R}_+} \|\nabla^2{v}(x,\cdot,t)\|_{L^2(\mathbbm{R})}\|\nabla^2{v}_y(x,\cdot,t)\|_{L^2(\mathbbm{R})}  \,\mathrm{d}x \int_{\mathbbm{R}} \|\nabla^2{v}(\cdot,y,t)\|_{L^2(\mathbbm{R}_+)}\|\nabla^2{v}_x(\cdot,y,t)\|_{L^2(\mathbbm{R}_+)} \,\mathrm{d}y\\
    \le&C\| \nabla^2{v}(t) \|_{L^2(\mathbbm{R}_+^2)}^2\| \nabla^3{v}(t) \|_{L^2(\mathbbm{R}_+^2)}^2\\
    \le& CM_0^2\| \nabla^3{v}(t) \|_{L^2(\mathbbm{R}_+^2)}^2.
  \end{aligned}
\end{equation}
Substituting $\eqref{eq-v-H1-3}$ into $\eqref{eq-v-H1-2}$, we obtain
\begin{equation}\label{eq-v-H1-4}
  \begin{aligned}[b]
    &\frac{\mathrm{d}}{\mathrm{d}t}\| v_{yy}(t) \|_{H^1(\mathbbm{R}_+^2)}^2+ \| \sqrt{U_x}v_{yy}(t) \|_{L^2(\mathbbm{R}_+^2)}^2 +\frac{\alpha}{2} \| \sqrt{U_x}\nabla v_{yy}(t) \|_{L^2(\mathbbm{R}_+^2)}^2+\| \nabla{v}_{yy}(t) \|_{L^2(\mathbbm{R}_+^2)}^2\\
    \le& C \| \nabla{v}(t) \|_{H^1(\mathbbm{R}_+^2)}^2+C(M_0+\delta')\| \nabla^3{v}(t) \|_{L^2(\mathbbm{R}_+^2)}^2.
  \end{aligned}
\end{equation}
Integrating $\eqref{eq-v-H1-4}$ over $[0,t]$, we can get
\begin{equation*}
  \begin{aligned}[b]
    \| v_{yy}(t) \|_{H^1(\mathbbm{R}_+^2)}^2&+\int_0^t ( \| \sqrt{U_x}v_{yy}(\tau) \|_{L^2(\mathbbm{R}_+^2)}^2 + \| \sqrt{U_x}\nabla v_{yy}(\tau) \|_{L^2(\mathbbm{R}_+^2)}^2) \,\mathrm{d}\tau\\
    & +\int_0^t \| \nabla{v}_{yy}(\tau) \|_{L^2(\mathbbm{R}_+^2)}^2 \,\mathrm{d}\tau\le C M_0^2+C(M_0+\delta')\int_0^t \| \nabla^3{v}(\tau) \|_{L^2(\mathbbm{R}_+^2)}^2 \,\mathrm{d}\tau.
  \end{aligned}
\end{equation*}
Furthermore, the equation $\eqref{eq-rdfc-vp}_2$ gives
\begin{equation*}
  |\nabla{\rm{div}}{p}_{yy}|^2+|p_{yy}|^2+2({\rm{div}}{p}_{yy})^2-2 {\rm{div}}\{ {\rm{div}}{p}_{yy}p_{yy}\}=|\nabla{v}_{yy}|^2,
\end{equation*}
which implies
\begin{equation*}
  \| {\rm{div}}{p}_{yy}(t) \|_{H^1(\mathbbm{R}_+^2)}^2+\| p_{yy}(t) \|_{L^2(\mathbbm{R}_+^2)}^2 \le \| \nabla v_{yy}(t) \|_{L^2(\mathbbm{R}_+^2)}^2,
\end{equation*}
owing to ${\rm{div}}{p}_{yy}(0,y,t)=0$.

Integrating the above inequality over $[0,t]$, the desired estimates $\eqref{eq-vyy-H1-jg}$ can be obtained, which completes the proof of Lemma \ref{lem-vyy-H1-jg}.
 $ \hfill\Box $
\begin{lem}\label{lem-vtH1-jg} Under the same assumptions of Proposition \ref{prop-priori-vp}, there exists a positive constant $C$ such that for $0\le t\le T$,
\begin{equation}\label{eq-vt-H2-jg}
  \begin{aligned}[b]
    \| v_t(t) \|_{H^1(\mathbbm{R}_+^2)}^2&+\| {\rm{div}}{p}(t) \|_{H^1(\mathbbm{R}_+^2)}^2+ \|p(t) \|_{L^2(\mathbbm{R}_+^2)}^2\\
    &+\int_0^t (\| v_t(\tau) \|_{H^1(\mathbbm{R}_+^2)}^2+\| {\rm{div}}{p}_t(\tau) \|_{H^1(\mathbbm{R}_+^2)}^2 +\| {p}_t(\tau) \|_{L^2(\mathbbm{R}_+^2)}^2) \,\mathrm{d}\tau
    \le CM_0^2.
  \end{aligned}
\end{equation}
\end{lem}
 {\it\bfseries Proof.}
We can get from \eqref{eq-rdfc-vp}$_1$ and $\nabla\eqref{eq-rdfc-vp}_1$ that
\begin{equation}\label{eq-vt-H2-1}
  \begin{aligned}[b]
    \iint_{\mathbbm{R}_+^2} v_t^2 \,\mathrm{d}x\mathrm{d}y\le& C \iint_{\mathbbm{R}_+^2} \left(v_x^2+U_x^2v^2+v_y^2+({\rm{div}}{p})^2\right) \,\mathrm{d}x\mathrm{d}y ,\\
    \iint_{\mathbbm{R}_+^2} |\nabla{v}_t|^2 \,\mathrm{d}x\mathrm{d}y\le& C \iint_{\mathbbm{R}_+^2} (|\nabla{v}|^4+U_x^2|\nabla{v}|^2+U_x^4v^2+|\nabla^2{v}|^2+U_{xx}^2v^2+|\nabla{\rm{div}}{p}|^2) \,\mathrm{d}x\mathrm{d}y.
  \end{aligned}
\end{equation}
Integrating \eqref{eq-vt-H2-1} over $[0,t]$, by Lemma \ref{lem-nabv}-\ref{lem-nab2v}, we have
\begin{equation*}
	\begin{aligned}[b]
		\int_0^t \| v_t(\tau) \|_{L^2(\mathbbm{R}_+^2)}^2  \,\mathrm{d}\tau\le& C \int_0^t (\| \nabla{v}(\tau) \|_{L^2(\mathbbm{R}_+^2)}^2+\| \sqrt{U_x}v(\tau) \|_{L^2(\mathbbm{R}_+^2)}^2+\| {\rm{div}}{p}(\tau) \|_{L^2(\mathbbm{R}_+^2)}^2)  \,\mathrm{d}\tau\le CM_0^2,\\
		\int_0^t \| \nabla{v}_t(\tau) \|_{L^2(\mathbbm{R}_+^2)}^2  \,\mathrm{d}\tau \le& C \int_0^t (\| \nabla{v}(\tau) \|_{H^1(\mathbbm{R}_+^2)}^2+\| \sqrt{U_x}v(\tau) \|_{L^2(\mathbbm{R}_+^2)}^2+\| \nabla{\rm{div}}{p}(\tau) \|_{L^2(\mathbbm{R}_+^2)}^2+\| U_{xx}(\tau) \|_{L^\infty(\mathbbm{R}_+)}^2 )  \,\mathrm{d}\tau\\
		\le& CM_0^2.
	\end{aligned}
\end{equation*}
On the other hand, rewriting \eqref{eq-rdfc-vp}$_2$ in the form $\nabla{\rm{div}}{p}-p=\nabla{v}$ and squaring the resulting equation, we can get
\begin{equation*}
	|\nabla{\rm{div}}{p}|^2+2({\rm{div}}{p})^2+|p|^2-2 {\rm{div}}\{{\rm{div}}{p}p\}=|\nabla{v}|^2.
\end{equation*}
From ${\rm{div}}{p}(0,y,t)=0$, we further get
\begin{equation}\label{eq-divpH1}
	\begin{aligned}
		\| {\rm{div}}{p}(t) \|_{H^1(\mathbbm{R}_+^2)}^2+ \|p(t) \|_{L^2(\mathbbm{R}_+^2)}^2\le \| \nabla{v}(t) \|_{H^1(\mathbbm{R}_+^2)}^2\le CM_0^2,
	\end{aligned}
\end{equation}
it follows from \eqref{eq-vt-H2-1} that
\begin{equation}
	\begin{aligned}[b]
		\| v_t(t) \|_{H^1(\mathbbm{R}_+^2)}^2\le C \| v(t) \|_{H^2(\mathbbm{R}_+^2)}^2+ C \| {\rm{div}}{p}(t) \|_{H^1(\mathbbm{R}_+^2)}^2\le CM_0^2.
	\end{aligned}
\end{equation}
Similar to \eqref{eq-divpH1}, from ${\rm{div}}{p}_t(0,y,t)=0$, we also get
\begin{equation}\label{eq-divpt-H1}
	\| {\rm{div}}{p}_t(t) \|_{H^1(\mathbbm{R}_+^2)}^2+ \|p_t(t) \|_{L^2(\mathbbm{R}_+^2)}^2\le \| \nabla{v}_t(t) \|_{H^1(\mathbbm{R}_+^2)}^2.
\end{equation}
Integrating \eqref{eq-divpt-H1} over $[0,t]$, the desired estimate \eqref{eq-vt-H2-jg} can be obtained. Hence, we complete the proof of Lemma \ref{lem-vtH1-jg}.
$ \hfill\Box $

\begin{lem}Under the same assumptions of Proposition \ref{prop-priori-vp}, there exists a positive constant $C$ such that for $0\le t\le T$,
\begin{equation}\label{eq-nabdelpv-jg}
  \begin{aligned}[b]
    \| \nabla^3 v(t) \|_{L^2(\mathbbm{R}_+^2)}^2+\int_0^t \| \sqrt{U_x} \nabla\Delta v(t) \|_{L^2(\mathbbm{R}_+^2)}^2  \,\mathrm{d}\tau+\int_0^t (\| \nabla^3{\rm{div}}{p}(t) \|_{L^2(\mathbbm{R}_+^2)}^2 +\| \nabla^3 v(t) \|_{L^2(\mathbbm{R}_+^2)}^2) \,\mathrm{d}\tau\le C M_0^2.
  \end{aligned}
\end{equation}
\end{lem}
 {\it\bfseries Proof.}
We can get from $\nabla\Delta\eqref{eq-rdfc-vp}_1\cdot \nabla\Delta v-\nabla {\rm{div}}\eqref{eq-rdfc-vp}\cdot\nabla\Delta {\rm{div}}{p}$ that
\begin{equation}\label{eq-nabdelv-1}
  \begin{aligned}[b]
    \frac{1}{2} \frac{\mathrm{d}}{\mathrm{d}t}|\nabla\Delta{v}|^2&+ \nabla\Delta(f(v+U)-f(U))_x\cdot \nabla\Delta{v}+\nabla\Delta(g(v+U)_y)\cdot \nabla\Delta{v}\\
    &+|\nabla\Delta {\rm{div}}{p}|^2+(\Delta {\rm{div}}{p})^2- {\rm{div}}\{\Delta {\rm{div}}{p} \nabla{\rm{div}}{p}\}=0.
  \end{aligned}
\end{equation}
Integrating $\eqref{eq-nabdelv-1}$ over $\mathbbm{R}_+^2$, from Lemma \ref{lem-fv} $\mathrm{(C6)}$ and $\mathrm{(C12)}$, we have
\begin{equation}\label{eq-nabdelv-2}
  \begin{aligned}[b]
     &\frac{\mathrm{d}}{\mathrm{d}t}\| \nabla\Delta{v}(t) \|_{L^2(\mathbbm{R}_+^2)}^2+\frac{\alpha}{2}\|\sqrt{U_x} \nabla \Delta{v}(t)\|_{L^2(\mathbbm{R}_+^2)}^2+\| \nabla\Delta {\rm{div}}{p}(t) \|_{L^2(\mathbbm{R}_+^2)}^2 +\|\Delta {\rm{div}}{p}(t)\|_{L^2(\mathbbm{R}_+^2)}^2\\
     \le& C |\nabla\Delta{v}(0,y,t)|^2+C\| \nabla{v}(t) \|_{H^1(\mathbbm{R}_+^2)}^2 +C(\| \nabla{v} \|_{L^\infty(\mathbbm{R}_+^2)} +\delta')\| \nabla^3{v}(t) \|_{L^2(\mathbbm{R}_+^2)}^2+\| U_{x}v(t) \|_{L^2(\mathbbm{R}_+^2)}^2 \\
     &+\| U_{xx} \|_{L^\infty(\mathbbm{R}_+)}^4 \| v(t) \|_{L^2(\mathbbm{R}_+^2)}^2  +C \iint_{\mathbbm{R}_+^2} |\nabla^2{v}|^4 \,\mathrm{d}x \mathrm{d}y+C \int_{\mathbbm{R}_+} |U_{xxxx}|^2 \,\mathrm{d}x \int_{\mathbbm{R}} \sup\limits_{x\in\mathbbm{R}_+}|v|^2 \,\mathrm{d}y.
  \end{aligned}
\end{equation}
The terms on right-hand side of $\eqref{eq-nabdelv-2}$ can be estimated as follows. From Lemma \ref{lem-bjv-gj},
\begin{equation}\label{eq-nabdelv-3}
  \begin{aligned}[b]
    |\nabla\Delta{v}(0,y,t)|^2
    &=v_{xxx}^2(0,y,t)+v_{xxy}^2(0,y,t)\\
    &\le C ({\rm{div}}{p}_{x}(0,t))^2+C ({\rm{div}}{p}_{xt}(0,y,t))^2+C ({\rm{div}}{p}_{xy}(0,y,t))^2+C ({\rm{div}}{p}_{xx}(0,y,t))^2\\
    &\le \frac{1}{8}\| \nabla^3{\rm{div}}{p}(t) \|_{L^2(\mathbbm{R}_+^2)}^2+C \| \nabla {\rm{div}}{p}(t) \|_{H^1(\mathbbm{R}_+^2)}^2 +C ({\rm{div}}{p}_{xt}(0,y,t))^2\\
    &\le \frac{1}{4}(\| \nabla\Delta{\rm{div}}{p}(t) \|_{L^2(\mathbbm{R}_+^2)}^2+\| \nabla{\rm{div}}{p}_{yy}(t) \|_{L^2(\mathbbm{R}_+^2)}^2)+C \| \nabla {\rm{div}}{p}(t) \|_{H^1(\mathbbm{R}_+^2)}^2 +C ({\rm{div}}{p}_{xt}(0,y,t))^2.
  \end{aligned}
\end{equation}
To treat the last term of \eqref{eq-nabdelv-3}, we obtain from $\partial_t\eqref{eq-rd2-dengjia}_1$ that
\begin{equation*}
	({\rm{div}}{p}_{xt}(0,y,t))^2=|p_{1t}(0,y,t)|^2\le \| p_{1t}(t) \|_{L^\infty(\mathbbm{R}_+^2)}^2\le C(\| p_{1xt}(t) \|_{L^2(\mathbbm{R}_+^2)}^2+\| p_{1t}(t) \|_{L^2(\mathbbm{R}_+^2)}^2 ).
\end{equation*}
 Differentiating \eqref{eq-rd2-dengjia}$_1$ with respect to $t$ and then squaring this equation, we get $({\rm{div}}{p}_{xt}-p_{1t})^2=v_{xt}^2$, which implies that
\begin{equation}\label{eq-divpxt}
	({\rm{div}}{p}_{xt})^2+|p_{1t}|^2-2\{{\rm{div}}{p}_tp_{1t}\}_x+2p_{1xt}^2=v_{xt}^2-2p_{1xt}p_{2yt}.
\end{equation}
Integrating \eqref{eq-divpxt} over $\mathbbm{R}_+^2$, from ${\rm{div}}{p}_t(0,y,t)=0$, we have
\begin{equation*}
\begin{aligned}[b]
	 &\iint_{\mathbbm{R}_+^2} ({\rm{div}}{p}_{xt})^2 \,\mathrm{d}x\mathrm{d}y+2\iint_{\mathbbm{R}_+^2} (p_{1xt})^2 \,\mathrm{d}x\mathrm{d}y +\iint_{\mathbbm{R}_+^2} (p_{1t})^2 \,\mathrm{d}x\mathrm{d}y \\
	\le& \iint_{\mathbbm{R}_+^2} v_{xt}^2 \,\mathrm{d}x\mathrm{d}y+\iint_{\mathbbm{R}_+^2} (p_{1xt}+p_{2yt})^2 \,\mathrm{d}x\mathrm{d}y \\
	\le& \iint_{\mathbbm{R}_+^2} v_{xt}^2 \,\mathrm{d}x\mathrm{d}y+\iint_{\mathbbm{R}_+^2} ({\rm{div}}{p}_t)^2 \,\mathrm{d}x\mathrm{d}y.	
\end{aligned}
\end{equation*}
Therefore, it follows from \eqref{eq-nabdelv-3} that
\begin{equation}\label{eq-nabdelv-bjgj}
\begin{aligned}[b]
	|\nabla\Delta{v}(0,y,t)|^2
	\le& \frac{1}{4}(\| \nabla\Delta{\rm{div}}{p}(t) \|_{L^2(\mathbbm{R}_+^2)}^2+\| \nabla{\rm{div}}{p}_{yy}(t) \|_{L^2(\mathbbm{R}_+^2)}^2)\\
	&+C (\| \nabla {\rm{div}}{p}(t) \|_{H^1(\mathbbm{R}_+^2)}^2+ \| \nabla v_t(t) \|_{L^2(\mathbbm{R}_+^2)}^2 + \|  {\rm{div}}{p}_t(t) \|_{L^2(\mathbbm{R}_+^2)}^2).
\end{aligned}
\end{equation}
Similar to \eqref{eq-v-H1-3}, we can deduce that
\begin{equation}\label{eq-nabdelv-4}
    \iint_{\mathbbm{R}_+^2} |\nabla^2{v}|^4 \,\mathrm{d}x \mathrm{d}y
    \le CM_0^2 \| \nabla^3{v}(t) \|_{L^2(\mathbbm{R}_+^2)}^2.
\end{equation}
From \eqref{eq-U-sjgj},
\begin{equation}\label{eq-nabdelv-5}
  \begin{aligned}[b]
     \int_{\mathbbm{R}_+} |U_{xxxx}|^2 \,\mathrm{d}x \int_{\mathbbm{R}} \sup\limits_{x\in\mathbbm{R}_+}|v|^2 \,\mathrm{d}y
     \le& \int_{\mathbbm{R}_+} |U_{xxxx}|^2 \,\mathrm{d}x \iint_{\mathbbm{R}_+^2} (v^2+v_x^2) \,\mathrm{d}x \mathrm{d}y \\
    \le& C M_0^2 \| U_{xxxx}(t) \|_{L^2(\mathbbm{R}_+)}^2\\
    \le&  C M_0^2(1+t)^{-\frac{3}{2} }\mathrm{log}^{10}(2+t).
  \end{aligned}
\end{equation}
Substituting \eqref{eq-nabdelv-bjgj}-\eqref{eq-nabdelv-5} into \eqref{eq-nabdelv-2}, for sufficiently small $\delta'$ we have
\begin{equation}\label{eq-nabdelv-6}
  \begin{aligned}[b]
     &\frac{\mathrm{d}}{\mathrm{d}t}\| \nabla\Delta{v}(t) \|_{L^2(\mathbbm{R}_+^2)}^2+\|\sqrt{U_x} \nabla \Delta{v}(t)\|_{L^2(\mathbbm{R}_+^2)}^2+\| \nabla\Delta {\rm{div}}{p}(t) \|_{L^2(\mathbbm{R}_+^2)}^2 +\|\Delta {\rm{div}}{p}(t)\|_{L^2(\mathbbm{R}_+^2)}^2\\
     \le& C\| \nabla{v}(t) \|_{H^1(\mathbbm{R}_+^2)}^2 +C(M_0+\delta')\| \nabla^3{v}(t) \|_{L^2(\mathbbm{R}_+^2)}^2+C M_0^2(1+t)^{-\frac{3}{2} }\mathrm{log}^{10}(2+t)+C\| \nabla v_t(t) \|_{L^2(\mathbbm{R}_+^2)}^2\\
     &\!+\!\frac{1}{4}(\| \nabla\Delta{\rm{div}}{p}(t) \|_{L^2(\mathbbm{R}_+^2)}^2\!+\!\| \nabla{\rm{div}}{p}_{yy}(t) \|_{L^2(\mathbbm{R}_+^2)}^2)\!+\!C \| \nabla {\rm{div}}{p}(t) \|_{H^1(\mathbbm{R}_+^2)}^2 \!+\!C \|  {\rm{div}}{p}_t(t) \|_{L^2(\mathbbm{R}_+^2)}^2\!+\!\| \sqrt{U_x}v(t) \|_{L^2(\mathbbm{R}_+^2)}^2 .
  \end{aligned}
\end{equation}
Integrating \eqref{eq-nabdelv-6} over $[0,t]$, substituting \eqref{eq-vt-H2-jg} and \eqref{eq-vyy-H1-jg} into the resulting inequality, we can get
\begin{equation}\label{eq-nabdivp-gj}
  \begin{aligned}[b]
    &\| \nabla\Delta{v}(t) \|_{L^2(\mathbbm{R}_+^2)}^2+\int_0^t (\|\sqrt{U_x} \nabla \Delta{v}(\tau)\|_{L^2(\mathbbm{R}_+^2)}^2+\| \nabla\Delta {\rm{div}}{p}(\tau) \|_{L^2(\mathbbm{R}_+^2)}^2 +\|\Delta {\rm{div}}{p}(\tau)\|_{L^2(\mathbbm{R}_+^2)}^2) \,\mathrm{d}\tau \\
    \le& CM_0^2+C(M_0+\delta')\int_0^t \| \nabla^3{v}(\tau) \|_{L^2(\mathbbm{R}_+^2)}^2 \,\mathrm{d}\tau.
  \end{aligned}
\end{equation}
From \eqref{eq-rdfc-vp}$_2$, combining \eqref{eq-vyy-H1-jg}, the last term can be estimated as
\begin{equation*}
  \begin{aligned}[b]
    \int_0^t \| \nabla^3{v}(\tau) \|_{L^2(\mathbbm{R}_+^2)}^2 \,\mathrm{d}\tau
    \le& C \int_0^t \| \nabla \Delta{v}(\tau) \|_{L^2(\mathbbm{R}_+^2)}^2 \,\mathrm{d}\tau+C \int_0^t \| \nabla {v}_{yy}(\tau) \|_{L^2(\mathbbm{R}_+^2)}^2 \,\mathrm{d}\tau\\
    \le& C \int_0^t (\| \nabla \Delta {\rm{div}}{p}(\tau) \|_{L^2(\mathbbm{R}_+^2)}^2+\| \nabla{\rm{div}}{p}(\tau) \|_{L^2(\mathbbm{R}_+^2)}^2) \,\mathrm{d}\tau\\
    &+C \int_0^t (\| \nabla {\rm{div}}{p}_{yy}(\tau) \|_{L^2(\mathbbm{R}_+^2)}^2+\| {p}_{yy}(\tau) \|_{L^2(\mathbbm{R}_+^2)}^2) \,\mathrm{d}\tau\\
    \le& CM_0^2+C(M_0+\delta')\int_0^t \| \nabla^3{v}(\tau) \|_{L^2(\mathbbm{R}_+^2)}^2 \,\mathrm{d}\tau,
  \end{aligned}
\end{equation*}
which implies that
\begin{equation}\label{eq-nav3v-cjf}
  \int_0^t \| \nabla^3{v}(\tau) \|_{L^2(\mathbbm{R}_+^2)}^2 \,\mathrm{d}\tau\le CM_0^2,
\end{equation}
for some small $M_0$ and $\delta'$. Combining \eqref{eq-nabdivp-gj}, \eqref{eq-vyy-H1-jg} and \eqref{eq-nav3v-cjf}, the desired estimate \eqref{eq-nabdelpv-jg} can be directly obtained.
 $ \hfill\Box$

\begin{lem}Under the same assumptions of Proposition \ref{prop-priori-vp}, there exists a positive constant $C$ such that for $0\le t\le T$,
	\begin{equation}\label{eq-p-djf-jg}
		\begin{aligned}[b]
			\| \nabla^2{\rm{div}}{p}(t) \|_{H^1(\mathbbm{R}_+^2)}^2+\| \nabla p(t) \|_{H^2(\mathbbm{R}_+^2)}^2+\int_0^t (\| \nabla^2{\rm{div}}{p}(\tau) \|_{H^1(\mathbbm{R}_+^2)}^2+\| \nabla p(\tau) \|_{H^2(\mathbbm{R}_+^2)}^2) \,\mathrm{d}\tau \le CM_0^2.
		\end{aligned}
	\end{equation}
\end{lem}
 {\it\bfseries Proof.}
Making use of \eqref{eq-rdfc-vp}$_2$, we can deduce that
\begin{equation*}
	\begin{aligned}
		&\| \nabla^2{\rm{div}}{p}(t) \|_{L^2(\mathbbm{R}_+^2)}^2
		= \| \Delta{\rm{div}}{p}(t) \|_{L^2(\mathbbm{R}_+^2)}^2
		\le C(\| {\rm{div}}{p}(t) \|_{L^2(\mathbbm{R}_+^2)}^2+\| \Delta v(t) \|_{L^2(\mathbbm{R}_+^2)}^2),\\
		&\| \nabla\Delta{\rm{div}}{p}(t) \|_{L^2(\mathbbm{R}_+^2)}^2
		\le C(\| \nabla{\rm{div}}{p}(t) \|_{L^2(\mathbbm{R}_+^2)}^2+\| \nabla\Delta v(t) \|_{L^2(\mathbbm{R}_+^2)}^2),\\
        &\| \nabla^3{\rm{div}}{p}(t) \|_{L^2(\mathbbm{R}_+^2)}^2
        \le \| \nabla{\rm{div}}{p}_{yy}(t) \|_{L^2(\mathbbm{R}_+^2)}^2+\| \nabla\Delta{\rm{div}}{p}(t) \|_{L^2(\mathbbm{R}_+^2)}^2,
	\end{aligned}
\end{equation*}
and
\begin{equation*}
	\| \nabla{\rm{div}}{p}_{yy}(t) \|_{L^2(\mathbbm{R}_+^2)}^2+2\| {\rm{div}}{p}_{yy}(t) \|_{L^2(\mathbbm{R}_+^2)}^2+\| {p}_{yy}(t) \|_{L^2(\mathbbm{R}_+^2)}^2\le \| \nabla v_{yy}(t) \|_{L^2(\mathbbm{R}_+^2)}^2.
\end{equation*}
Therefore, it holds that
\begin{equation}\label{eq-nab2divp-djf}
	\| \nabla^2{\rm{div}}{p}(t) \|_{H^1(\mathbbm{R}_+^2)}^2+\int_0^t \| \nabla^2{\rm{div}}{p}(\tau) \|_{H^1(\mathbbm{R}_+^2)}^2 \,\mathrm{d}\tau \le CM_0^2.
\end{equation}
By utilizing \eqref{eq-rdfc-vp}$_2$ again, we have
\begin{equation*}
	\begin{aligned}[b]
		&\| \nabla{p}(t) \|_{L^2(\mathbbm{R}_+^2)}^2\le C(\| \nabla^2{\rm{div}}{p}(t) \|_{L^2(\mathbbm{R}_+^2)}^2+\| \nabla^2v(t) \|_{L^2(\mathbbm{R}_+^2)}^2),\\
		&\| \nabla^2{p}(t) \|_{L^2(\mathbbm{R}_+^2)}^2\le C(\| \nabla^3{\rm{div}}{p}(t) \|_{L^2(\mathbbm{R}_+^2)}^2+\| \nabla^3v(t) \|_{L^2(\mathbbm{R}_+^2)}^2).
	\end{aligned}
\end{equation*}
In order to show the $L^2$-estimate of $\nabla^3 p$, we know the definition of $\nabla^3 p$ is that
 \begin{equation*}
 	\nabla^3p=
 	\left(
  \begin{array}{cccc}
    p_{1xxx} & p_{1xxy} & p_{1xyy} & p_{1yyy}\\
    p_{2xxx} & p_{2xxy} & p_{2xyy} & p_{2yyy}\\
  \end{array}
\right).
 \end{equation*}
 From $\eqref{eq-p1p2}$, we note that
	\begin{align}
	    &{\rm{div}}{p}_{xx}-p_{1xxx}=p_{2xxy}=p_{1xyy}={\rm{div}}{p}_{yy}-p_{2yyy},\label{eq-p1xxx}\\
		&p_{2xxx}=p_{1xxy}={\rm{div}}{p}_{xy}-p_{2xyy}={\rm{div}}{p}_{xy}-p_{1yyy}.\label{eq-p2xxx}
	\end{align}
Therefore, we just need to show the $L^2$-estimate of $p_{2xyy}$ and $p_{2yyy}$. Once we get these two estimates, combined with the obtained estimates of $\nabla^2 {\rm{div}}{p}$, we can deduce the estimate of $\nabla^3 p$ immediately by using \eqref{eq-p1xxx}-\eqref{eq-p2xxx}. Now, we estimate $p_{2xyy}$ and $p_{2yyy}$.

We can get from $\partial_x\partial_y$\eqref{eq-rd2-dengjia}$_2$ that
 \begin{equation}\label{eq-p-H3-1}
 	({\rm{div}}{p}_{xyy})^2+p_{2xy}^2+2p_{2xyy}^2-2\left\{ {\rm{div}}{p}_{xy}p_{2xy} \right\}_y=v_{xyy}^2-2p_{1xxy}p_{2xyy}.
 \end{equation}
Integrating $\eqref{eq-p-H3-1}$ over $\mathbbm{R}_+^2$, we have
\begin{equation*}
	\begin{aligned}
		\iint_{\mathbbm{R}_+^2} ({\rm{div}}{p}_{xyy})^2 \,\mathrm{d}x \mathrm{d}y+\iint_{\mathbbm{R}_+^2} p_{2xy}^2 \,\mathrm{d}x \mathrm{d}y+2\iint_{\mathbbm{R}_+^2} p_{2xyy}^2 \,\mathrm{d}x \mathrm{d}y
		\le&  \iint_{\mathbbm{R}_+^2} v_{xyy}^2 \,\mathrm{d}x \mathrm{d}y+\iint_{\mathbbm{R}_+^2} (p_{1xxy}+p_{2xyy})^2 \,\mathrm{d}x \mathrm{d}y\\
		\le& \iint_{\mathbbm{R}_+^2} v_{xyy}^2 \,\mathrm{d}x \mathrm{d}y+\iint_{\mathbbm{R}_+^2} ({\rm{div}}{p}_{xy})^2 \,\mathrm{d}x \mathrm{d}y.
	\end{aligned}
\end{equation*}
Similarly, we can get from $\partial_y^2$\eqref{eq-rd2-dengjia}$_2$ that
 \begin{equation*}
 	({\rm{div}}{p}_{yyy})^2+p_{2yy}^2+2p_{2yyy}^2-2\left\{ {\rm{div}}{p}_{yy}p_{2yy} \right\}_y=v_{yyy}^2-2p_{1xyy}p_{2yyy}.
 \end{equation*}
 It follows that
\begin{equation*}
	\begin{aligned}
		\iint_{\mathbbm{R}_+^2} ({\rm{div}}{p}_{yyy})^2 \,\mathrm{d}x \mathrm{d}y+\iint_{\mathbbm{R}_+^2} p_{2yy}^2 \,\mathrm{d}x \mathrm{d}y+2\iint_{\mathbbm{R}_+^2} p_{2yyy}^2 \,\mathrm{d}x \mathrm{d}y
		\le& \iint_{\mathbbm{R}_+^2} v_{yyy}^2 \,\mathrm{d}x \mathrm{d}y+\iint_{\mathbbm{R}_+^2} ({\rm{div}}{p}_{yy})^2 \,\mathrm{d}x \mathrm{d}y.
	\end{aligned}
\end{equation*}
Combining the \eqref{eq-nabdelpv-jg} and \eqref{eq-nab2divp-djf}, we can obtain the desired estimate \eqref{eq-p-djf-jg}.
 $ \hfill\Box $
\subsubsection{Large-time behavior}\label{sec-3-3}
The combination of the existence and uniqueness of the local solution and the \emph{a priori} estimates can extend the local solution for problem \eqref{eq-rdfc-vp}-\eqref{eq-boundary-vp} globally, that is
\begin{equation*}
\begin{cases}
	v\in C^0([0,\infty);H^3(\mathbbm{R}_+^2))\cap C^1([0,\infty);H^2(\mathbbm{R}_+^2)), \quad \nabla v \in L^2(0,\infty;H^2(\mathbbm{R}_+^2)),\\[1mm]
	p\in C^0([0,\infty);H^3(\mathbbm{R}_+^2))\cap L^2(0,\infty;H^3(\mathbbm{R}_+^2)), \quad  {\rm{div}}{p}\in C^0([0,\infty);H^3(\mathbbm{R}_+^2))\cap L^2(0,\infty;H^3(\mathbbm{R}_+^2)),\\[1mm]
	 v_t \in L^2(0,\infty;H^1(\mathbbm{R}_+^2)), \quad  p_t\in L^2(0,\infty;L^2(\mathbbm{R}_+^2)),  \quad  {\rm{div}}{p}_t\in L^2(0,\infty;H^1(\mathbbm{R}_+^2)).
\end{cases}
\end{equation*}
It follows that for $t\ge 0$,
\begin{equation}\label{eq-bigtime}
\begin{aligned}[b]
	&\| v(t) \|_{H^3(\mathbbm{R}_+^2)}^2+\| v_t(t) \|_{H^1(\mathbbm{R}_+^2)}^2+\| p(t) \|_{H^3(\mathbbm{R}_+^2)}^2+\| {\rm{div}}{p}(t) \|_{H^3(\mathbbm{R}_+^2)}^2+\int_0^\infty (\| \nabla{v}(t) \|_{H^2(\mathbbm{R}_+^2)}^2+\| v_t(t) \|_{H^1(\mathbbm{R}_+^2)}^2) \,\mathrm{d}\tau\\
	& \ \ +\int_0^\infty (\| p(t) \|_{H^3(\mathbbm{R}_+^2)}^2+\| {\rm{div}}{p}(t) \|_{H^3(\mathbbm{R}_+^2)}^2+\| p_t(t) \|_{L^2(\mathbbm{R}_+^2)}^2+\| {\rm{div}}{p}_t(t) \|_{H^1(\mathbbm{R}_+^2)}^2 )  \,\mathrm{d}\tau <\infty,\\
\end{aligned}
\end{equation}
In order to show the Large-time behavior \eqref{eq-largebehavior-v} in Theorem \ref{thm-vp-main}. From \eqref{eq-bigtime}, we just need to proof that
\begin{equation}\label{eq-bigtime1}
	\begin{aligned}[b]
		&\int_0^\infty \left|\frac{\mathrm{d}}{\mathrm{d}t}( \|\nabla{v}(t)\|_{L^2(\mathbbm{R}_+^2)}^2\!+\!\|\nabla^2{v}(t)\|_{L^2(\mathbbm{R}_+^2)}^2) \right| \,\mathrm{d}t<\infty,\\
		&\int_0^\infty \left|\frac{\mathrm{d}}{\mathrm{d}t}(\| \nabla{\rm{div}}{p} (t)\|_{L^2(\mathbbm{R}_+^2)}^2+\| {\rm{div}}{p}(t) \|_{L^2(\mathbbm{R}_+^2)}^2  +\|p(t)\|_{L^2(\mathbbm{R}_+^2)}^2 )\right| \,\mathrm{d}t<\infty.
	\end{aligned}
\end{equation}
Once we get the above two inequalities, combining \eqref{eq-bigtime}, we have
$$\| \nabla{v}(t) \|_{L^2(\mathbbm{R}_+^2)}, ~ ~ \| \nabla^2{v}(t) \|_{L^2(\mathbbm{R}_+^2)}, ~~ \| \nabla{\rm{div}}{p} (t)\|_{L^2(\mathbbm{R}_+^2)}, ~~ \| {\rm{div}}{p}(t) \|_{L^2(\mathbbm{R}_+^2)}, ~~ \|p(t)\|_{L^2(\mathbbm{R}_+^2)}\rightarrow 0, \quad \text{as } t\rightarrow \infty.$$
 By utilizing Gagliardo-Nirenberg inequality
$$\| f \|_{L^\infty(\mathbbm{R}_+^2)}\le C \| f \|_{L^2(\mathbbm{R}_+^2)}^\frac{1}{2} \| \nabla^2 f\|_{L^2(\mathbbm{R}_+^2)}^\frac{1}{2}, ~~~~ \| \nabla f \|_{L^\infty(\mathbbm{R}_+^2)}\le C \| f \|_{L^2(\mathbbm{R}_+^2)}^\frac{1}{3} \| \nabla^3 f\|_{L^2(\mathbbm{R}_+^2)}^\frac{2}{3},$$
we can deduce that as $t\rightarrow \infty$,
\begin{equation*}
	\begin{aligned}[b]
	&\|\nabla^k v(t) \|_{L^\infty(\mathbbm{R}_+^2)}\le C \| \nabla^k v(t) \|_{L^2(\mathbbm{R}_+^2)}^{1/2}\| \nabla^{k+2} v(t) \|_{L^2(\mathbbm{R}_+^2)}^{1/2} \rightarrow 0, ~~~~ k=0,1 \\
    &\| p(t) \|_{L^\infty(\mathbbm{R}_+^2)}\le C \| p(t) \|_{L^2(\mathbbm{R}_+^2)}^{1/2}\| \nabla^2{p}(t) \|_{L^2(\mathbbm{R}_+^2)}^{1/2} \rightarrow 0, \\
    &\| \nabla p(t) \|_{L^\infty(\mathbbm{R}_+^2)}\le C \| p(t) \|_{L^2}^{1/3}\| \nabla^3{p}(t) \|_{L^2(\mathbbm{R}_+^2)}^{2/3} \rightarrow 0,\\
    &\| \nabla {\rm{div}}{p}(t) \|_{L^\infty(\mathbbm{R}_+^2)} \le C \| \nabla {\rm{div}}{p}(t) \|_{L^2(\mathbbm{R}_+^2)}^{1/2}\| \nabla^{3} {\rm{div}}{p}(t) \|_{L^2(\mathbbm{R}_+^2)}^{1/2} \rightarrow 0.
    \end{aligned}
\end{equation*}
In order to prove \eqref{eq-bigtime1}, the key is to estimate $\nabla^2{v}_t$. We can get from $\Delta$\eqref{eq-rdfc-vp}$_1$ that
\begin{equation}
	\begin{aligned}[b]
		\| \Delta{v}_t(t) \|_{L^2(\mathbbm{R}_+^2)}^2\le C(\| \nabla{v}(t) \|_{H^2(\mathbbm{R}_+^2)}^2+\| \sqrt{U_x}v(t) \|_{L^2(\mathbbm{R}_+^2)}^2+\| U_{xxx} \|_{L^\infty(\mathbbm{R}_+)}^2 \| v(t) \|_{L^2(\mathbbm{R}_+^2)}^2+\| \Delta {\rm{div}}{p}(t) \|_{L^2(\mathbbm{R}_+^2)}^2 ).
	\end{aligned}
\end{equation}
Similar to $\eqref{eq-delv-zh}$, we can deduce that
$$ \| \Delta v_t(t) \|_{L^2(\mathbbm{R}_+^2)}^2=\| \nabla^2v_t(t) \|_{L^2(\mathbbm{R}_+^2)}^2. $$
It follows that
\begin{equation}
	\begin{aligned}[b]
		\int_0^{\infty} \| \nabla^2 v_t(t) \|_{L^2(\mathbbm{R}_+^2)}^2 \,\mathrm{d}t <\infty.
	\end{aligned}
\end{equation}
Then, we can obtain
\begin{equation*}
\begin{aligned}
	&\int_0^\infty \left|\frac{\mathrm{d}}{\mathrm{d}t}( \|\nabla{v}(t)\|_{L^2(\mathbbm{R}_+^2)}^2\!+\!\|\nabla^2{v}(t)\|_{L^2(\mathbbm{R}_+^2)}^2) \right| \,\mathrm{d}t\\
	\le& C \int_0^\infty (\|\nabla{v}(t)\|_{L^2(\mathbbm{R}_+^2)}^2+\|\nabla{v}_t(t)\|_{L^2(\mathbbm{R}_+^2)}^2+\|\nabla^2{v}(t)\|_{L^2(\mathbbm{R}_+^2)}^2+\|\nabla^2{v}_t(t)\|_{L^2(\mathbbm{R}_+^2)}^2)  \, \mathrm{d}t\\
	 <& \infty,	
\end{aligned}
\end{equation*}
and
\begin{equation*}
	\begin{aligned}[b]
		&\int_0^\infty \left|\frac{\mathrm{d}}{\mathrm{d}t}(\| \nabla{\rm{div}}{p} (t)\|_{L^2(\mathbbm{R}_+^2)}^2+\| {\rm{div}}{p}(t) \|_{L^2(\mathbbm{R}_+^2)}^2  +\|p(t)\|_{L^2(\mathbbm{R}_+^2)}^2 )\right| \,\mathrm{d}t\\
		\le&C \int_0^\infty (\| \nabla{\rm{div}}{p}(t) \|_{L^2(\mathbbm{R}_+^2)}^2+\| {\rm{div}}{p}(t) \|_{L^2(\mathbbm{R}_+^2)}^2+\|p(t) \|_{L^2(\mathbbm{R}_+^2)}^2)  \,\mathrm{d}t \\
		&+C \int_0^\infty (\| \nabla{\rm{div}}{p}_t(t) \|_{L^2(\mathbbm{R}_+^2)}^2+\| {\rm{div}}p_t (t) \|_{L^2(\mathbbm{R}_+^2)}^2+\|p_t(t) \|_{L^2(\mathbbm{R}_+^2)}^2)  \,\mathrm{d}t\\
		<& \infty.
	\end{aligned}
\end{equation*}
Therefore, the proof of Theorem \ref{thm-vp-main} is completed. Combined with the proved Theorem \ref{thm-V-main}, we prove the main Theorem \ref{thm-main}.

 \vspace{6mm}

\noindent {\bf Acknowledgements:}
The research was supported by the National Natural Science Foundation of China $\#$12171160, 11771150, $\#$11831003 and Guangdong Basic and Applied Basic Research Foundation $\#$2020B1515310015.

\end{document}